\documentclass[preprint,11pt]{elsarticle}

\usepackage[utf8]{inputenc}
\usepackage{fullpage,times}
\usepackage{amsmath,amssymb, mathtools}
\usepackage{graphicx}
\usepackage{algorithm}
\usepackage{algorithmicx}
\usepackage{subfig}
\usepackage[font=small,skip=-1pt]{caption}
\usepackage{url}
\usepackage[numbers]{natbib} 
\usepackage{placeins}
\usepackage{siunitx}
\usepackage{url}
\usepackage{xspace}
\usepackage[table,xcdraw]{xcolor}

\usepackage[titletoc,title]{appendix}
\usepackage{xpatch}
\let\appendixpagenameorig\appendixpagename
\renewcommand{\appendixpagename}{\sffamily\appendixpagenameorig}

\usepackage{comment}

\usepackage{hyperref}
\hypersetup{
colorlinks={true},
linkcolor={Blue3},
urlcolor={Blue3},
citecolor={Blue3},
pdfstartview={FitH},
pdfauthor={Author},
pdftitle={Title},
pdfsubject={Subject}
plainpages=false
}

\usepackage{enumitem}

\usepackage{booktabs} 
\usepackage{multirow}
\usepackage{multicol}
\usepackage{makecell}

\usepackage{setspace}

\definecolor{wisconsin-red}{rgb}{0.6,0,0}
\definecolor{darkgreen}{rgb}{0.2,0.6,0.2}
\definecolor{maroon}{rgb}{0.5, 0.0, 0.0}
\definecolor{violet}{rgb}{0.75, 0.0, 1.0}
\definecolor{lightgray}{gray}{0.9}
\definecolor{navyblue}{rgb}{0.0, 0.0, 0.5}
\definecolor{darkmidnightblue}{rgb}{0.0, 0.2, 0.4}
\definecolor{midnightblue}{rgb}{0.0,0.4,0.85}
\definecolor{Gray}{gray}{0.75}
\definecolor{darkgreen}{rgb}{0,0.5,0}
\definecolor{apricot}{rgb}{0.98, 0.81, 0.69}

\newcommand{\setItemSep}{\setlength\itemsep}

\newcolumntype{C}[1]{>{\centering\arraybackslash}p{#1}}
\newcolumntype{P}[1]{>{\raggedright\arraybackslash}p{#1}}
\newcolumntype{L}[1]{>{\raggedleft\arraybackslash}p{#1}}

\usepackage{upgreek}
\newcommand{\nl}{\ensuremath{TN_{l}}}

\newcommand{\nlq}{\ensuremath{TN_{l'}}}
\newcommand{\ns}{\ensuremath{N}}

\newcommand{\tp}{\ensuremath{TP}}
\newcommand{\tph}{\ensuremath{TP_{Hf}}}
\newcommand{\tpl}{\ensuremath{TP_{Lf}}}
\newcommand{\tpqq}{\ensuremath{CL^{n}}}
\newcommand{\tpqqfh}{\ensuremath{CL^{n}_{Hf}}}
\newcommand{\tpqqfl}{\ensuremath{CL^{n}_{Lf}}}

\newcommand{\hlq}{\ensuremath{h_{l'}}}
\newcommand{\ctw}{\ensuremath{c^{{tw}}}}

\newcommand{\hl}{\ensuremath{h_{l}}}

\newcommand{\pl}{\ensuremath{P_{l}}}
\newcommand{\ql}{\ensuremath{Q_{l'}}}

\newcommand{\plpn}{\ensuremath{Q^-_{l'}}}

\newcommand{\hlmax}{\ensuremath{h^{max}_l}}
\newcommand{\hlmin}{\ensuremath{h^{min}_l}}

\newcommand{\muupkk}{\ensuremath{{\mu^{{(k+1)ns}}_{lp}}}}
\newcommand{\muup}{\ensuremath{{\mu^{{kns}}_{lp}}}}
\newcommand{\muupz}{\ensuremath{{\mu^{{(0)ns}}_{lp}}}}

\newcommand{\GBDEw}{\ensuremath{\text{GBDEw}^{sn}_{l'q}}}
\newcommand{\GBDEwp}{\ensuremath{\text{GBDEw}^{sn}_{lp}}}

\newcommand{\TotalBDEw}{\ensuremath{\text{TotalBDEw}^{sn}_{l'q}}}
\newcommand{\TotalBDEwp}{\ensuremath{\text{TotalBDEw}^{sn}_{lp}}}

\newcommand{\plpqn}{\ensuremath{Q^-_{l'}}}
\newcommand{\plpq}{\ensuremath{Q_{l'}}}
\newcommand{\plpqs}{\ensuremath{Q^*_{l'}}}
\newcommand{\tb}{\ensuremath{{\text{dtb}^{n}_{l'q}}}}

\newcommand{\tbb}{\ensuremath{{\text{ptb}^{n}_{l'q}}}}

\newcommand{\twait}{\ensuremath{{\text{Twait}^{n}_{lp{l'}q}}}}

\newcommand{\ndwait}{\ensuremath{{\text{NTwait}^{n}_{lp{l'}q}}}}

\newcommand{\ndwaitsss}{\ensuremath{{\text{NTwait}^{sn}_{lp{l'}q^{*}}}}}

\newcommand{\zy}{\ensuremath{{\text{ZY}^{n}_{lp{l'}q}}}}

\newcommand{\cvt}{\ensuremath{c^{vt}}}
\newcommand{\cvh}{\ensuremath{c^{dt}}}

\newcommand{\qS}{\ensuremath{q^*_{l'}}}

\newcommand{\prob}{\ensuremath{p^{s}}}
\newcommand{\ttds}{\ensuremath{ttd_{lp}^{{s(n_{prev}-n)}}}}
\newcommand{\aws}{\ensuremath{{awt}^{sn}_{ll'}}}
\newcommand{\tds}{\ensuremath{{td}^{sn}_{lpl'}}}
\newcommand{\ssps}{\ensuremath{{sp^{s(n_{prev}{-n})}_{l'q}}}}
\newcommand{\TldLs}{\ensuremath{{LTld}^{sn}_{l'q}}}
\newcommand{\alds}{\ensuremath{{ad}^{sn}_{l'q}}}

\newcommand{\lams}{\ensuremath{\lambda^{sn}_{l'}}}

\newcommand{\Pdepps}{\ensuremath{{\text{Pdep}^{ksn}_{lp}}}}
\newcommand{\Pdeppsk}{\ensuremath{{\text{Pdep}^{snk}_{lp}}}}
\newcommand{\Pdeppsnk}{\ensuremath{{\text{Pdep}^{sn}_{lp}}}}

\newcommand{\Pdeppskk}{\ensuremath{{\text{Pdep}^{sn(k+1)}_{lp}}}}
\newcommand{\Pdeppskz}{\ensuremath{{\text{Pdep}^{sn(0)}_{lp}}}}
\newcommand{\Pdeppsa}{\ensuremath{{\overline{\text{Pdep}}^{n{k_{prev}}}_{lp}}}}

\newcommand{\Pdeppsakk}{\ensuremath{{\overline{\text{Pdep}}^{n(k+1)}_{lp}}}}
\newcommand{\Pdeppsakz}{\ensuremath{{\overline{\text{Pdep}}^{n(0)}_{lp}}}}
\newcommand{\Pdeppn}{\ensuremath{{\text{Pdep}^{n}_{l(p-1)}}}}
\newcommand{\Pdep}{\ensuremath{{\text{Pdep}^{n}_{lp}}}}
\newcommand{\Pdepq}{\ensuremath{{\text{Pdep}^{n}_{l'q}}}}
\newcommand{\Adep}{\ensuremath{{\text{Adep}^{sn}_{l'q}}}}

\newcommand{\Adepq}{\ensuremath{{\text{Adep}^{sn}_{lp}}}}
\newcommand{\Adepqp}{\ensuremath{{\text{Adep}^{sn}_{l(p-1)}}}}
\newcommand{\Adepp}{\ensuremath{{\text{Adep}^{s(n_{prev})}_{lp}}}}

\newcommand{\Aarr}{\ensuremath{{\text{Aarr}^{sn}_{lp}}}}
\newcommand{\Aarrq}{\ensuremath{{\text{Aarr}^{sn}_{l'q}}}}

\newcommand{\servos}{{\ensuremath{\text{serv1}^{sn}_{l'q}}}}
\newcommand{\serves}{{\ensuremath{\text{serv5}^{sn}_{l'q}}}}
\newcommand{\servts}{{\ensuremath{\text{serv2}^{sn}_{l'q}}}}
\newcommand{\servtsq}{{\ensuremath{\text{serv2q}^{sn}_{l'q}}}}
\newcommand{\servthsq}{{\ensuremath{\text{serv3q}^{sn}_{l'q}}}}
\newcommand{\servthsg}{{\ensuremath{\text{serv3g}^{sn}_{l'q}}}}
\newcommand{\servtsg}{{\ensuremath{\text{serv2g}^{sn}_{l'q}}}}
\newcommand{\servths}{{\ensuremath{\text{serv3}^{sn}_{l'q}}}}
\newcommand{\servfs}{{\ensuremath{\text{serv4}^{sn}_{l'q}}}}
\newcommand{\tbs}{\ensuremath{{\text{dtb}^{sn}_{l'q}}}}
\newcommand{\ibs}{\ensuremath{{\text{ib}^{sn}_{l'q}}}}

\newcommand{\api}{\ensuremath{{\text{ap}^{sn}_{l'q}}}}
\newcommand{\ipp}{\ensuremath{{\text{pi}^{sn}_{l'q}}}}
\newcommand{\pipa}{\ensuremath{{\text{pipa}^{sn}_{l'q}}}}
\newcommand{\ippp}{\ensuremath{{\text{pi}^{sn}_{lp}}}}
\newcommand{\Esthii}{\ensuremath{{\text{Esthii}^{sn}_{l'q}}}}
\newcommand{\TYE}{\ensuremath{{\text{TYE}^{sn}_{l'q}}}}
\newcommand{\STTYE}{\ensuremath{{\text{STTYE}^{sn}_{l'q}}}}
\newcommand{\ttbos}{\ensuremath{{\text{tb}^{sn}_{l'q}}}}
\newcommand{\pttbos}{\ensuremath{{\text{ptb}^{sn}_{l'q}}}}
\newcommand{\tbbs}{\ensuremath{{\text{ptb}^{sn}_{l'q}}}}
\newcommand{\tbbso}{\ensuremath{{\text{ptb1}^{sn}_{l'q}}}}
\newcommand{\tbbst}{\ensuremath{{\text{ptb2}^{sn}_{l'q}}}}
\newcommand{\tbbsth}{\ensuremath{{\text{ptb3}^{sn}_{l'q}}}}
\newcommand{\TIs}{\ensuremath{{\text{TI}^{sn}_{l'q}}}}
\newcommand{\TIse}{\ensuremath{{\text{TIE}^{sn}_{l'q}}}}

\newcommand{\nudemand}{\nu^{n}_{l'q}}
\newcommand{\yys}{\ensuremath{{\text{Y}^{sn}_{lp{l'}q}}}}
\newcommand{\yyns}{\ensuremath{{\text{Y}^{sn}_{lp{l'}(q-1)}}}}
\newcommand{\zys}{\ensuremath{{\text{ZY}^{sn}_{lp{l'}q}}}}
\newcommand{\sys}{\ensuremath{{\text{SY}^{sn}_{lp{l'}}}}}

\newcommand{\yyyys}{\ensuremath{{\text{YY}^{sn}_{lp{l'}q}}}}
\newcommand{\tys}{\ensuremath{{\text{TY}^{sn}_{lp{l'}q}}}}
\newcommand{\ttys}{\ensuremath{{\text{TY}^{sn}_{lp{l'}{q^*}}}}}
\newcommand{\twaits}{\ensuremath{{\text{Twait}^{sn}_{lp{l'}q}}}}
\newcommand{\iwaits}{\ensuremath{{\text{iwait}^{sn}_{lp{l'}q}}}}
\newcommand{\ndwaitss}{\ensuremath{{\text{NTwait}^{sn}_{lp{l'}q}}}}

\newcommand{\ivdsp}{\ensuremath{{\text{ivd}^{sn}_{lp}}}}

\newcommand{\bdos}{\ensuremath{\text{GBD1}^{sn}_{{l'q}}}}

\newcommand{\bdts}{\ensuremath{\text{GBD2}^{sn}_{{l'q}}}}
\newcommand{\bdths}{\ensuremath{\text{GBD3}^{sn}_{{l'q}}}}
\newcommand{\bdthsq}{\ensuremath{\text{GBD3q}^{sn}_{{l'q}}}}
\newcommand{\bdthsg}{\ensuremath{\text{GBD3g}^{sn}_{{l'q}}}}
\newcommand{\bdq}{\ensuremath{\text{GBD2q}^{sn}_{{l'q}}}}
\newcommand{\bdg}{\ensuremath{\text{GBD2g}^{sn}_{{l'q}}}}
\newcommand{\bdfs}{\ensuremath{\text{GBD4}^{sn}_{{l'q}}}}
\newcommand{\bdes}{\ensuremath{\text{GBD5}^{sn}_{{l'q}}}}

\newcommand{\Tbdns}{\ensuremath{{\text{Tbd}^{sn}_{l'q}}}}
\newcommand{\Tbdnsp}{\ensuremath{{\text{Tbd}^{sn}_{lp}}}}

\newcommand{\STo}{\ensuremath{{\text{ST1}^{sn}_{l'q}}}}
\newcommand{\STt}{\ensuremath{{\text{ST2}^{sn}_{l'q}}}}
\newcommand{\STth}{\ensuremath{{\text{ST3}^{sn}_{l'q}}}}

\newcommand{\STot}{\ensuremath{{\text{ST12}^{sn}_{l'q}}}}
\newcommand{\STtt}{\ensuremath{{\text{ST23}^{sn}_{l'q}}}}
\newcommand{\deltahh}{{\Delta^{sn}_{l'q}}}
\newcommand{\betao}{{\beta^{s}_{l'}}}
\newcommand{\betat}{{(1-\gamma^{s}_{l'})}}
\newcommand{\betatt}{{\gamma^{s}_{l'}}}

\newcommand{\tpqqh}{\ensuremath{CL^{n}_{Hf}}}
\newcommand{\tpqql}{\ensuremath{CL^{n}_{Lf}}}

\newcommand{\arpi}{\ensuremath{{\text{pa}^{sn}_{l'q}}}}
\newcommand{\alti}{\ensuremath{{\text{ai}^{sn}_{l'q}}}}

\newcommand{\dwb}{\ensuremath{{\text{dwb}^{sn}_{l'q}}}}

\newcommand{\dwbozo}{\ensuremath{{\text{dwb1z1}^{sn}_{l'q}}}}
\newcommand{\dwbozto}{\ensuremath{{\text{dwb2z1}^{sn}_{l'q}}}}
\newcommand{\dwboztho}{\ensuremath{{\text{dwb3z1}^{sn}_{l'q}}}}
\newcommand{\dwbozthf}{\ensuremath{{\text{dwb3z12}^{sn}_{l'q}}}}
\newcommand{\dwboztht}{\ensuremath{{\text{dwb3z2}^{sn}_{l'q}}}}
\newcommand{\dwbozt}{\ensuremath{{\text{dwb1z2}^{sn}_{l'q}}}}
\newcommand{\dwbozth}{\ensuremath{{\text{dwb1z3}^{sn}_{l'q}}}}
\newcommand{\dwbiozo}{\ensuremath{{\text{dwbi1z1}^{sn}_{l'q}}}}
\newcommand{\dwbiozto}{\ensuremath{{\text{dwbi2z1}^{sn}_{l'q}}}}
\newcommand{\dwbioztho}{\ensuremath{{\text{dwbi3z1}^{sn}_{l'q}}}}
\newcommand{\dwbiozt}{\ensuremath{{\text{dwbi1z2}^{sn}_{l'q}}}}
\newcommand{\dwbiozth}{\ensuremath{{\text{dwbi1z3}^{sn}_{l'q}}}}

\newcommand{\dwbi}{\ensuremath{{\text{dwbi}^{sn}_{l'q}}}}
\newcommand{\tbop}{\ensuremath{{\text{tbO}^{sn}_{l'q}}}}
\newcommand{\tbopozo}{\ensuremath{{\text{tbO1z1}^{sn}_{l'q}}}}
\newcommand{\tbopoztot}{\ensuremath{{\text{tbO2}^{sn}_{l'q}}}}
\newcommand{\tbopoztoth}{\ensuremath{{\text{tbO3}^{sn}_{l'q}}}}
\newcommand{\tbopozt}{\ensuremath{{\text{tbO1z2}^{sn}_{l'q}}}}
\newcommand{\tbopozth}{\ensuremath{{\text{tbO1z3}^{sn}_{l'q}}}}

\newcommand{\Ewait}{\ensuremath{{\text{Ewait}^{sn}_{{l'}q}}}}
\newcommand{\Ewaitp}{\ensuremath{{\text{Ewait}^{sn}_{{l}p}}}}
\newcommand{\TEwait}{\ensuremath{{\text{TEwait}^{sn}_{{l'}q}}}}
\newcommand{\TEwaitp}{\ensuremath{{\text{TEwait}^{sn}_{{l}p}}}}

\newcommand{\EwaitII}{\ensuremath{{\text{EwaitII}^{sn}_{{l'}q}}}}

\newcommand{\alightb}{{\ensuremath{\text{alightb}^{sn}_{l'q}}}}
\newcommand{\alightbi}{{\ensuremath{\text{alightbi}^{sn}_{l'q}}}}
\newcommand{\alightE}{{\ensuremath{\text{alightE}^{sn}_{l'q}}}}
\newcommand{\dwIzo}{{\ensuremath{\text{dwtI3z1}^{sn}_{l'q}}}}
\newcommand{\dwIE}{{\ensuremath{\text{dwtE}^{sn}_{l'q}}}}

\newcommand{\dwIzt}{{\ensuremath{\text{dwtI3z2}^{sn}_{l'q}}}}
\newcommand{\dwIzth}{{\ensuremath{\text{dwtI3z3}^{sn}_{l'q}}}}

\newcommand{\servELf}{{\ensuremath{\text{servELf}^{sn}_{l'q}}}}

\newcommand{\aldnns}{\ensuremath{{ad}^{sn}_{l'q}}}

\newcommand{\dwsi}{{\ensuremath{\text{dwtI}^{sn}_{l'q}}}}

\newcommand{\Eso}{\ensuremath{{\text{Es1}^{sn}_{{l'}q}}}}
\newcommand{\Esoi}{\ensuremath{{\text{Es1i}^{sn}_{{l'}q}}}}

\newcommand{\Est}{\ensuremath{{\text{Es2}^{sn}_{{l'}q}}}}
\newcommand{\Estzo}{\ensuremath{{\text{Es2z1}^{sn}_{{l'}q}}}}
\newcommand{\Esthzot}{\ensuremath{{\text{Es3z12}^{sn}_{{l'}q}}}}
\newcommand{\Esth}{\ensuremath{{\text{Es3}^{sn}_{{l'}q}}}}

\newcommand{\Esthzoth}{\ensuremath{{\text{Es3z3}^{sn}_{{l'}q}}}}
\newcommand{\Esthzothi}{\ensuremath{{\text{Es3iz3}^{sn}_{{l'}q}}}}
\newcommand{\Estzoi}{\ensuremath{{\text{Es2iz1}^{sn}_{{l'}q}}}}
\newcommand{\Esthzoti}{\ensuremath{{\text{Es3iz12}^{sn}_{{l'}q}}}}
\newcommand{\Estzt}{\ensuremath{{\text{Es2z2}^{sn}_{{l'}q}}}}
\newcommand{\Estzti}{\ensuremath{{\text{Es2iz2}^{sn}_{{l'}q}}}}
\newcommand{\yqo}{\ensuremath{{\text{Q1}^{sn}_{lp{l'}q}}}}
\newcommand{\ygo}{\ensuremath{{\text{G1}^{sn}_{lp{l'}q}}}}
\newcommand{\yqot}{\ensuremath{{\text{Q12}^{sn}_{lp{l'}q}}}}
\newcommand{\ygot}{\ensuremath{{\text{G12}^{sn}_{lp{l'}q}}}}
\newcommand{\yqt}{\ensuremath{{\text{Q2}^{sn}_{lp{l'}q}}}}
\newcommand{\ygt}{\ensuremath{{\text{G2}^{sn}_{lp{l'}q}}}}
\newcommand{\yyso}{\ensuremath{{\text{Y1}^{sn}_{lp{l'}q}}}}
\newcommand{\yyst}{\ensuremath{{\text{Y2}^{sn}_{lp{l'}q}}}}
\newcommand{\yysthzo}{\ensuremath{{\text{Y3z1}^{sn}_{lp{l'}q}}}}
\newcommand{\yysthzt}{\ensuremath{{\text{Y3z2}^{sn}_{lp{l'}q}}}}
\newcommand{\tbso}{\ensuremath{{\text{dtb1}^{sn}_{l'q}}}}
\newcommand{\yyyyso}{\ensuremath{{\text{YY1}^{sn}_{lp{l'}q}}}}

\newcommand{\yyyyst}{\ensuremath{{\text{YY2}^{sn}_{lp{l'}q}}}}
\newcommand{\yyyysthzo}{\ensuremath{{\text{YY3z1}^{sn}_{lp{l'}q}}}}
\newcommand{\yyyysthzt}{\ensuremath{{\text{YY3z2}^{sn}_{lp{l'}q}}}}
\newcommand{\yyyystc}{\ensuremath{{\text{YY2c}^{sn}_{lp{l'}q}}}}
\newcommand{\yyyysthc}{\ensuremath{{\text{YY3c}^{sn}_{lp{l'}q}}}}

\newcommand{\ttboso}{\ensuremath{{\text{tb1}^{sn}_{l'q}}}}
\newcommand{\ttbosop}{\ensuremath{{\text{ptb1}^{sn}_{l'q}}}}
\newcommand{\ttbost}{\ensuremath{{\text{tb2}^{sn}_{l'q}}}}
\newcommand{\ttbostp}{\ensuremath{{\text{ptb2}^{sn}_{l'q}}}}
\newcommand{\ttbosth}{\ensuremath{{\text{tb3}^{sn}_{l'q}}}}
\newcommand{\ttbosthp}{\ensuremath{{\text{ptb3}^{sn}_{l'q}}}}
\newcommand{\zyso}{\ensuremath{{\text{ZY1}^{sn}_{lp{l'}q}}}}
\newcommand{\yynso}{\ensuremath{{\text{Y1}^{sn}_{lp{l'}(q-1)}}}}
\newcommand{\tyso}{\ensuremath{{\text{TY1}^{sn}_{lp{l'}q}}}}

\newcommand{\syso}{\ensuremath{{\text{SY1}^{sn}_{lp{l'}}}}}

\newcommand{\servLo}{{\ensuremath{\text{servL1}^{sn}_{l'q}}}}
\newcommand{\servLt}{{\ensuremath{\text{servL2}^{sn}_{l'q}}}}

\newcommand{\pddo}{{\ensuremath{\text{pdd}^{sn}_{l'q}}}}
\newcommand{\pddop}{{\ensuremath{\text{pdd}^{sn}_{lp}}}}
\newcommand{\vtdo}{{\ensuremath{\text{vtd}^{sn}_{l'q}}}}
\newcommand{\vtdop}{{\ensuremath{\text{vtd}^{sn}_{lp}}}}

\newcommand{\ADiff}{{\ensuremath{\text{ADiff}^{sn}_{l'q}}}}
\newcommand{\RDiff}{{\ensuremath{\text{RDiff}^{sn}_{l'q}}}}
\newcommand{\RDiffp}{{\ensuremath{\text{RDiff}^{sn}_{lp}}}}
\newcommand{\Rths}{{\ensuremath{\text{Rths}^{sn}_{l'q}}}}
\newcommand{\Aths}{{\ensuremath{\text{Aths}^{sn}_{l'q}}}}
\newcommand{\Rthsi}{{\ensuremath{\text{RthsI}^{sn}_{l'q}}}}
\newcommand{\Athsi}{{\ensuremath{\text{AthsI}^{sn}_{l'q}}}}

\newcommand{\Ivdd}{{\ensuremath{\text{Ivdd}^{sn}_{l'q}}}}
\newcommand{\Ivddp}{{\ensuremath{\text{Ivdd}^{sn_{prev}}_{l'q}}}}

\newcommand{\floor}[1]{\lfloor #1 \rfloor}

\makeatletter
\def\ps@pprintTitle{%
  \let\@oddhead\@empty
  \let\@evenhead\@empty
  \def\@oddfoot{\reset@font\hfil\thepage\hfil}
  \let\@evenfoot\@oddfoot
}
\makeatother

\allowdisplaybreaks
\begin{document}\sloppy
\setlength{\parindent}{2em}

\begin{frontmatter}

\title{A Comprehensive Stochastic Programming Model for Transfer Synchronization in Transit Networks} 

\author[1]{Zahra Ansarilari}
\ead{zahra.ansarilari@mail.utoronto.ca}

\author[2]{Merve Bodur\corref{cor1}
}
\ead{merve.bodur@ed.ac.uk}

\author[1]{Amer Shalaby}
\ead{amer.shalaby@utoronto.ca}

\cortext[cor1]{Corresponding author}
\address[1]{Department of Civil and Mineral Engineering, University of Toronto, Canada}
\address[2]{School of Mathematics, University of Edinburgh, Edinburgh, UK}

\begin{abstract}
We investigate the stochastic transfer synchronization problem, which seeks to synchronize the timetables of different routes in a transit network to reduce transfer waiting times, delay times, and unnecessary in-vehicle times. We present a sophisticated two-stage stochastic mixed-integer programming model that takes into account variability in passenger walking times between bus stops, bus running times, dwell times, and demand uncertainty. Our model incorporates new features related to dwell time determination by considering passenger arrival patterns at bus stops which have been neglected in the literature on transfer synchronization and timetabling. 
We solve a sample average approximation of our model using a problem-based scenario reduction approach, and the progressive hedging algorithm.
As a proof of concept, our computational experiments on two single transfer nodes in the City of Toronto, with a mixture of low- and high-frequency routes, demonstrate the potential advantages of the proposed model. Our findings highlight the necessity and value of incorporating stochasticity in transfer-based timetabling models.
\end{abstract}
\begin{keyword}
Transfer synchronization \sep Timetabling \sep Bus running time and passenger demand uncertainty \sep Passenger arrival patterns \sep Dwell time determination \sep Stochastic programming 
\end{keyword}
\end{frontmatter}
\vspace{-8pt}
\section{Introduction}
To deliver fast, convenient, and reliable transfers to transit users, special efforts must be made in both schedule development and operational control due to the inherent stochastic nature of transit systems, the disutility of transferring and the high penalty of missing a connecting bus. A growing number of studies have tackled the \emph{transfer synchronization problem} \citep{liu2021review} but mostly in its deterministic setting; and the few that tackled the stochastic version of the problem have not considered in detail the impacts of transit systems' \emph{stochasticity} on successful transfer synchronization. 
In real-world practice, transit agencies commonly build timetables without accounting adequately for stochastic operations and with slight or no consideration of synchronization among intersecting routes at transfer nodes. 
Given the importance of transfers for many transit users, it is imperative to develop seamlessly synchronized timetables that take into account the stochastic nature of transit operations.

In contrast to unpredictable and typically non-recurrent sources of uncertainties in a transit system, such as major road accidents or severe weather conditions, the \emph{uncertainty of bus running times and passenger demand} due to recurrent perturbations can be estimated based on historical data, and thus proactively considered at the scheduling stage.
Although stochastic bus running times have been considered in previous transit scheduling research, the variability of passenger demand has received little attention in transfer synchronization studies \citep{liu2021review}.
Two primary approaches have been used in practice and previous studies to tackle variable running times: (1) taking the average or arbitrary percentile values of running times, and (2) incorporating slack times into timetables to increase the possibility of on-time performance of buses. Nevertheless, both strategies have drawbacks. In the first approach, many trips may not meet planned timetables, particularly when running times deviate substantially from their average. In the second approach, long slack times could result in longer passenger travel time and higher transit operating cost. 

Moreover, most previous research has overlooked \emph{essential modelling components} regarding the influence of different stochastic factors in transit systems that should be jointly considered when developing timetables. For instance, bus \emph{dwell time determination} should be properly and explicitly formulated based on the stochastic components involved in service time for boarding and alighting passengers. Otherwise, timetables generated with given dwell times would be susceptible to failure. This time interval, which is not accounted for during timetabling, may be small for a bus at a single transfer point, but it might lead to significant delays in the vehicle's arrival times at subsequent stops due to the stochastic nature of transit operations such as varying bus running times. Hence, established timetables and intended successful transfers may not occur if an accurate computation of bus dwell time is not integrated during the scheduling phase. Moreover, in stochastic transfer synchronization models, the suitable treatment of dwell time is crucial when slack times are incorporated in timetables. 

To address the above-mentioned gaps in timetabling and primarily transfer synchronization, we propose a two-stage stochastic mixed-integer programming (MIP) model to develop synchronized timetables that account for demand variability and passenger arrival patterns at bus stops 
as well as variability in passenger walking times between transit stops, bus running times, and dwell times. Our main contributions can be summarized as follows:

\vspace*{-0.3cm}
\begin{itemize}[leftmargin=0.32cm]
    \item We propose a novel stochastic bus timetabling model focusing on transfer synchronization. In particular:  
    \begin{itemize}
    \item We present a detailed formulation for determining bus dwell times in a stochastic setting based on the number of alighting and boarding passengers with various arrival patterns of local passengers at bus stops (random arrivals for high-frequency lines and distribution-based arrivals for low-frequency lines). 
    \item We introduce new formulation for different types of waiting time with their associated penalized demand: (1) transfer waiting time, (2) unnecessary in-vehicle time when a bus remains at a stop without passenger service till its timetable departure time, and (3) delay time due to the late departure time of a bus compared to its timetabled time.
    \end{itemize}
    \item We propose a heuristic solution method using a problem-based scenario reduction approach and Progressive Hedging (PH) algorithm to efficiently reach near-optimal solutions. The scenario reduction approach substantially alleviates the main challenge of solving our sophisticated stochastic MIP model. 
    \item We conduct several experiments to demonstrate: (1) the benefit of detailed transfer synchronization formulation, (2) the value of our stochastic transfer synchronization model, 
    and (3) the efficiency of our solution method.
\end{itemize}

{{The rest of this paper is organized as follows. In Section \ref{sec:SLiteratureReview}, we review the transfer synchronization literature in detail focusing on stochastic models. In Section~\ref{sec:SOBJ and CONS}, we provide an overview of our modelling framework, assumptions, and the transfer synchronization problem with regards to different types of successful transfers based on passengers' arrival patterns followed by our mathematical formulation in Section \ref{SMIPS}.
In Section \ref{sec:SMethodologyS}, we present our solution method to solve the sample average approximation of our proposed stochastic programming model, based on a problem-based scenario reduction approach and the PH algorithm. In Section \ref{sec:Experiment}, we present and discuss the results of our numerical examples examining different aspects of our model formulation and its benefits.
Finally, in Section \ref{sec:SConclusion}, we conclude the paper by summarizing our main findings, takeaway messages, and future research recommendations. }}

\section{Literature Review}
\label{sec:SLiteratureReview}
Over the past few decades, numerous studies have undertaken transfer time analysis and optimization at both the scheduling and operational stages. At the scheduling stage, the majority of studies \citep{abdolmaleki2020transit, ansarilari4043348novel, ansarilari2022transfer, ceder2001creating, chu2019models, fouilhoux2016valid, ibarra2014integrated, ibarra2012synchronization, wu2016multi} considered deterministic conditions as the main underlying assumption. In other words, passenger demand and their route choice (i.e., transferring and non-transferring passengers) and time-related components (such as bus running times, dwell times, and passenger walking times between the platforms) are considered fixed and given as inputs to optimization models. 
However, deterministic approaches to modelling transit systems with stochastic elements have serious limitations
\cite{bookbinder1992transfer}, 
specifically for the transferring process. 
Ideally, the construction of synchronized timetables should take into account the effects of predictable operational variations in the system to safeguard against failed connections. To do so, some studies have focused on the inherent stochasticity in the transfer process, and how it is affected by variations in traffic conditions, running times, dwell times, service frequencies and demand distributions. \citet{gkiotsalitis2021stop} review timetabling studies in public transport, discussing these considerations in detail.
Below, we review the most relevant literature on stochastic transfer synchronization problems.

\subsection{Analytical Models for Stochastic Transfer Synchronization}
Due to the complexity of stochastic transfer synchronization models, some studies used approximations via analytical modelling approaches in which solutions are mostly obtained in closed form by solving a system of equations and doing basic calculations. As one of the first studies, \citet{lee1991optimal} formulated a stochastic model for a {{Time Transfer System (TTS)}} design with rail-bus transfer synchronization to determine optimal slack times. Slack times are applied to compensate for the variability in bus running times and to reduce the probability of missed transfers. The proposed approach is useful only when the standard deviation of bus arrivals is greater than a certain threshold. The considered objective  included three main costs: scheduled delay cost (a penalty for drivers and non-transferring passengers, resulting from holding a bus/train due to the scheduled slack time), missed connections from bus to train or from train to bus, and delay incurred by transferring passengers if bus/train is late. 
Similarly, \citet{chowdhury2002intermodal} defined a cost function including both user and supplier costs to determine route headways and slack times for feeder-bus and train synchronization, accounting for variable bus arrival times. Three degrees of transfer synchronization with different associated costs were investigated: full, partial, and none.
The authors noted the benefit of using the model when transit systems have low-frequency lines with low standard deviation in bus arrival times and a large volume of transferring passengers. 
\citet{sivakumaran2012cost} explored headway synchronization of feeder and trunk routes through Pareto-improving approach in which bus user cost and operational cost are minimized. 

Using integer-ratio headways is one of the common approaches in analytical models to generate synchronized timetables. The determined headways are integer multipliers of the base cycle time, i.e., round fraction of 60 minutes.
\citet{ting2005schedule} showed the benefit of using integer-ratio headways for a multiple hub transit network of low-frequency lines with large headway variances. 
\citet{aksu2014transit} formulated a cost function including operational cost as well as in-vehicle, transfer and waiting time, and illustrated the benefit of integer-ratio headways. 
\citet{kim2014integration} compared common headways and integer-ratio headways in a combined fixed- and flexible-route transit with TTS design. The considered cost function included extra transfer waiting time cost due to slack time, missed connection, delayed connection, and inter-cycle waiting times. The decision variables were service type, vehicle size, number of zones, line headways, fleet size and slack times optimized through genetic algorithms. 

\citet{ting2005schedule} developed an optimization model to determine the optimal headways and slack times for multiple hubs in a transit network with TTS design. 
The inputs included vehicle arrival time distribution, uniformly distributed passenger arrival time, and fixed origin-destination matrix for a predetermined network. This is one of the very few studies which considered passenger arrival patterns explicitly. 
The transfer cost function included the cost of missed transfers by considering the joint arrival distribution of feeder and connecting buses. Additionally, the transfer cost function considered slack time related components: extra in-vehicle time for on-board passengers 
and extra vehicle operating cost. Due to the model's complexity, they 
combined an analytical solution with a step-by-step heuristic. 
Recently, \citet{yang2020schedule} developed a two-phase solution algorithm merging an adaptive genetic algorithm and an analytical optimization model of transfer synchronization while considering demand heterogeneity patterns in a feeder-trunk transfer node. 

Analytical transfer synchronization models 
are incapable of precisely representing cost components, especially user costs. 
This disadvantage affects the outcomes more adversely when models are solved in a stochastic setting involving slack time variables. 
To address this issue, another group of studies applied mathematical programming models. 

\subsection{Mathematical Programming Models for Stochastic Transfer Synchronization}
\citet{bookbinder1992transfer} combined a simulation platform and an optimization model to incorporate bus travel time stochasticity in transfer synchronization. 
They defined an interval for successful transferring by considering a joint distribution of feeder and connecting arrival times. The objective function minimized the impacts of bus running time randomness on mean disutility, where longer waiting times with higher probability are given higher objective weight. 
The mean disutility measure should be carefully selected, e.g., 
the authors stated that mean transfer waiting time can be a misleading measure as it does not consider a transfer's reliability due to penalizing both short and long transfer waiting times with the same weight.  
\citet{desilets1992syncro} showed the benefit of using fixed headways instead of flexible headways for transfer synchronization with variable bus arrival times through two stochastic optimization models, the first minimizing the expected transfer waiting time while the second minimizing the standard deviation of transfer waiting times,  
which are solved via a local search algorithm. 

\citet{teodorovic2005schedule} developed a model minimizing total transfer waiting times based on the deterministic model by \citet{klemt1988schedule} considering the number of transferring passengers as random variables. 
A solution method combining an ant colony system and fuzzy logic was applied. \citet{nair2013large} designed a two-stage stochastic programming model minimizing expected transferring passenger, which reduced waiting times by 26.38\% for a large-scale transit network in Washington, D.C., based on probabilistic information on transfer demand. However, the model did not account for the uncertainty in bus running times. In fact, it used the existing timetables, including pre-defined slack times, as input. The authors only proposed an optimal temporal shift to improve transferring passenger waiting times. 

\cite{lee2022path} developed a MIP model incorporating time-dependent travel times and addressing the complexities of transfers in multi-modal mobility framework. They introduced a novel concept of path-oriented scheduling and applied their model to three Copenhagen lines using simulation data. However, the model's assumptions include the independence of dwell times from passenger flow. Although, they considered three passenger arrival types (early, on-time, and late) in determining waiting time, the formulation lacks the acknowledgment that arrival patterns depend on line frequency. Furthermore, the study overlooks variations in waiting, delay, and in-vehicle times corresponding to different demand levels, limiting its scope. In fact, the objective function focuses solely on assigning different weights to waiting times.
 
\citet{dou2016time} proposed a multi-objective nonlinear MIP scheduling model to reduce transfer disutility and increase schedule adherence at time control points. The authors solved a robust equivalent of their model through an approximate deterministic model with Monte Carlo simulation for three bus lines and the mass rapid train transit system in Singapore. The model minimized the transfer waiting time from bus to train while considering stochastic bus running times by assigning appropriate slack time to each segment of a bus route with predetermined timing control points. 
\citet{gkiotsalitis2019robust} 
developed network-wide synchronized bus timetables via incorporating fluctuations of bus travel times and dwell times during daily operations and also considering the service regularity as the main performance indicator. Using a robust formulation approach, they minimized the worst-case scenario of missing transfers.

\citet{wu2016designing} introduced a new modelling approach by considering both real-time transfer coordination control and planned slack time at once in a stochastic MIP model. The model included a holding strategy designed to work as follows: if a feeder vehicle is late to arrive at a transfer node according to the published timetable but its predicted arrival time is within a given safety margin, a ready-to-leave connecting vehicle can be held, resulting in a departure time later than the timetable. However, the authors assumed that the held vehicle could make up the holding time before arriving at the next transfer node. The authors considered a scheduled-based system in which line headways are synchronized through common headways. Therefore, average transfer waiting times are approximated by a closed form formulation. 
\citet{wu2019stochastic} developed an extended transfer synchronization model considering not only the travel time variability of buses but also passengers' rerouting, both pre-trip and on-trip, and stochastic demand assignment in the system. They applied a bi-level programming model in which headways and slack times are determined simultaneously with passengers' route choices through two travel strategies: non-adaptive and adaptive routing. The novelty of the study is the consideration of passenger route modifications because of missed transfers, which affect the total system performance. 
Their model maximized transfer reliability which is the expected probability of successful transferring.

\citet{ansarilari4043348novel} proposed a MIP model which incorporates key details of the transfer process such as pre-planned holding times, vehicle capacities, and most notably dwell time determination considering the actual number of successful boarding passengers for the first time in the literature. They developed a Lagrangian-based heuristic solution method. Via detailed numerical experiments, they demonstrated that more realistic modelling of the transfer process can  significantly improve the quality of the obtained timetables. Motivated by this important fact, our work aims to propose a much more realistic mathematical model for the transfer synchronization problem incorporating main sources of uncertainty among other features. More detailed conceptual comparisons of our work with \cite{ansarilari4043348novel} are provided in the subsequent sections, and a numerical comparison is presented in our numerical results section. 

Our review identifies two significant gaps in the literature.
First, most studies did not explicitly formulate transfer synchronization; instead, they used closed-form formulas or the TTS design to assign slack times to mitigate the variability in bus running times. As such, the proposed objective functions 
did not accurately represent the costs associated with passengers and agencies, e.g., passengers penalized for delayed bus departure times or unnecessary in-vehicle time due to excess slack times have not been properly accounted for in the models. Second, passenger demand uncertainty and passenger arrival patterns at bus stops, as well as their dual impact on the variability of bus dwell times and the reliability of bus arrival/departure times, have not been adequately examined. To address these limitations, we devise a comprehensive transfer synchronization model, where most notably, (i) we determine bus dwell times based on the arrival patterns of local passengers and the arrival time of transferring passengers at bus stops relative to uncertain bus arrival times and timetable departure times, (ii) determine the number of affected passengers with different types of waiting times, namely, transfer waiting, delayed departure, and unnecessary in-vehicle times, and (iii) track different types of passengers for low- and high-frequency lines separately.

\section{Problem Description}
\label{sec:SOBJ and CONS}


We consider a pre-designed transit network with a set of transfer nodes \textit{N} and a set of bus lines \textit{L} with known headways, which are typically determined based on ridership levels. The main output of our model is bus timetable departure times from terminals and transfer nodes, which would be published and made available to passengers. In our two-stage stochastic MIP model, these constitute the first-stage, i.e., deterministic, decisions.
All the other variables such as bus arrival times, bus dwell times, and waiting times would be second-stage decisions and thus be allowed to take different values for each  future scenario (i.e., the realization of uncertain parameters). Thus, any stated variable other than timetable departure times implies we are referring to a specific scenario. Our objective is to minimize transfer waiting times, unnecessary in-vehicle times, and delay times with their associated demand while taking into account the stochasticity of passenger walking times between bus stops, bus running times, and bus dwell times in addition to passengers' arrival patterns at bus stops and demand variability.

\begin{figure}[htb]
\begin{center}
\centerline{\includegraphics[width=0.75\linewidth]{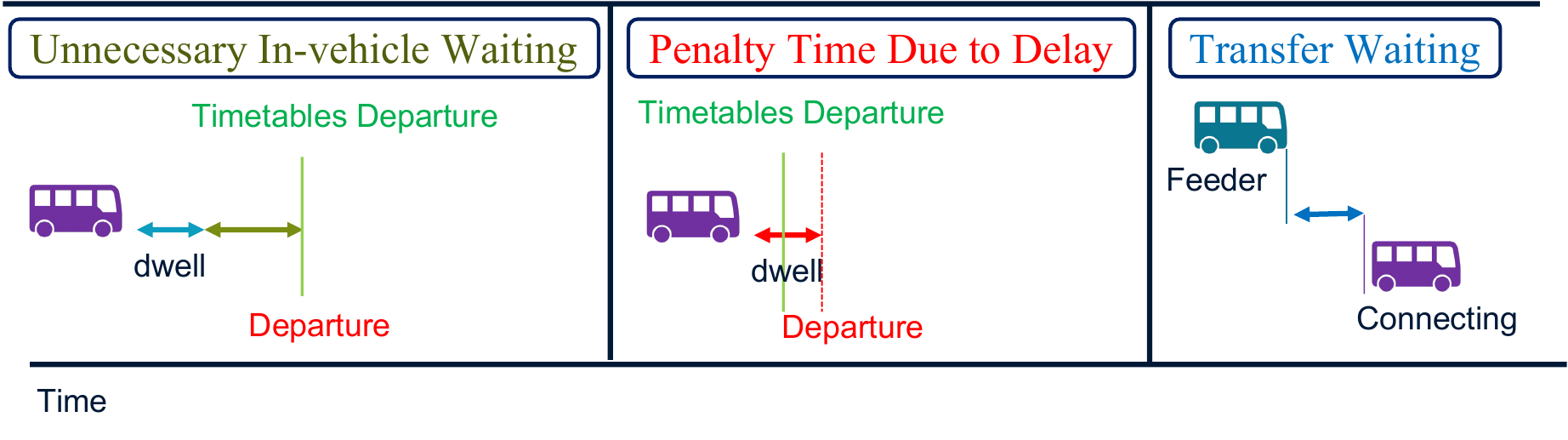}}
\caption{Different types of waiting times.}
\label{fig:Wait}
\end{center}
\end{figure}

Three types of waiting times can occur as shown in Figure \ref{fig:Wait}.
When a bus is held without any service, i.e., being ready to depart but must wait until its timetable departure time, passengers are subjected to 
unnecessary in-vehicle time. On the other hand, if a bus requires a longer dwell time to provide service to passengers and thus passes its timetable departure time, or if a bus arrives later than its timetable departure time even with a short dwell time, the bus departure time would be later than its timetable departure time, resulting in a delay penalty for passengers.
Finally, transfer waiting time is the time interval between transferring passengers' arrival at their connecting bus stop and the arrival of their bus. Hence, there is no transfer waiting time if transferring passengers arrive at the bus stop when the connecting bus is already there.

The number of people boarding and alighting determines the required dwell time for a bus. Boarding passengers include both transferring and local passengers. The number of passengers and their relative arrival at a bus stop compared to the arrival of the connecting bus and its timetable departure time impact how the required dwell time of a bus and consequently the bus departure time and different waiting times are derived in our model. In what follows, we discuss our assumptions and the key concepts of our modelling approach regarding the types of successful transfers, different types of passengers' waiting times, and different arrival patterns of local passengers.

\subsection{Assumptions}
\label{sec:SAssumptions}
We make some assumptions based on the availability of detailed historical data. The following inputs are known at the scheduling stage and assumed identical and fixed during the planning horizon for all scenarios: 
\begin{itemize}[leftmargin=1.2em]\setItemSep{-0.1cm}
  \item[$\diamond$] Line headways; however, they could vary within a specified given range (\hlmin~ and \hlmax~ for line $l$) chosen by the operating agency.
  \item[$\diamond$] Number of bus trips over the planning horizon, and bus sizes.
  \end{itemize}
The following inputs are given for each individual scenario, whose values are derived from historical data and a statistically independent probability distribution associated with each input.
  \begin{itemize}[leftmargin=1.2em]\setItemSep{-0.1cm}
  \item[$\diamond$] The number of in-vehicle passengers on the buses departing from the terminals of each line. (For the other nodes, the number of in-vehicle passengers are calculated in our model based on the counts of alighting and boarding passengers at previous nodes).
  \item[$\diamond$] 
  The total net number of boarding and alighting passengers at stops between the terminal and the first transfer node, as well as between two consecutive transfer nodes.
  \item[$\diamond$] 
  The number of alighting passengers and transferring passengers when the bus arrives at a stop.
  \item[$\diamond$] The number of local passengers (i.e., the number of passengers walking to a transfer node) for each bus at each transfer node for a low-frequency line;
  \item[$\diamond$] The arrival rates of local passengers for trips of high-frequency lines.
  \item[$\diamond$] The passenger walking time between bus stops (from a feeder bus to a connecting bus) and bus running times from node to node.
  \end{itemize}
Moreover, we assume that buses have enough capacity to satisfy demand. The boarding and alighting doors are separate. We consider a schedule-based system, thus a bus cannot leave earlier than its timetable departure time. There is no uncertainty at terminals and buses depart from the terminals based on their timetable departure times. Lastly, real-time information on bus operation (e.g., delayed buses) is unavailable, and passengers do not adjust their path/mode choice. 
  
\subsection{Our New Approach of Transfer Synchronization Formulation in Stochastic Setting}
\label{subsec:situationsinstochastic}
As mentioned above, our main contributions are the new approach of  successful transfer formulation and bus dwell time determination for high- and low-frequency lines in a stochastic setting. In what follows, we first explain the three main types of successful transfers. Then, we separately explain bus dwell time determination based on passengers' arrival patterns, including local and transferring passengers, for high- and low-frequency lines.

\subsubsection{Types of Successful Transfer}
\label{sec:Types}
\begin{figure}[!ht]
\begin{center}
\centerline{\includegraphics[width=0.75\linewidth]{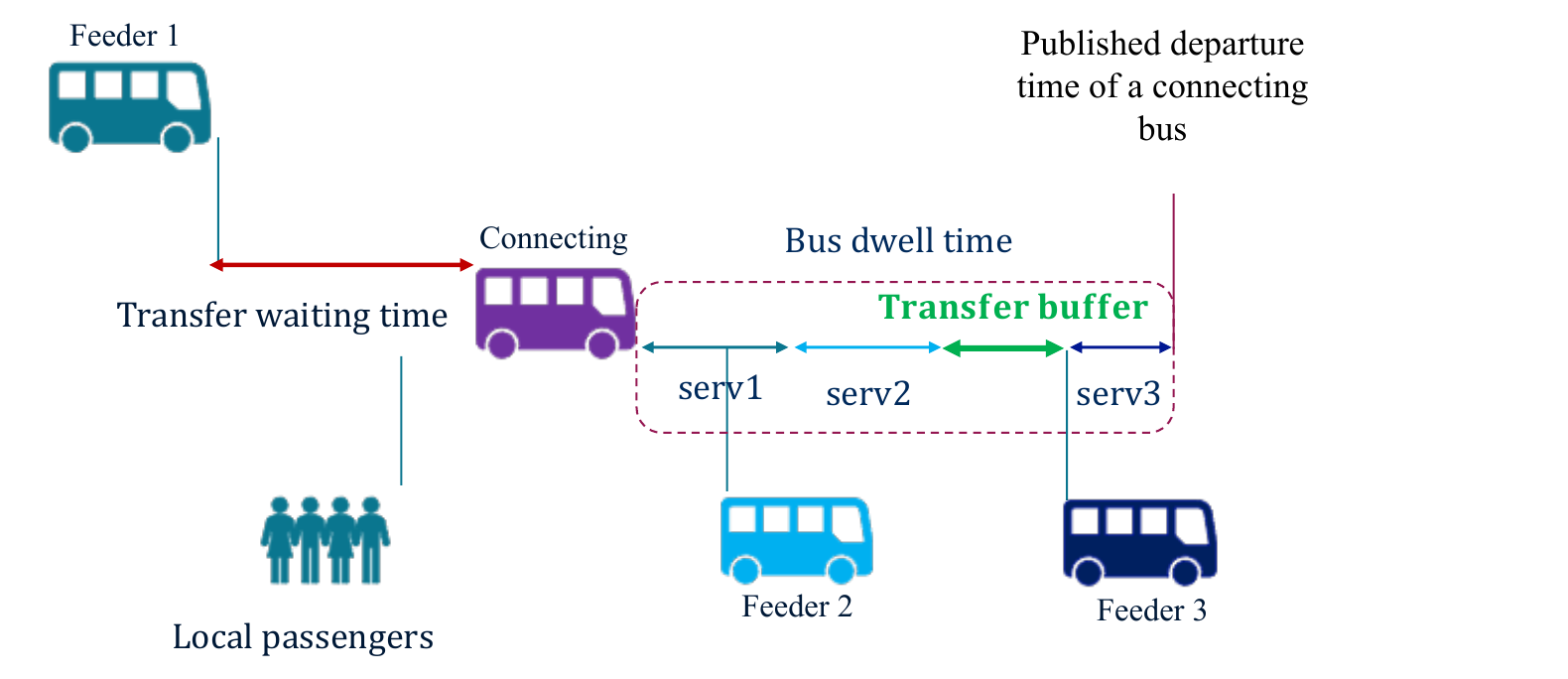}}
\caption{The process of three types of a successful transfer.}
\label{fig:Transfer}
\end{center}
\end{figure}

A successful transfer is achieved if the arrival time of transferring passengers at a bus stop is before the departure time of the connecting bus. The approach to formulating successful transfer synchronization in this study is based on deterministic transfer waiting time optimization work in  \citep{ansarilari4043348novel}, as shown in Figure \ref{fig:Transfer}. 

\begin{itemize}\setItemSep{-0.2cm}
    \item Successful transfer Type 1 occurs if the arrival time of transferring passengers at the connecting bus stop is earlier than or at the arrival time of their connecting bus. In this case, the transfer waiting time, denoted by NTwait, is the duration between the arrival time of transferring passengers at the bus stop and the arrival time of their connecting bus. 
    \item Successful transfer Type 2 occurs when transferring passengers arrive during the boarding time, denoted by serv1, of the passengers who arrived at the stop beforehand, i.e., local or successful transferring passengers of Type 1 (from another feeder). We denote by serv2  the boarding time of passengers with successful transfer Type 2. 
    \item Successful transfer Type 3 happens after all previously arriving passengers at the stop, including local and successful transferring passengers (Types 1 and 2), have already boarded the connecting bus. It can occur either immediately after serv1 plus serv2 or with the aid of transfer buffer time, denoted by tb. It is important to note that transfer buffer time is 
    added to stopping time, i.e., amounting to a pre-planned holding time at the scheduling stage, to enhance successful transfers while considering the in-vehicle time penalty induced by passengers already onboard. 

\end{itemize}

That being said, in the deterministic model in \citep{ansarilari4043348novel}, the transfer buffer time was not introduced for safeguarding against stochasticity. 
Although in this paper we use a similar approach to consider successful transfer synchronization as in \citep{ansarilari4043348novel}, due to uncertainty in bus running times and different arrival patterns of local passengers, significant modifications need to be applied, specifically for transferring to low-frequency lines. In fact, the three types of successful transfer can occur several times for a low-frequency line if a bus arrives early. 

\subsubsection{Dwell Time Determination Based on Passengers' Arrival Pattern at Bus Stops}
\label{subsubsec:zonedefns}
Determining bus dwell times is essential, particularly in stochastic timetabling. The reliability of the generated bus departure times increases with how well the dwell time values are specified in the mathematical model. Moreover, the detailed dwell time formulation in a model improves the accuracy of waiting time costs and the corresponding affected demand. To calculate bus dwell times, we estimate the number of boarding passengers, both local and transferring, considering their arrival times at the bus stop relative to the arrival time of their bus and its timetable departure time. 
The suggested timetable departure times as the model's primary outcome combined with the system's stochasticity contribute to different possible situations of passenger arrivals and dwell time determinations in each scenario, which we categorize in zones and describe them in the following for high- and low-frequency lines separately. 

\begin{figure}[!ht]
\begin{center}
\centerline{\includegraphics[width=0.9\linewidth]{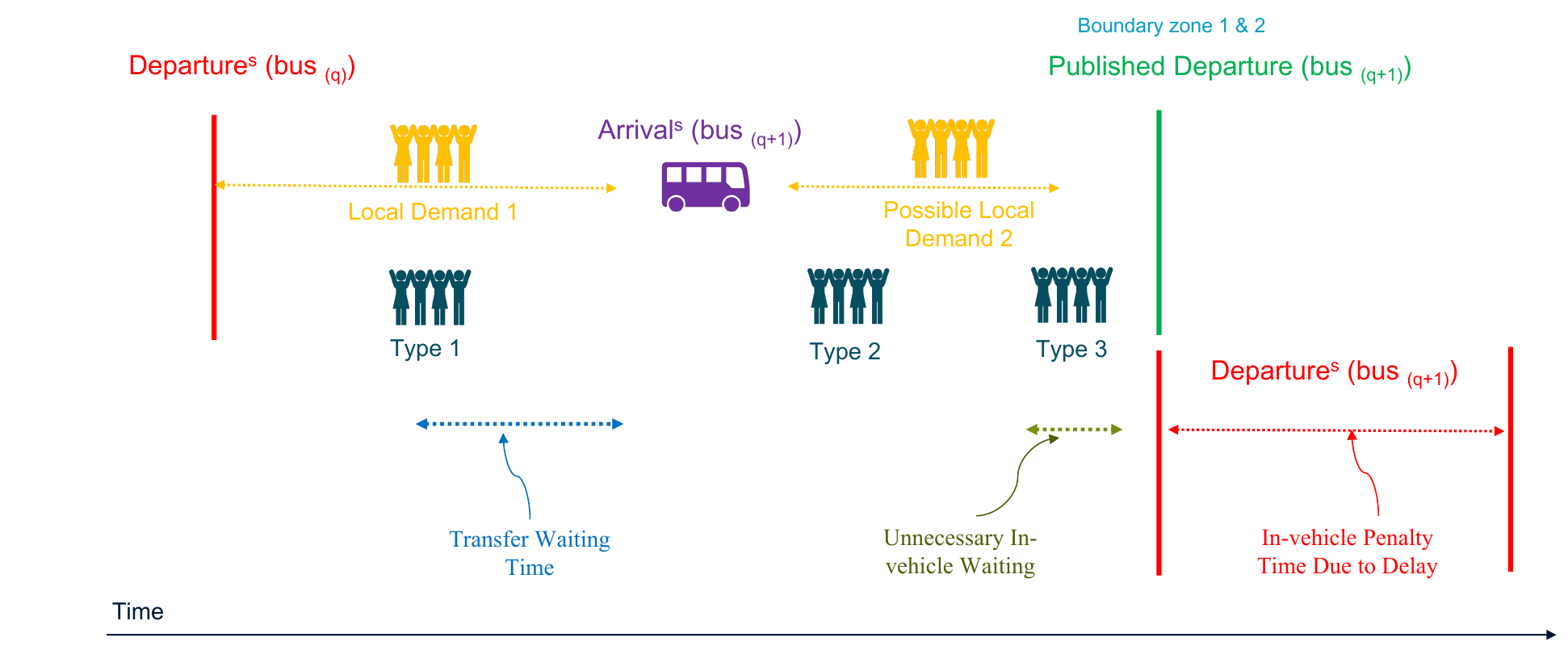}}
\caption{Local passengers arrival pattern for high-frequency line, Zone 1.}
\label{fig:HighZ1}
\end{center}
\end{figure}

\begin{figure}[!ht]
\begin{center}
\centerline{\includegraphics[width=0.85\linewidth]{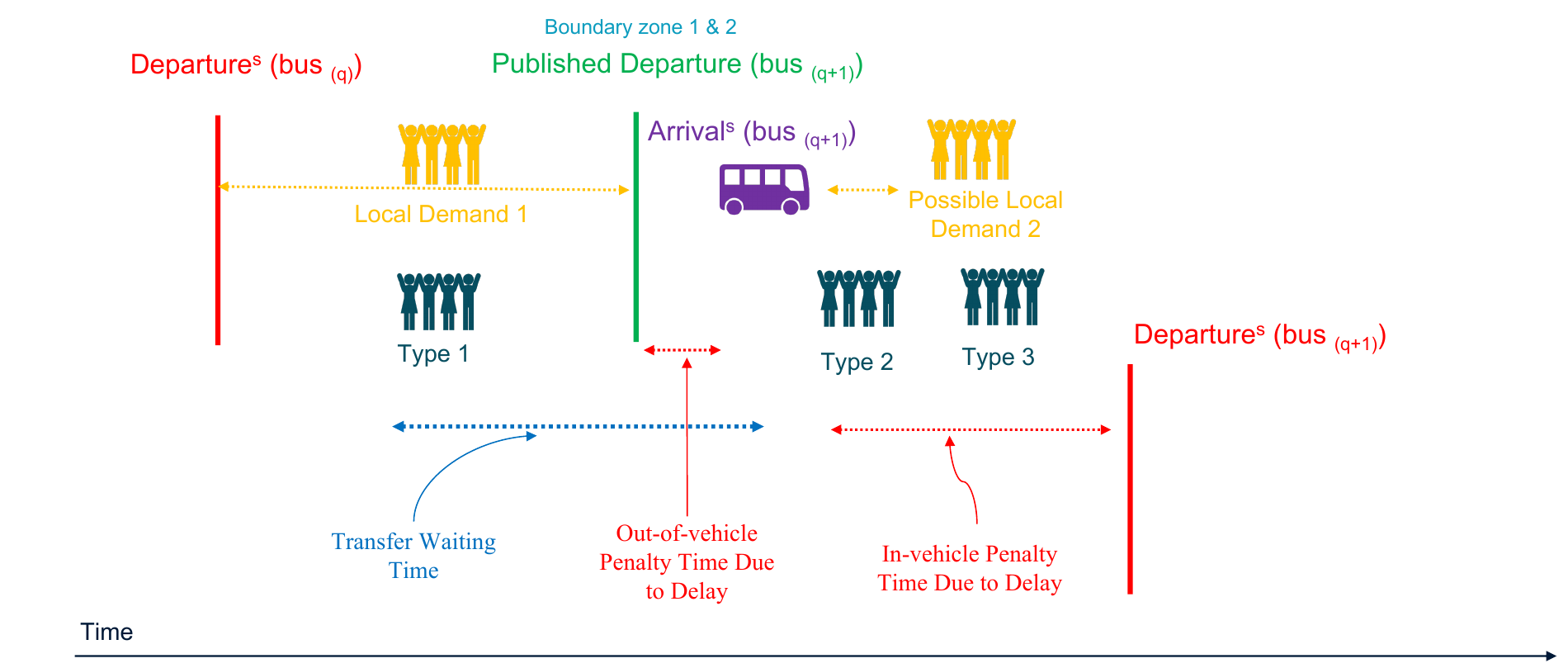}}
\caption{Local passengers arrival pattern for high-frequency line, Zone 2.}
\label{fig:HighZ2}
\end{center}
\end{figure}

\paragraph{High-frequency lines}
The arrival process of local passengers for a high-frequency line is commonly assumed to be random and independent of bus arrivals, even when published or real-time information is available \citep{chen2013implementation,delgado2012much,hounsell2012data}. As such, the number of local passengers for high-frequency lines could be calculated based on the mean arrival rate of passengers and the time difference between the previous bus departure time and the arrival time of the next bus. 
However, if a bus arrives early or late compared to its timetable departure time, more local passengers are expected to arrive, and their number would depend on the extent of the earliness or lateness of the bus.  
This is because our stochastic model considers variable arrival and departure times in various scenarios resulting in different numbers of local passengers and their required dwell times. We establish two zones for high-frequency lines to model the number of passengers and the bus dwell time properly: (1) when the bus arrives before its timetable departure time, Figure \ref{fig:HighZ1}, and (2) when the bus arrives after its timetable departure time, Figure \ref{fig:HighZ2}. Regardless of a bus arrival zone, the first group of local passengers arrives before the arrival of the bus. If a bus arrives in Zone 1, a second group of local passengers might also arrive during the bus stopping time until the timetable departure time. Additionally, if the bus arrives after its timetable departure time, Zone 2, a second group of local passengers would also arrive during the total service time of the first group of local passengers and successfully transfer passengers.

Each type of successful transfer can occur only once for a high-frequency line.
We define a new auxiliary variable, called dwell/stopping time budget, as the difference between the timetable departure time and the bus arrival time. If a bus arrives within Zone 1, based on the number of local passengers and successful transfer passengers of Types 1 and 2, and the second group of local passengers, we determine the required boarding time, serv1, serv2, and serv4, respectively. If the dwell time budget is longer, then the remaining dwell time budget would be considered for possible transfer buffer time, tb. If the transfer buffer time is assigned and we have a third group of successful transfers, with the service time of serv3, then we can determine the total boarding times plus transfer buffer time, compare it with the alighting time of passengers, and finally calculate the bus dwell time for such instance. It should be noted that if the bus arrives very close to its timetable departure time in Zones 1 or 3, successful transfer Type 3 can still occur, but with the value of a zero transfer buffer. 
Moreover, if the bus arrives at Zone 2, the delay penalty is considered in two ways: out-of-vehicle delay and in-vehicle delay. Therefore, we can determine the affected demand by the delay time more accurately and assign different weights as the perceived delay time is different if it is experienced in- or out-of-vehicle.

\begin{figure}[!ht]
\begin{center}
\centerline{\includegraphics[width=0.78\linewidth]{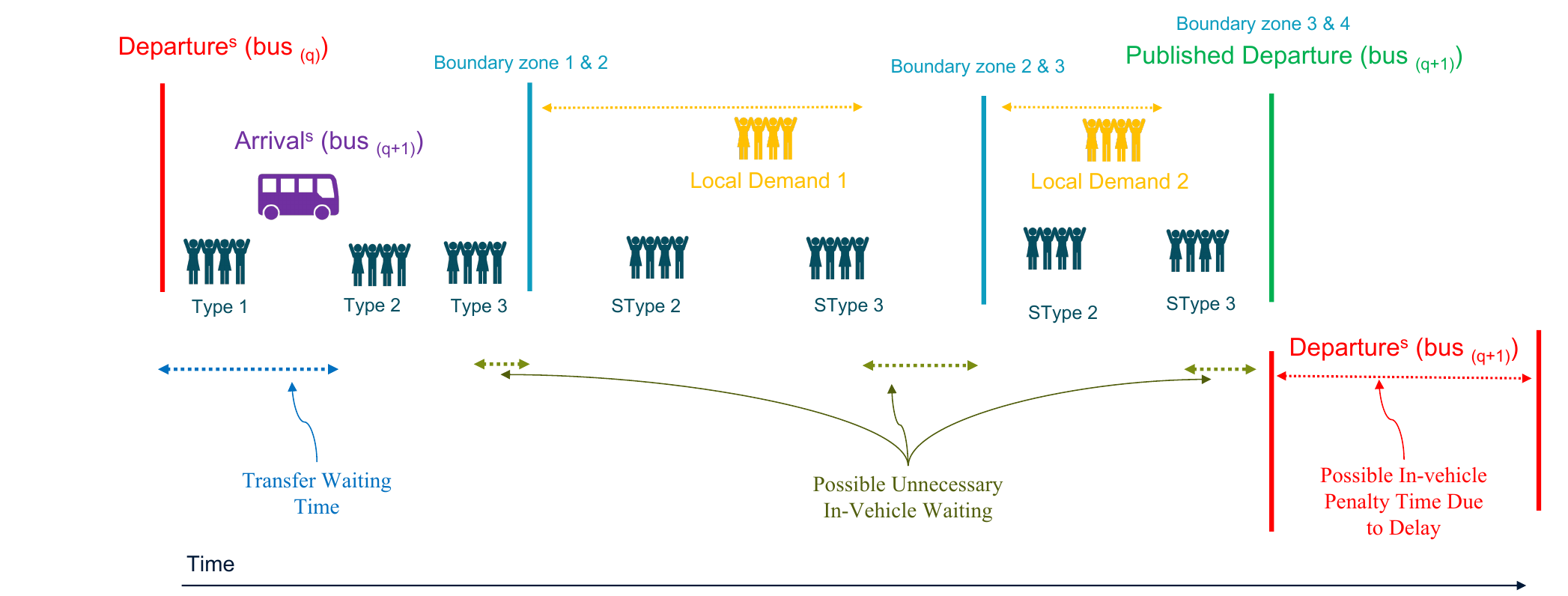}}
\caption{Local passengers' arrival pattern for low-frequency line, Zone 1.}
\label{fig:Low1}
\end{center}
\end{figure}

\begin{figure}[!ht]
\begin{center}
\centerline{\includegraphics[width=0.78\linewidth]{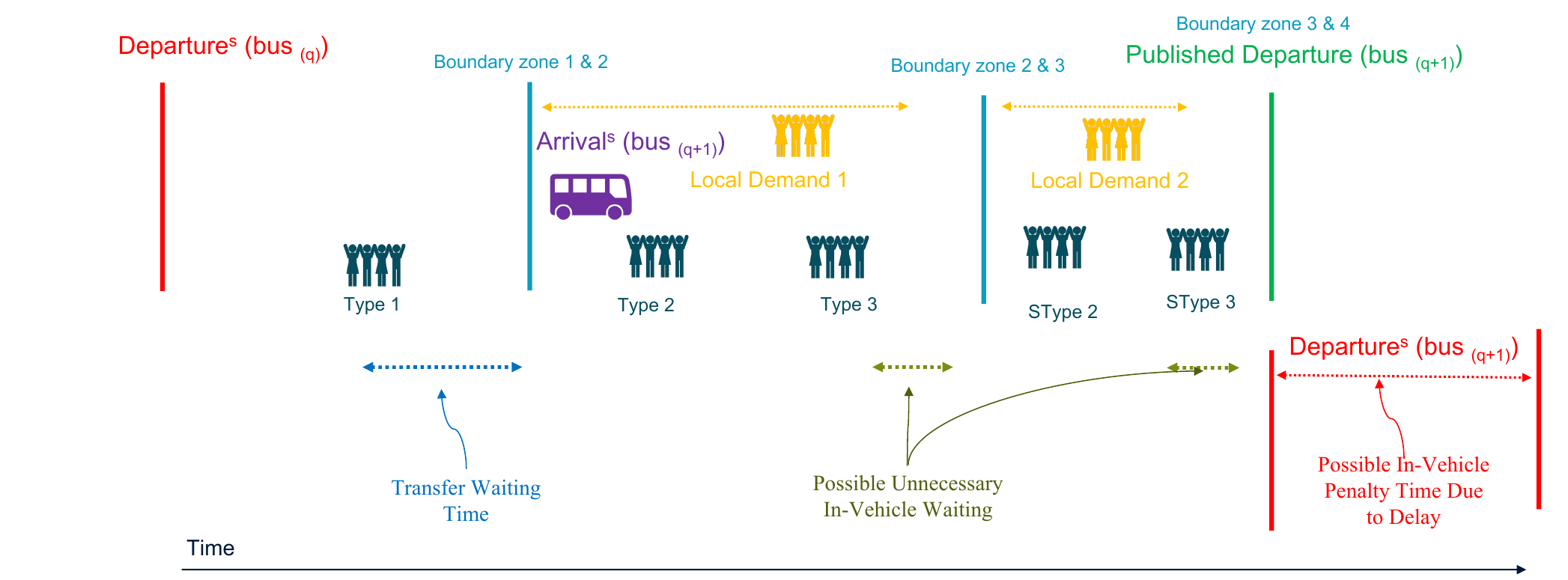}}
\caption{Local passengers' arrival pattern for low-frequency line, Zone 2.}
\label{fig:Low2}
\end{center}
\end{figure}

\begin{figure}[!ht]
\begin{center}
\centerline{\includegraphics[width=0.78\linewidth]{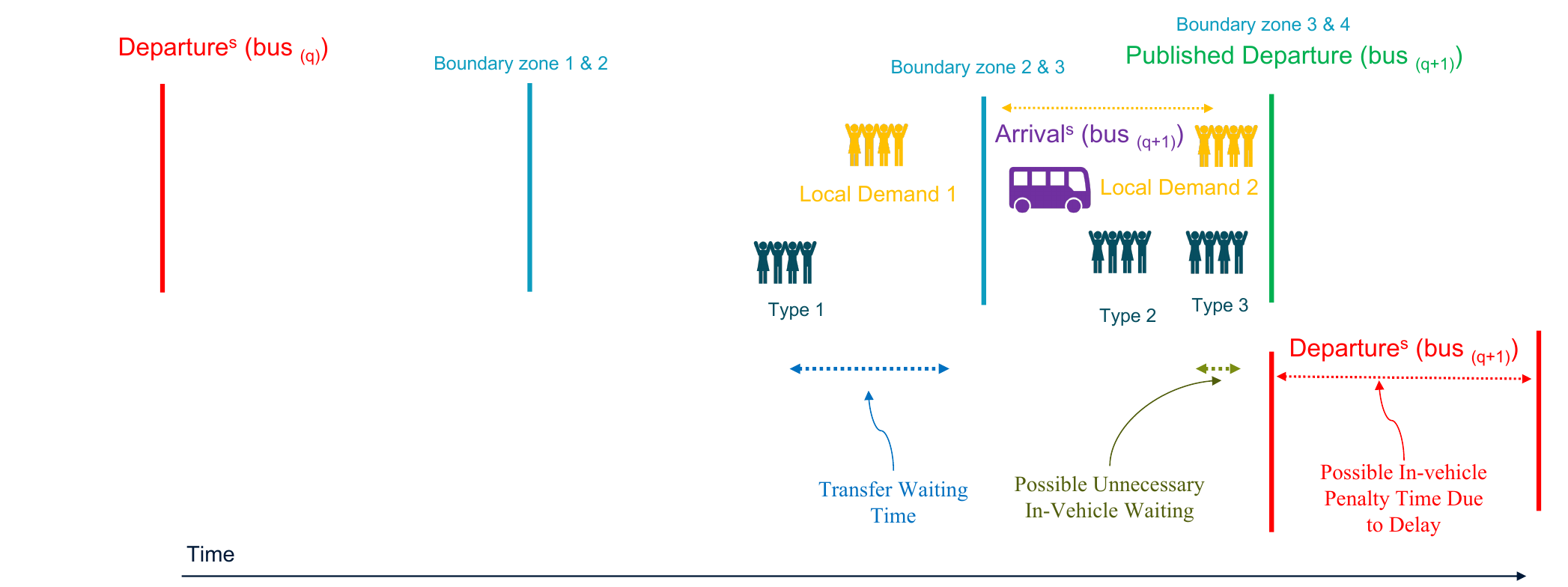}}
\caption{Local passengers' arrival pattern for low-frequency line, Zone 3.}
\label{fig:Low3}
\end{center}
\end{figure}

\begin{figure}[!ht]
\begin{center}
\centerline{\includegraphics[width=0.78\linewidth]{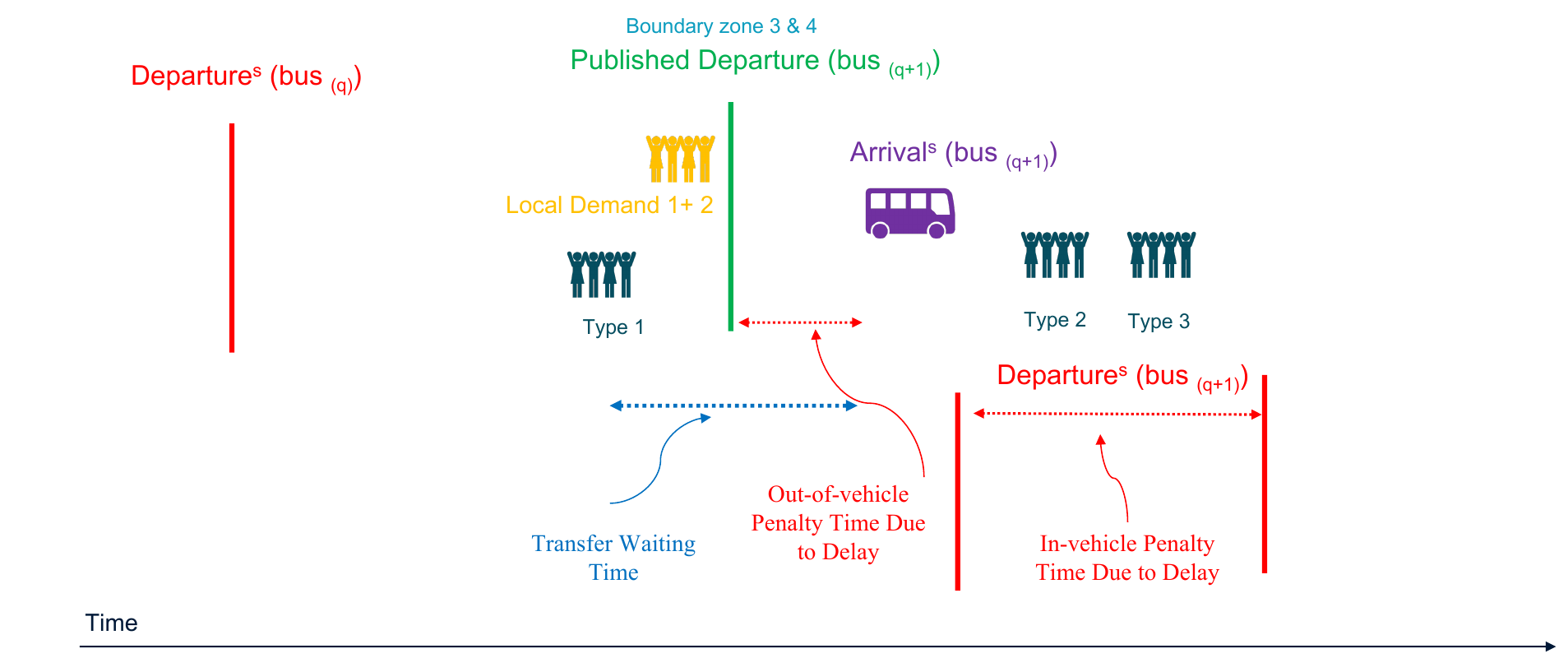}}
\caption{Local passengers' arrival pattern for low-frequency line, Zone 4.}
\label{fig:Low4}
\end{center}
\end{figure}

\paragraph{Low-frequency lines} Local passengers for low-frequency lines tend to check timetables and minimize their waiting times accordingly \citep{ansari2021waiting}. {{Thus, passenger arrivals resemble a skewed distribution towards timetable departure times \citep{ting2005schedule}.}} To simplify the distribution-based arrival pattern, we break down local passenger possible arrivals into four zones. Zone 1 denotes a case when a bus comes early and no local passengers are expected, Figure \ref{fig:Low1}. Zone 2 is when the majority of local passengers are expected to arrive, Figure \ref{fig:Low2}. Zone 3 corresponds to the arrival time of the second group of local passengers, who arrive close to the timetable departure time, Figure \ref{fig:Low3}. Zone 4 refers to a case when a bus arrives beyond its timetable departure time, so all the local passengers are at the bus stop, Figure \ref{fig:Low4}.

Regardless of the bus arrival zones, all three types of successful transfers can occur for a connecting bus. However, the possibility of successful transfer Type 3 differs in each arrival zone as will be discussed in the following. Let us explain how successful transfers can occur for a low-frequency connecting line by assuming a bus arriving at Zone 1 where there are no local passengers. If a group of transferring passengers arrives before the arrival of the bus, successful transfer Type 1 can occur for which the boarding time is serv1. During serv1, another group of transferring passengers could arrive, denoting successful transfer Type 2 with a boarding time of serv2. The model would evaluate to determine if the resulting time marker, i.e., bus arrival time plus serv1 and serv2, is currently at the starting boundary of Zone 2 or not. If not, the first transfer buffer time can be assigned, leading to a possible successful transfer Type 3 with a boarding time of serv3. The model would again determine the time marker. If the sum of serv1, serv2, transfer buffer, and serv3 passes the starting boundary of Zone 2, excess dwell time needs to be considered, Es1. Otherwise, the time marker would start again at the starting boundary of Zone 2, where the first group of local passengers would arrive. From now on, only successful transfer Types 2 and 3 can occur. Successful transfer Type 2 would be when a group of transferring passengers arrive during the service time of the first group of local passengers, servL1, plus the first excessive dwell time from Zone 1, Es1, if any. {If successful transfer Types 2 and 3 occur more than one time, we call them semi-Type 2 and 3, shown as SType in Figures \ref{fig:Low1} - \ref{fig:Low3}.}
This pattern would be formulated similarly for the starting boundary of Zones 3 and 4. 
This formulation logic, so far  presented for a bus that has arrived within Zone 1, 
will also be applied to buses arriving 
at the other zones. Specifically, in the case of arrival Zone 2, serv1 includes any possible successful transfer Type 1 and the first group of local passengers. Similarly, excessive dwell time and possible successful transfer Types 2 and 3 would be considered when the time marker of this specific bus reaches the starting boundary of Zone 3. If a bus arrives at Zone 3, the formulation for considering successful transfers is very similar to high-frequency lines. In this case, serv1 includes all the local passengers, and the upper bound of transfer buffer time would be determined based on serv1 and serv2 and the duration between the arrival of the bus and its timetable departure time. Finally, if a bus arrives at Zone 4, all the successful transfer types can occur, but transfer buffer time would be fixed to zero since the bus should depart as soon as possible.

{{It is important to clarify the similarities and differences between our formulation and the slack time variable used in the literature for mitigating stochasticity in timetables. In our model, we do not explicitly define slack time variables. However, by setting timetable departure times to simultaneously minimize transfer waiting time, unnecessary in-vehicle time, and delay, we assess both the benefits and drawbacks of assigning slack times more thoroughly compared to previous studies. 
One of the key roles of slack time variables is to safeguard against missing a connection due to the late arrival of a feeder bus. So when a feeder bus arrives late, the connecting bus has not yet left since the estimated lateness of transferring passengers has been factored in timetables by slack times. Thus, passengers do not miss their connections. However, in terms of drawbacks, when a bus is being held until its timetable departure time without providing any service and any possible successful transfer, the onboard passengers are penalized with no benefit. This issue is more critical for transferring to low-frequency lines where missing a transfer has a higher cost, thus slack time values are often overestimated to reduce the cost of missing transfers while not properly considering the cost of adding excessive in-vehicle times to onboard passengers. 

{{In our model, the timetable departure times are the first-stage decisions, determined by examining their effect on the second-stage decisions (e.g., the arrival and departure times of buses) in different scenarios. The difference between the timetable departure time and the bus arrival time in a scenario, i.e., dwell time budget, influences how successful transfer occurs and whether the bus can depart based on its timetable departure time. Thus, the dwell time budget also affects how the concept of slack times, i.e., adding extra time into timetable departure times, is used in our model. Any duration in which a bus is waiting until its timetable departure time is considered unnecessary in-vehicle time in our model. For instance, as shown in Figure \ref{fig:Low1}, if a low-frequency bus arrives at Zone 1, the unnecessary in-vehicle time can occur at Zones 1, 2, and 3. 

In summary, the bus dwell time, transfer waiting time, unnecessary in-vehicle time, and delay time are determined in our model based on three factors: (1) the relative bus arrival time to its timetable departure time, (2) the arrival time of passengers, either local or transferring, relative to the arrival of the bus and its timetable departure time, and (3) the number of passengers and their required boarding and alighting time. As a result, the model designs timetable departure times to not only mitigate stochasticity and improve successful transferring and reliability of bus departure times but also to minimize different types of waiting times experienced by affected passengers.}} 
Our model can quantify the trade-off between costs and benefits of having slack times for different groups of passengers in more detail compared to previous studies.}}

\section{{Mathematical Formulation: Stochastic MIP Model}}
\label{SMIPS}
In this section, we progressively introduce the notation used in our formulation and explain the specifics of our proposed stochastic MIP model's objective function and constraints. (For completeness, tables \ref{t:notation} and \ref{t:notationP} describe all the sets and parameters, whereas the other notation tables, Tables \ref{t:variablesAll} -
\ref{t:variablesDelayy}, present the time-related and demand-related variables.) Unless otherwise specified, all decision variables are non-negative and continuous.
We present the extensive form for the sample average approximation of our stochastic program built with the scenario set of $S$. A summary of our full model is provided in Table \ref{tab:Model Summary}, which also summarizes the successful transfer types and dwell time determination zones for high- and low-frequency lines defined in Section \ref{subsec:situationsinstochastic}. Table \ref{tab:Model Summary} can be used as a reference guide while reviewing the details of our model. 

\begin{table}[hbt!]
\small
    \caption{Model Summary}
    \label{tab:Model Summary}
    \centering
\resizebox{\textwidth}{!}{  
\begin{tabular}{r p{0.42\textwidth} || r p{0.448\textwidth}}
     \toprule
     \multicolumn{2}{l}{} & \multicolumn{2}{l}{\it Dwell Time Determination Cases} \\
     \cmidrule(lr){3-4}
     \multicolumn{2}{l}{\it Successful Transfers} & {\it High-frequency lines} & {\it Low-frequency lines} \\
     \cmidrule(lr){1-2}
     \cmidrule(lr){3-4}
     \multicolumn{2}{l}{Type 1: Arrive earlier than the bus} & \multicolumn{1}{l}{Zone 1: Bus arrives early} & \multicolumn{1}{l}{Zone 1: Bus arrives early, no locals expected} \\
     \multicolumn{2}{l}{Type 2: Arrive during boarding} & \multicolumn{1}{l}{Zone 2: Bus arrives late} & \multicolumn{1}{l}{Zone 2: Bus arrives early, majority of locals expected} \\
     \multicolumn{2}{l}{Type 3: After Types 1\&2 boarded, thanks to transfer buffer} & \multicolumn{1}{l}{(early/late means before/after} & \multicolumn{1}{l}{Zone 3: Bus arrives early, second group of locals expected} \\
     \multicolumn{2}{l}{} & \multicolumn{1}{l}{\ timetable departure time)} & \multicolumn{1}{l}{Zone 4: Bus arrives late} \\
     \midrule
     \midrule
     \multicolumn{4}{c}{Objective: Minimize \eqref{Eq1a}+\eqref{Eq1b}+\eqref{Eq1c}+\eqref{Eq1d}} \\
     \multicolumn{4}{c}{\eqref{Eq1a}: Transfer waiting time, \eqref{Eq1b}: In-vehicle waiting time, \eqref{Eq1c}-\eqref{Eq1d}: Delay penalty} \\
     \midrule
     \multicolumn{4}{c}{Timetable departure times: $\Pdep$} \\
     \multicolumn{4}{c}{\eqref{Eq2}-\eqref{Eq4}: Link $\Pdep$ to headway ranges} \\
     \midrule
     \multicolumn{4}{c}{Actual arrival/departure times: $\Aarr$/$\Adepq$} \\
     \multicolumn{4}{c}{\eqref{Eq5}: Link $\Adepq$ and $\Pdep$ for terminals} \\
     \multicolumn{4}{c}{\eqref{Eq6}: Link $\Aarr$ and $\Adepq$, via consecutive transfer nodes} \\
     \midrule
     \midrule
     \multicolumn{2}{c||}{\bf High-frequency Lines} & \multicolumn{2}{c}{\bf Low-frequency Lines} \\
     \cmidrule(lr){1-2}
     \cmidrule(lr){3-4}
     \multicolumn{2}{c||}{\it Determine successful transfer types} & \eqref{Eq83body}-\eqref{Eq90body}: & Specify zone based on bus arrival w.r.t. published \\
     \multicolumn{2}{c||}{\it Assign connecting buses} & & departure time and local arrivals (via ST) \\
     \multicolumn{2}{c||}{\it Calculate wait times} & \eqref{Eq91}: & Y1=0 $\Leftrightarrow$ Types 1\&2 cannot occur \\
     \cmidrule{1-2}
     \eqref{Eq7}: & $Y=0 \Leftrightarrow $ Types 1\&2 cannot occur & \eqref{Eq93}-\eqref{Eq97}: & Calculate \# Type 1 \& first group of Type 2  \\
     \eqref{Eq8}: & $YY=1 \Leftrightarrow $ Type 3 occurs & & (GBD1, GBD2) and their service times (serv1, serv2) \\
     \eqref{Eq9}: & Bound dtb by serv1 (based on Type 2 definition) & \eqref{Eq98}: & YY=1 $\Leftrightarrow$ Type 3 occurs \\
     \eqref{Eq11}-\eqref{Eq12}: & Find first eligible connecting bus for Types 1\&2 (via ZY) & \eqref{Eq100}-\eqref{Eq133}: & Calculate \# Type 3 (GBD3) and their service time (serv3) considering bus arrival at different zones \\
     \eqref{Eq13}-\eqref{Eq15}: & Type 3 $\Rightarrow$ not Type 1, nor Type 2 & \eqref{Eq134}-\eqref{Eq138}: & Determine excess time passing Zone 2 boundary \\
      & (SY=0 $\Rightarrow$ Type 3 or via next-horizon bus) &\eqref{Eq140}-\eqref{Eq149}: & Determine \# second group of Type 2 (GBD2q) and \\
      & (TY=1 if Type 1 of Type 2) & & their service time (serv2q)    \\
      \eqref{Eq16}: & Each feeder bus assigned to a connecting bus, if  & \eqref{Eq150}: & Calculate first group of locals service time (servL1) \\
      & needed from the next horizon   & \eqref{Eq151}-\eqref{Eq169}: & Calculate \# second group of Type 3 (GBD3q) and \\
      \eqref{Eq17}-\eqref{Eq18}: & Calculate waiting time, NTwait  (zero for Types & & their service time (serv3q) \\
      & 2\&3, large value for next-horizon connection) & \eqref{Eq170}-\eqref{Eq186}: & Determine excess time passing Zone 3 boundary \\
      \eqref{Eq19}-\eqref{Eq22}: & Type 1 or 2 indicated (TY = 1) $\Rightarrow$ decide which & \eqref{Eq187}-\eqref{Eq201}: & Calculate \# last group of Type 2 (GBD2g) and their \\
      & one (calculating Twait, comparing to NTwait)  & & service time (serv2g) \\
      & (iwait = 1 $\Rightarrow$ Type 1) & \eqref{Eq202}-\eqref{Eq227}: & Calculate \# last group of Type 3 (GBD3g) and their \\
      & (iwait = 0 $\Rightarrow$ Type 2) & & service time (serv3g) \\
      \cmidrule{1-2}
      \multicolumn{2}{c||}{\it Calculate demand and service times}  & \eqref{Eq230}-\eqref{Eq235}: & Pick a connecting bus for each successful transfer \\
      \cmidrule{1-2}
      \eqref{Eq23}-\eqref{Eq27}: & Calculate \# Types 1,2,3 (GBD1, GBD2, GBD4) & \eqref{Eq236}-\eqref{Eq254}: & Determine excess time passing Zone 4 boundary \\
      \eqref{Eq28}: & Calculate their service times (serv1, serv2, serv4) & \eqref{Eq255}-\eqref{Eq267}: & Determine actual departure times \\
      \cmidrule{3-4}
      \eqref{Eq30}-\eqref{Eq31}: & Bus arrives late $\Leftrightarrow$ pa = 0, Lateness = $\Delta$ & \multicolumn{2}{c}{\bf Unnecessary In-vehicle Time and Delay Penalty for All Lines} \\
      \cmidrule{3-4}
      \eqref{Eq35}-\eqref{Eq36}: & Calculate \# second group of locals (GBD3) and & \eqref{Eq271}-\eqref{Eq272}: &  Calculate total \# boarding \\
      & their service time (serv3) & \eqref{Eq282}: & Determine \# in-vehicle when a bus leaves a node \\
      \cmidrule{1-2}
      \multicolumn{2}{c||}{\it Determine transfer buffer range} & \eqref{Eq283}: & Calculate unnecessary in-vehicle time \\
      \cmidrule{1-2}
      \eqref{Eq37}-\eqref{Eq41}: & Determine stopping time budget = dtw & \eqref{Eq336}-\eqref{Eq346}: &  Dwell time used in \eqref{Eq283} for low-frequency lines \\
      & = max$\{$Pdep-Aarr, alighting time$\}$ & \eqref{Eq400}-\eqref{Eq402}: & Compare unnecessary in-vehicle time with the given   \\
      \eqref{Eq43}-\eqref{Eq49}: & Determine transfer buffer opportunity, set it as  & & threshold to penalize excess \\
      & upper bound on transfer buffer time & \multicolumn{2}{c}{\it If bus departs later than timetabled time, extra penalty for certain passengers} \\
      \eqref{Eq50}: & Set transfer buffer lower bound (negative) & \eqref{Eq403}-\eqref{Eq418}: & Case 1: Bus arrives early but leaves late \\
      \eqref{Eq40}: & Record positive transfer buffer time, if any & \eqref{Eq419}-\eqref{Eq427}: & Case 2: Bus arrives late and departs late \\
      \cmidrule{1-2}
      \multicolumn{2}{c||}{\it Calculate actual departure times} \\
      \cmidrule{1-2}
      \eqref{Eq51}-\eqref{Eq58}: & Determine initial dwell time to see when bus is ready to depart (pi=0 $\Rightarrow$ after timetable departure)\\
      \eqref{Eq63}-\eqref{Eq67}: & pi=0 $\Rightarrow$ Calculate additional locals (GBD5) and their service time (serv5) \\
      \eqref{Eq60}-\eqref{Eq62}: & Determine actual departure times  \\
     \bottomrule
\end{tabular}
}
\end{table}

We use a few notational standards to aid comprehension and readability of our model. We choose upright letters for decision variables and italic letters for sets, indices, and parameters. In addition, in transfer-related components, $l$ and $p$ are indices associated to a feeder bus, whereas $q$ and $l'$ are indices relating to a connecting bus. Nonetheless, in general cases, such as bus arrival and departure times, $l$ and $p$ are indices for both feeder and connecting buses. If only the $q$ and $l'$ indices are used, it only applies to connecting buses.
Furthermore, certain low-frequency line variables have two numbers in their names. The first number denotes the zone for which the variables are used, and the second number indicates the zone in which the bus would arrive. For example, if we refer to the transfer buffer time upper limit established in Zone 2, but for a bus arriving in Zone 1, the variable is tbO2z1. If, on the other hand, the variable's name contains only one number, the definition should be reviewed. For example, tb1 displays the first transfer buffer that might be allocated to a bus regardless of  its arrival Zone. ST1 is another example, which refers to the arrival of a low-frequency line at Zone 1. {{Note that if any variable contains an $s$ index or if any constraint contains a variable with an $s$ index, it means that we are explicitly considering scenario $s$.}}

Our model's objective function, given in \eqref{Eq1}, aims to minimize expected total cost, where each scenario cost contributes proportional to its occurrence probability, $\prob$. It includes three main components for various types of passenger waiting times. The first component, \eqref{Eq1a}, refers to transfer waiting time and the associated passengers. The second component, \eqref{Eq1b}, accounts for the penalty of unnecessary in-vehicle time for passengers who are already on board as a result of transfer buffer times or the bus being held until its timetable departure time. The third and forth components, \eqref{Eq1c} and \eqref{Eq1d}, consider the delay penalty of passengers when the bus is departing after its timetable departure time.   
\begin{subequations}
\label{Eq1}
    \begin{alignat}{2}
    \text{min} \ & \sum_{s \in S} \sum_{(n,l,l') \in \tp}{ \hspace*{0cm}} \sum_{p \in \pl} \sum_{q \in \plpqs}{(\prob)}{(\ctw)}{(\tds)}{(\ndwaitss)} +   & {\footnotesize \color{gray} // \text{transfer waiting} \ } \label{Eq1a}\\
    & \sum_{s \in S} \sum_{(n,l) \in \tpqq}\hspace*{0cm}\sum_{p \in {\pl}}{(\prob)}{(\cvt/2)}{(\vtdo) (\RDiffp)} + & {\footnotesize \color{gray} // \text{in-vehicle waiting} \ } \label{Eq1b}\\
    & \sum_{s \in S} \sum_{(n,l) \in \tpqq}\hspace*{0cm}\sum_{p \in {{\pl}}}{(\prob)}{(\cvh)}{ (\TotalBDEwp + \GBDEwp) (\TEwaitp)} + & {\footnotesize \color{gray} // \text{delay (local \& Type 1)} \ }  \label{Eq1c}\\
    & \sum_{s \in S} \sum_{(n,l) \in \tpqq}\hspace*{0cm}\sum_{p \in {{\pl}}}{(\prob)}{((\cvh + \cvt)/2)}{(\pddop)}{(\TEwaitp)}  & \hspace*{-0.48cm} {\footnotesize \color{gray} // \text{delay (all through \& boarding) }}  \label{Eq1d}
    \end{alignat} 
\end{subequations}
\eqref{Eq1a} multiplies the given transfer demand, $td$, with the corresponding transfer waiting time, NTwait. In \eqref{Eq1b}, vtd refers to the number of onboard passengers, which is multiplied by their unnecessary in-vehicle time, RDiff. We multiply RDiff by half to use the average duration because it is hard to exactly determine the number of passengers and their actual unnecessary in-vehicle time. So, alternatively, we assume the total onboard passengers being penalized for the half of total unnecessary duration of in-vehicle waiting time. Further details are provided when vtd and RDiff are explained.
However, the delay penalty can be determined more accurately.
To consider the duration of the delay time for each group of penalized passengers which depends on the bus arrival zone, different variables need to be defined; the common TEwait is the excessive wait time. \eqref{Eq1c} represents the delay penalty when a bus arrives after its timetable departure time. Thus, local passengers, TotalBDEw, and successfully transferred passengers Type 1, GBDEw, who arrive before the arrival of their connecting bus might be penalized. The Type 1 successful transfer passengers need to be penalized if they arrive before the arrival of the bus and also before the timetable departure time. In this case, the extra waiting time after the timetable departure time needs to be penalized with a different weight compared to normal transfer waiting time. The other potential successful transfers, however, should not be penalized. In fact, they intended to get the subsequent bus as their connection, but because this bus arrived late, they were able to board without having to wait for the transfer. 
The other case, \eqref{Eq1d}, is when a bus arrives before its timetable departure time yet leaves after its timetable departure time. In this instance, all through and boarding passengers, pdd, would be penalized. 




We consider the penalty if the unnecessary in-vehicle time or the delay time exceeds their predetermined threshold. Also, each type of waiting time is additionally weighted with an appropriate cost coefficient, $c$ terms, which can be derived from the literature on passengers' perceptions of different waiting times. {{For \eqref{Eq1d} we use the average weight of unnecessary in-vehicle time and delay time. The reason is that, delay time penalties mostly refer to the case when passengers experience the delay waiting out of the vehicle. However, \eqref{Eq1d} is when passengers are delayed in the vehicle.}}

As previously noted, the purpose of employing two-stage stochastic programming is to provide first-stage scheduling decisions to be fixed in the second (operational) stage, in our case the timetable departure times, that minimize the expected operational cost of the system, i.e., different waiting times in our model, across all possible scenarios. Thus, the timetable departure times variables, Pdep, do not have the scenario index, whereas the actual departure times, Adep, do. 
\eqref{Eq2} ensures that the departure time of each line's first trip from its terminal is less than the line's maximum allowable headway. 
\eqref{Eq3} and \eqref{Eq4}, respectively, adjust the subsequent trips of a line to depart sequentially between the minimum and maximum allowable headway, both from terminals and the transfer nodes. Assuming 
no uncertainty at terminals, buses would depart on time based on the timetable schedules. 
\eqref{Eq5} implies that the departure times under each scenario for all trips equal their respective timetable departure times from terminals. However, depending on the running times between the terminal and the first transfer node as well as those between consecutive transfer nodes in each scenario, the arrival times of the same trips in different scenarios would not be identical, as specified by 
\eqref{Eq6}. 
\begin{alignat}{2}
    & \hspace*{-0.52cm}\Pdep \leq \hlmax  &&\  \hspace*{0.50cm}l \in L , n = n_t(l), p = 1 \label{Eq2}\\
    & \hspace*{-0.52cm}\Pdep - \Pdeppn \geq \hlmin &&\ \hspace*{0.50cm}  l\in L, n \in TN_{(lt)}, p \in \pl, p\neq 1 \label{Eq3}\\  
    & \hspace*{-0.52cm}\Pdep - \Pdeppn \leq \hlmax &&\ \hspace*{0.50cm} l\in L, n \in TN_{(lt)}, p \in \pl, p\neq 1 \label{Eq4}\\
    &\hspace*{-0.52cm}\Adepq = \Pdep &&\ \hspace*{0.50cm} l\in L, n\in \nl, n = n_t(l), p \in \pl, s \in S \label{Eq5}\\
    &\hspace*{-0.52cm}\Aarr = \Adepp + \ttds &&\ \hspace*{0.50cm} l\in L, n\in \nl, p \in \pl, s \in S \label{Eq6}
\end{alignat}
In the following, we formulate the tedious model components, for each scenario $s$,   separately for high- and low-frequency lines. We note that many binary variables are necessarily defined to indicate different cases of successful transfer and dwell time determination based on passenger and bus arrivals, as explained in Section \ref{subsec:situationsinstochastic}, and in turn to  construct logical constraints. 

\subsection{High-frequency Connecting Lines}
\label{subsec:highfreq}
\eqref{Eq7} is defined to check the possible occurrence of successful transfer Type 1. If transferring passengers of feeder bus $p$ of line $l$ arrive at their connecting bus stop later than connecting bus $q$ of line $l'$ plus its dwell transfer buffer time of dtb, then successful transfer Types 1 and 2 cannot occur, i.e., Y = 0.
It should be noted that 
\eqref{Eq7} is only imposed to identify the eligibility of connecting buses. In other words, if Y equals one, it does not necessarily imply that Type 1 or Type 2 transfer was successful. 
\eqref{Eq8} on the other hand, is defined to verify if a successful transfer Type 3 is made with the aid of a transfer buffer, tb. If YY equals one in 
\eqref{Eq8}, it means that successful transfer Type 3 would occur for transferring passengers of bus $p$ of line $l$ transferring to bus $q$ of line $l'$. When YY equals one, the arrival time of transferring passengers equals the arrival time of connecting bus $q$ of line $l'$ plus serv1, serv2, and the specified transfer buffer time. Notice that serv1 refers to the service time for local passengers who arrive before bus $q$ and successfully transfer passengers Type 1. Similarly, serv2 refers to the service time for successful transfer Type 2. 
\eqref{Eq28} determines these  service times. Note that based on the specification for the successful transfer Type 2, the upper limit for dwell transfer buffer, dtb, is serv1, again in 
\eqref{Eq9}.

\medskip
For all $(n, l,l')\in \tph, p \in \pl, q \in \ql$:
\vspace*{-0.25cm}
\begin{alignat}{2}
  &\hspace*{-0.7cm} (\Aarr + \aws) - (\Aarrq + \tbs) \leq 0 \Leftrightarrow \yys = 1 
  \label{Eq7}\\
  & \hspace*{-0.7cm}(\Aarr + \aws) = (\Aarrq + \servos + \servts +\ttbos ) \Leftrightarrow \yyyys = 1 
  \label{Eq8} \\
  & \hspace*{-2.05cm} \text{For all } (n,l') \in \tpqqh, q \in \ql: \nonumber \\
  & \hspace*{-0.7cm} \tbs \leq \servos 
  \label{Eq9}
\end{alignat}


\eqref{Eq11} and \eqref{Eq12} find the first eligible connecting bus for a successful transfer Type 1 and 2 through variable ZY. 
\eqref{Eq13}-\eqref{Eq15} ensure not to assign successful transfer Type 1 or 2 if the correct type of successful transfer is in fact Type 3. It is worth noting that the second summation in 
\eqref{Eq13} eliminates the connecting lines' last trip. This is because if all ZY variables are zero while the successful transfer occurs through the last connecting trip with successful transfer Type 3, then SY would be negative. In 
\eqref{Eq13}, if SY equals zero, two possible situations can appear: (1) both summations equal zero, indicating that all ZY and YY variables are zero. In this case, the successful transfer can take place by the last connecting bus via YY, Type 3, or via the first connecting bus of the next planning horizon; (2) the two summations equal one, and the successful transfer can only happen via Type 3. This is because the first summation equal to one only reveals that bus $(q+1)$ is the first eligible connecting bus for successful transfer Types 1 and 2, but, in fact,  the successful transfer would occur with the aid of bus $q$'s transfer buffer. As a result, if SY equals zero, TY in 
\eqref{Eq14} and \eqref{Eq15} is also zero, indicating that successful transfer Types 1 and 2 are not possible. If SY equals one, TY equals one as well, suggesting that successful transfer would occur through Types 1 or 2. 
\eqref{Eq16} guarantees that each feeder bus is assigned to a connecting bus and has a successful transfer.  Passengers' transfer waiting time is determined by constraints \eqref{Eq17}. Positive transfer waiting time can occur only if transferring passengers arrive at the bus stop before their connecting bus arrives. For successful transfer Types 2 and 3, transfer waiting time equals zero. 
\eqref{Eq18} are defined as if necessary, to assign a connecting bus from the next planning horizon. If no feasible connection exists during the planning horizon to satisfy 
\eqref{Eq16}, the bus $q^*$ from the next horizon is selected, 
\eqref{Eq18} will assign a long waiting time for the first bus with from the next horizon.

\medskip
For all $(n, l,l')\in \tph, p \in \pl$:
\vspace*{-0.25cm}
\begin{alignat}{2}
  &\hspace*{-0.28cm} \zys = \yys  && 
  q = 1
  \label{Eq11} \\
  &\hspace*{-0.28cm} \zys = \yys - \yyns  &&
  q \in \ql \setminus \{1\}
  \label{Eq12} \\
  &\hspace*{-0.28cm} \sys = \sum_{q \in \ql}{(\zys)}-\sum_{q \in \plpn}{(\yyyys)} &&
  \label{Eq13}\\
  &\hspace*{-0.28cm} \tys \leq \sys  &&
  q \in \ql
  \label{Eq14} \\
  &\hspace*{-0.28cm} \tys \geq (\sys - 1) + \zy  &&
  q \in \ql
  \label{Eq15} \\
  & \hspace*{-0.28cm}\sum_{q \in \ql}{(\tys + \yyyys) + \ttys} = 1 &&
  \label{Eq16}\\
  &\hspace*{-0.28cm} \tys = 1 \Rightarrow \Aarrq - (\Aarr + \aws) \leq \ndwaitss 
  \qquad &&  
  q \in \ql
  \label{Eq17} \\
  &\hspace*{-0.28cm}
  \ttys = 1 \Rightarrow \ndwaitsss \geq longwait
  \label{Eq18} 
\end{alignat}
To decide precisely when successful transfers 1 and 2 would occur, we define 
\eqref{Eq19}-\eqref{Eq22}. As explained previously, TY cannot identify which of Types 1 or 2 would occur. Yet, this information is essential so that we can count the number of passengers and their required boarding time for each successful transfer. 
\eqref{Eq19} 
calculates another form of waiting time to detect successful transfer Types 1 and 2. If Twait = Ntwait + dtb, the transfer was accomplished via Type 1. Then, 
\eqref{Eq22} compares NTwait with Twait. If iwait = 1, the difference between NTwait and Twait is less than or equal to dtb, indicating a successful transfer Type 1 would occur. If, otherwise, iwait equals 0, it shows that successful transfer Type 2 would occur.

\medskip
For all $(n, l,l')\in \tph, p \in \pl, q \in \ql$:
\vspace*{-0.25cm}
\begin{alignat}{2}
  &\hspace*{-0.7cm} \tys = 1 \Rightarrow \twaits = (\Aarrq + \tbs) - (\Aarr + \aws) 
  \label{Eq19}\\
  & \hspace*{-0.7cm} (\twaits - \ndwaitss) \leq \tbs  \Leftrightarrow \iwaits = 1  
  \label{Eq22}
\end{alignat}
So far, we have identified the three possible successful transfers for transferring to a high-frequency line. Now, we can calculate the corresponding demand and service times. 
The 
number of local passengers and passengers successfully transferred via Type 1, GBD1, is calculated in . 
\eqref{Eq23}.
The second set of passengers is successfully transferred through Type 2, i.e., TY = 1 and iwait = 0, whose number GBD2 is determined in  
\eqref{Eq25}. 
Similarly, the number of passengers successfully transferred via Type 3, GBD4, 
is calculated by 
\eqref{Eq27}.
Multiplying these numbers with the required boarding time per passenger, \eqref{Eq28} determines the associated service times.

\medskip
For all $(n,l') \in \tpqqfh, q \in \ql$:
\vspace*{-0.25cm}
\begin{alignat}{2}
    & \hspace*{-2cm} \bdos = \sum_{l\in L: (n,l,l') \in \tph}\sum_{p \in \pl}{\iwaits \tds} + (\Aarrq-\Adepq)\lams  
   \label{Eq23}\\
   & \hspace*{-2cm} \bdts = \sum_{l\in L: (n,l,l') \in \tph}\sum_{p \in \pl}{(\tys - \iwaits) \tds}  
   \label{Eq25} \\
  & \hspace*{-2cm} \bdfs = \sum_{l\in L: (n,l,l') \in \tph}\sum_{p \in \pl}{\yyyys \tds}  
  \label{Eq27}\\
  & \hspace*{-2cm} \servos = b^t  \bdos, \ \ \ \servts = b^t  \bdts, \ \ \ \servfs = b^t  \bdfs  
  \label{Eq28}
\end{alignat}
GBD1 and GBD2 do not depend on when the bus arrives with respect to its timetable departure time. However, the relative arrival of the bus to its timetable departure time influences the upper limit for transfer buffer time and hence the potential of successful transfer Type 3. Furthermore, for local passengers, we count the number of passengers who come to a bus stop between the previous bus departure time and the following bus arrival time. Nonetheless, because the bus should not leave before its timetable departure time, another group of local passengers may arrive at the bus stop. Therefore, we need to consider their required service time as well. To determine any probable second group of local passengers and the upper bound for transfer buffer time, we must first verify if a bus arrives before or after its published departure time. 
\eqref{Eq30}-\eqref{Eq31} decide if a bus arrives before its published departure time, pa = 0, and, if so, how much time passes in between, $\Delta$. The number of the second group of local passengers, GBD3, and their service time, serv3, are then calculated using 
\eqref{Eq35} and \eqref{Eq36}.

\medskip
For all $(n,l') \in \tpqqfh, q \in \ql$:
\vspace*{-0.25cm}
\begin{alignat}{2}
  & \hspace*{-1.3cm} \Pdepq - \Aarrq \leq 0 \Leftrightarrow \arpi = 1 
  \label{Eq30}\\
  & \hspace*{-1.3cm} \deltahh = \Pdepq - \Aarrq \Leftrightarrow \arpi = 0 
  \label{Eq31}\\
  & \hspace*{-1.3cm} \bdths = \sum_{l\in L: (n,l,l') \in \tph}\sum_{p \in \pl}{\deltahh \lams } 
  \label{Eq35}\\
  & \hspace*{-1.3cm} \servths = b^t  \bdths  
  \label{Eq36}
\end{alignat}
To determine the upper limit for transfer buffer, we define the {stopping time budget} for a bus. The stopping time budget is the maximum length of two time intervals: (1) the time interval between the bus arrival and its timetable departure time, and (2) the time that the bus has to allocate for alighting passengers. The stopping time budget, dwb, is determined through 
\eqref{Eq37}-\eqref{Eq41}. Then, we need to check if dwb is enough to cover the required time for boarding passengers, excluding the successfully transferred passengers via Type 3. We are excluding GBD3 because we want to check if there is enough time left for any transfer buffer time. 
Nevertheless, if any group of passengers successfully transfers via Type 3 with a transfer buffer time less than or equal to zero, their boarding time is also taken into account for total bus dwell time. The reason is that our model assumes no boarding denial is permitted. In other words, if a passenger comes at a bus stop, the driver must offer them service time. 
\eqref{Eq43}-\eqref{Eq47} define transfer buffer opportunity, tbO, which is set as the maximum limit for transfer buffer time, tb, in \eqref{Eq49}. 
Note that, if passengers arrive during serv2, tb can get a negative value. So, based on the specification for the successful transfer Type 3, the lower limit for transfer buffer, tb, is negative of serv2, defined by 
\eqref{Eq50}. 
Then, 
\eqref{Eq40} determines the positive transfer buffer time which is used for the bus dwell time determination.

\medskip
For all $(n,l') \in \tpqqh, q \in \ql$:
\vspace*{-0.25cm}
\begin{alignat}{2}
  & \hspace*{-0.7cm} (\aldnns) a^t - \deltahh \leq 0 \Leftrightarrow \alti = 1 
  \label{Eq37}\\
  & \hspace*{-0.7cm} \alti = 0 \Rightarrow \dwb = (\aldnns) a^t
  \label{Eq39}\\
  & \hspace*{-0.7cm} \alti = 1 \Rightarrow \dwb = \deltahh
  \label{Eq41}\\
  & \hspace*{-0.7cm} \dwb - (\servos + \servts + \servths) \leq 0 \Leftrightarrow \dwbi = 1 
   \label{Eq43}\\
  & \hspace*{-0.7cm} \tbop = \dwb - (\servos + \servts + \servths) \Leftrightarrow \dwbi = 1 
   \label{Eq45}\\
  & \hspace*{-0.7cm} \tbop = 0 \Leftrightarrow \dwbi = 0 
  \label{Eq47}\\
  & \hspace*{-0.7cm} \ttbos \leq \tbop 
  \label{Eq49}\\
  & \hspace*{-0.7cm} -\servts \leq  \ttbos 
  \label{Eq50}\\
   & \hspace*{-0.7cm} 
   \pttbos = \max \{0, \ttbos \} 
   \label{Eq40}
\end{alignat}
Now that we have identified all boarding passengers and their service time, we need to determine the bus's departure time and, thus, possible unnecessary in-vehicle time or delay penalty. 
\eqref{Eq51} compares the total boarding plus transfer buffer time and alighting time of the bus. The initial dwell time, dwtI, is then computed using 
\eqref{Eq53} and \eqref{Eq55}. We can use dwtI to determine when the bus is ready to depart and compare it to the timetable departure time through 
\eqref{Eq58}. If pi = 0, indicating that the bus would depart after its timetable departure time, we must add another group of local passengers. This is because in GBD3 we calculate the second group of local passengers presuming the bus will depart based on its published departure time. 
\eqref{Eq63} and \eqref{Eq65} determine GBD5, the third group of local passengers. Then, by adding the corresponding boarding time of GBD5 demand, serv5 in 
\eqref{Eq67}, we can calculate the bus departure time using 
\eqref{Eq60}. Note that if the initial dwell time plus the arrival time of the bus is less than the timetable departure time, the departure time of the bus equals its timetable departure time, ensured by 
\eqref{Eq62}.

\medskip
For all $(n,l') \in \tpqqfh, q \in \ql$:
\vspace*{-0.25cm}
\begin{alignat}{2}
  & \hspace*{-0.7cm} (\servos + \servts + \servths +  \pttbos + \servfs) - (\aldnns) a^t \leq 0 \Leftrightarrow \TIs = 1 
   \label{Eq51}\\
   &\hspace*{-0.7cm} \dwsi = (\servos + \servts + \servths +  \pttbos + \servfs) \Leftrightarrow \TIs = 1 
   \label{Eq53} \\
   & \hspace*{-0.7cm} \dwsi = (\aldnns) a^t  \Leftrightarrow \TIs = 0 
   \label{Eq55}\\
    & \hspace*{-0.7cm} (\Aarr + \dwsi) \leq \Pdep \Leftrightarrow \ipp = 1  
    \label{Eq58}\\
  & \hspace*{-0.7cm} 
  \ipp = 0 \Rightarrow \bdes = \sum_{l\in L: (n,l,l') \in \tph}\sum_{p \in \pl}{(\Aarr + \dwsi - \Pdep) \lams}  
  \label{Eq63}\\
  & \hspace*{-0.7cm} \ipp = 1 \Rightarrow \bdes = 0
  \label{Eq65}\\
   & \hspace*{-0.7cm} \serves = b^t  \bdes  
   \label{Eq67}\\
   & \hspace*{-0.7cm} \ipp = 0 \Rightarrow \Adep = (\Aarr + \dwsi + \serves) 
   \label{Eq60}\\
   & \hspace*{-0.7cm} \ipp = 1 \Rightarrow  \Adep = \Pdep 
   \label{Eq62}
\end{alignat}

Since determining unnecessary in-vehicle time and delay penalty is very similar for high- and low-frequency lines, we will explain them in Section \ref{BothUN}, after discussing successful transferring and local passengers service time for low-frequency lines.

\subsection{Low-frequency Connecting Lines}
As explained in Section \ref{subsubsec:zonedefns}, we establish four zones for the arrival of a low-frequency bus. 
Regarding the local passengers' arrival patterns used for two zones, based on the historical data, we assume that we can convert the distribution of local passengers' arrival at a bus stop into two parts. The first group of passengers would arrive $\betao \hlq$ before the timetable departure time, while the second group would arrive $\betat \hlq$ before the timetable departure time, where $\betat \hlq$ is smaller than $\betao \hlq$. These parameters can be selected independently of headways and based on available data for each line/trip. 
\eqref{Eq83body}-\eqref{Eq90body} are specified to distinguish the four zones. If ST1 equals zero in 
\eqref{Eq83body}, bus $q$ of line $l$ would arrive in Zone 1 when there are no local passengers. Otherwise, if ST1 equals one, it simply indicates that the bus is in Zone 2, 3, or 4. Likewise, if ST2 equals one, the bus would arrive in Zone 3 or 4. Thus, if ST1 equals one and ST2 equals zero, ST12 equals one indicating that the bus would arrive in Zone 2. If ST3 equals zero in 
\eqref{Eq87body} and ST2 equals one, then ST23 equals one in 
\eqref{Eq90body}, so the bus would arrive in Zone 3. Finally, if ST3 equals one, the bus would arrive in Zone 4 after its timetable departure time.

\medskip
For all $(n,l') \in \tpqql, q \in \ql$:
\vspace*{-0.25cm}
\begin{alignat}{2}
    & \hspace*{-0.7cm} (\Pdep - \Aarrq) \leq (\betao \hlq)  \Leftrightarrow \STo = 1 
    \label{Eq83body}\\
    & \hspace*{-0.7cm} (\Pdep - \Aarrq) \leq \betat \hlq \Leftrightarrow \STt = 1 
    \label{Eq85body}\\
    & \hspace*{-0.7cm} \Pdep \leq \Aarrq \Leftrightarrow \STth = 1 
    \label{Eq87body}\\
    & \hspace*{-0.7cm} \STot = \STo - \STt 
    \label{Eq89body}\\
    & \hspace*{-0.7cm} \STtt = \STt - \STth 
    \label{Eq90body}
\end{alignat}



Now that we know the arrival zone of a connecting bus of a low-frequency line, we need to determine the number of passengers while accounting for their relative arrival time in comparison to the arrival time of  bus $q$ on line $l'$ and its timetable departure time. In that regard, we investigate the three different waiting times and the bus departure time. The modeling approach to determining service times and identifying possible successful transfers for low-frequency lines is very similar to the approach we described above for high-frequency lines. The only differences are how we handle local passengers' dwell time, dwell transfer buffer upper limit and the transfer buffer upper limit, which vary depending on which zone a bus arrives in. As such, we provide the full model components for low-frequency lines in \ref{appsec:lowfreq}, and present a summary for them in Table \ref{tab:Model Summary}. 

As shown in Table \ref{tab:Model Summary}, we sequentially determine the number of passengers for each bus, starting from the earliest bus arrival time to the latest that fall into Zone 1 to Zone 4. More specifically, we calculate in order the number of Type 1 successful transfers (GBD1), first group of Type 2 successful transfers (GBD2), first group of Type 3 successful transfers (GBD3), along with their service times, determine the excess time passing Zone 2 boundary, calculate the second group of Type 2 successful transfers (GBD2q), whose service time combined the service time of first group of local passengers determine the number of second group of Type 3 successful transfers (GBD3q), accordingly determine the excess time passing Zone 3 boundary, calculate the number of last group of Type 2 successful transfers (GBD2g) as well as the number of last group of Type 3 successful transfers (GBD3g), accordingly determine the excess time passing Zone 4 boundary, and finally calculate the actual bus departure times. 

\subsection{Unnecessary In-vehicle Time and Delay Penalty for Both High- and Low-frequency Lines}
\label{BothUN}
To determine the number of passengers experiencing unnecessary in-vehicle time or delay penalty, we must track the number of passengers in each bus.
The number of total boarding can be determined for high- and low-frequency lines through 
\eqref{Eq271} and \eqref{Eq272}, respectively.
\eqref{Eq282} determines the number of passengers on board when a bus leaves a node. As shown in 
\eqref{Eq282}, the in-vehicle demand of a bus when it arrives at a node equals the number of in-vehicle demands when it departs its previous node plus the net number of boarding and alighting at the regular stops in between, $\Ivddp - \ssps$. Note that there is no capacity limitation and all the generated demand can board their bus. The number of in-vehicle passengers when the bus departs its terminal is given inputs. 
\begin{alignat}{2}
& \hspace*{-2.5cm} \text{For all } (n,l')  \in \tpqqh, q \in \ql: \nonumber \\
    & \hspace*{-0.7cm} \Tbdns =  \bdos + \bdts + \bdfs + \bdths + \bdes 
    \label{Eq271}\\
& \hspace*{-2.5cm} \text{For all } (n,l') \in \tpqql, q \in \ql : \nonumber \\ 
    & \hspace*{-0.7cm} \Tbdns =  \bdos + \bdts + \bdfs + \bdths + \bdq + \nonumber\\
    & \hspace*{0.78cm} \bdthsq + \bdg + \bdthsg 
    \label{Eq272}\\
& \hspace*{-2.5cm} \text{For all } (n,l') \in \tpqq, q \in \ql : \nonumber \\ 
    & \hspace*{-0.7cm} \Ivdd = \Ivddp - \ssps - \aldnns + \Tbdns 
    \label{Eq282}
\end{alignat}


Regardless of on-time or delayed departure time of a bus, it is possible that onboard passengers experience unnecessary in-vehicle time. The unnecessary in-vehicle time can occur several times for a bus based on its frequency. {{Note that if passengers are alighting during no boarding, there is a service occurring, and no penalty is needed. This is because unnecessary in-vehicle time is when onboard passengers experience they are being held with no service, either boarding or alighting.}} For a low- or high-frequency line, during transfer buffer time(s) and when the bus is ready to leave but has to wait until its timetable departure time, the onboard passengers experience unnecessary in-vehicle time if there is no alighting service.  Thus, 
the difference between the bus arrival time plus dwtI, the larger value of total boarding time and alighting time, and its departure time, is the unnecessary in-vehicle time determined through 
\eqref{Eq283}. 
\begin{alignat}{2}
  & \hspace*{-1.2cm} \text{For all } (n,l') \in \tpqq, q \in \ql: \nonumber \\ 
    & \hspace*{-0.4cm} \RDiff = \Adep - (\Aarr + \dwsi ) 
    \label{Eq283}
\end{alignat}
For a high-frequency bus, the initial dwell time, dwtI, is already computed using 
\eqref{Eq53} and \eqref{Eq55}. On the other hand, for a low-frequency bus, unnecessary in-vehicle time can occur in different zones. For instance, when the bus arrives early, at Zone 1, the through passengers would experience unnecessary in-vehicle time during different transfer buffer times of each zone or/and when there is no demand, i.e., no service of alighting or boarding, but the bus has to wait until its timetable departure time. For low-frequency, the total service time depends on the arrival zone of the bus which is determined through 
\eqref{Eq336}-\eqref{Eq343} for each zone separately. Then, the unnecessary in-vehicle time based on both boarding and alighting times is determined through 
\eqref{Eq344}-\eqref{Eq346}. 
Note that it is possible that the bus departs on-time or later than its timetable departure time, but unnecessary in-vehicle time still can occur, specially for low-frequency lines. 
\begin{alignat}{2}
  & \hspace*{-0.67cm} \text{For all } (n,l') \in \tpqqfl, q \in \ql: \nonumber \\ 
    & \hspace*{-0.4cm} \STo = 0 \Rightarrow \dwIzo = \servos + \servts + \servths + \servLo + \servLt +   &&\ \nonumber\\ & \hspace*{4cm} \servtsq + \servthsq + \servtsg +  \servthsg  
    \label{Eq336} \\
    & \hspace*{-0.4cm} \STo = 1 \Rightarrow \dwIzo = 0    
    \label{Eq337} \\
    & \hspace*{-0.4cm} \STot = 1 \Rightarrow \dwIzt = \servos + \servts + \servths +  \servLt + \servtsg + \servthsg 
    &&
    \label{Eq338} \\
    & \hspace*{-0.4cm} \STot = 0 \Rightarrow \dwIzt = 0   
    \label{Eq339} \\   
    & \hspace*{-0.4cm} \STtt = 1 \Rightarrow \dwIzth  = \servos + \servts + \servths    
   \label{Eq342} \\
   & \hspace*{-0.4cm} \STtt = 0  \Rightarrow \dwIzth = 0  
   \label{Eq343} \\
   & \hspace*{-0.4cm} (\dwIzo + \dwIzt + \dwIzth ) - (\aldnns) a^t \leq 0 \Leftrightarrow \TIs = 1 
   \label{Eq344} \\
   &\hspace*{-0.4cm} \TIs = 1 \Rightarrow \dwsi = (\dwIzo + \dwIzt + \dwIzth)   
   \label{Eq345} \\
   & \hspace*{-0.4cm} \TIs = 0 \Rightarrow \dwsi = (\aldnns) a^t   
   \label{Eq346}
\end{alignat}
After determining the unnecessary in-vehicle time, we compare it with a given threshold in 
\eqref{Eq400}, and if is exceeds the predetermined limit, the onboard passengers would be penalized in the objective function. Since it is not easy to consider the exact duration of unnecessary in-vehicle time and the associated demand, alternatively, we determine the total unnecessary in-vehicle time and assume, on average, the onboard demand would experience this penalty as half of the duration. 
The number of passengers who would be penalized in case of unnecessary in-vehicle time includes through passengers and boarding passengers, either local or transferring passengers, as shown in 
\eqref{Eq401}-\eqref{Eq402}.

\medskip
For all $(n,l') \in \tpqq, q \in \ql$:
\vspace*{-0.25cm}
\begin{alignat}{2}
    & \hspace*{-0.7cm} \Rths - \RDiff \leq 0 \Leftrightarrow \Rthsi = 1 
    \label{Eq400}\\
    & \hspace*{-0.7cm} \Rthsi = 1 \Rightarrow \vtdo = \Ivddp - \ssps -\aldnns + \Tbdns  
    \label{Eq401}\\
    & \hspace*{-0.7cm} \Rthsi = 0 \Rightarrow \vtdo = 0   
    \label{Eq402}
\end{alignat}
On the other hand, if a bus departs later than its timetable departure time, we determine how much later the bus would depart. A bus would depart after its timetable departure time regardless of whether it arrives before the timetable departure time or after. This distinction is essential because passengers perceive the delay, i.e., departing later than the timetable departure time, differently depending on whether they are inside or outside the bus. Moreover, this distinction guides the model to which passengers should be penalized in the objective function due to a delayed departure time. Thus, we discuss the two possible cases separately: (1) a bus arrives after its timetable departure time, and (2) the bus arrives before its timetable departure time yet departs later than its timetable departure time.

For the first case, the local passengers and the transferring passengers who arrive before the arrival of their connecting bus, i.e., successful transfer Type 1, might be the passengers who need to be considered for the out-of-vehicle cost of a delayed bus. The delay, i.e., the difference between the departure time and the timetable departure time of a high-frequency line, is determined through 
\eqref{Eq403}. Before determining the number of penalized passengers we have to check weather the delay exceeds the predetermined threshold in 
\eqref{Eq404}. If yes, the considered delay is defined in 
\eqref{Eq405} and \eqref{Eq4051}. 
The number of local passengers to be considered for the delay penalty is specified by 
\eqref{Eq406} and \eqref{Eq4061}. Note that this is not the total number of local passengers boarding the bus. In fact, the local passengers who arrive after the timetable departure time are supposed to get to the next bus, but since this bus is late, they will get on the bus, and their boarding time is considered as well. However, we do not include them as penalized passengers. So the number of local passengers is determined based on the difference between the timetable departure time of the current bus and the departure time of the last bus. Then we need to examine the arrival of successful transfer Type 1. If NTwait of a group of passengers is more than zero, it means that they arrive before the arrival of their connecting bus. In order to consider them for out-of-vehicle delay penalty, we compare the transfer waiting time and the delay duration and the threshold, constraints \eqref{Eq407}; if the NTwait is smaller, we know that although they arrive before their bus, they do not arrive early enough to be considered as penalized passengers. Otherwise, we have to consider them for the delay penalization, TYE = 0. Note that these constraints are valid for both high- and low-frequency lines. To make sure that we are considering the case of late departure time, i.e., ap = 1, STTYE is defined. If STYYE = 1, it means that the connecting bus arrives after its timetable departure time, ap = 1, and the successfully transferred passengers Type 1, TYE = 0, need to be penalized for the delay, ensured by 
\eqref{Eq408}. The associated demand is determined in 
\eqref{Eq409}. 
\begin{alignat}{2}
& \hspace*{-2cm} \text{For all } (n,l') \in \tpqqfh, q \in \ql: \nonumber \\ 
    & \hspace*{0cm} \Ewait = \Aarr - \Pdep \Leftrightarrow \arpi = 1 
    \label{Eq403}\\
    & \hspace*{0cm} \EwaitII= 1 \Rightarrow \TotalBDEw = (\Pdep - \Adep) \lams  
   \label{Eq406}\\
    & \hspace*{0cm} \EwaitII= 0 \Rightarrow \TotalBDEw = 0  
    \label{Eq4061}\\
    & \hspace*{0cm} \STTYE = \api - \TYE  
    \label{Eq408}\\
    & \hspace*{0cm} \STTYE= 1 \Rightarrow \GBDEw = \bdos - (\Aarrq-\Adepqp)\lams  
   \label{Eq409} \\
& \hspace*{-2cm} \text{For all } (n,l') \in \tpqql, q \in \ql: \nonumber \\ 
    & \hspace*{0cm} \Ewait -  \Aths \geq 0 \Leftrightarrow \EwaitII= 1 
    \label{Eq404}\\
    & \hspace*{0cm} \EwaitII= 1 \Rightarrow \TEwait = \Ewait  
    \label{Eq405}\\
    & \hspace*{0cm} \EwaitII= 0 \Rightarrow \TEwait = 0  
    \label{Eq4051}\\
    & \hspace*{0cm} \ndwaitss - \Ewait - \Aths \leq 0 \Leftrightarrow \TYE = 1 
   \label{Eq407} 
\end{alignat}
With similar logic, the constrains for a low-frequency bus that arrives and departs later than its timetable departure time are as follows. The only difference is that, we use variable ST3 = 1 denoting the late arrival of the bus. Also, the local demand being penalized for the late arrival of the bus is a given value for low-frequency lines.

\medskip
For all $(n,l') \in \tpqql, q \in \ql$:
\vspace*{-0.25cm}
\begin{alignat}{2}
    & \hspace*{-0.7cm} \STth = 1 \Rightarrow \Ewait = \Aarr - \Pdep  
    \label{Eq410}\\
    & \hspace*{-0.7cm}  \STth = 0 \Rightarrow \Ewait = 0
    \label{Eq411}\\
    & \hspace*{-0.7cm} \Ewait - \Aths \geq 0 \Leftrightarrow \EwaitII= 1 
    \label{Eq412}\\
    & \hspace*{-0.7cm} \EwaitII= 1 \Rightarrow \TEwait = \Ewait  
    \label{Eq413}\\
    & \hspace*{-0.7cm} \EwaitII= 1 \Rightarrow \TotalBDEw = D   
    \label{Eq414}\\
    & \hspace*{-0.7cm} \EwaitII= 0 \Rightarrow \TotalBDEw = 0  
    \label{Eq415}\\
    & \hspace*{-0.7cm} \STTYE = \STth - \TYE  
    \label{Eq416}\\
    & \hspace*{-0.7cm} \STTYE= 1 \Rightarrow \GBDEw = \bdos - D  
    \label{Eq417}\\
    & \hspace*{-0.7cm} \STTYE= 0 \Rightarrow \GBDEw = 0  
    \label{Eq418}
\end{alignat}
The second case is when the bus arrives before its timetable departure time yet departs with delay. To distinguish this case, we define two binary variables, pipa for high-frequency lines and Esthii for low-frequency lines. If a high-frequency bus arrives before its timetable departure time, pa = 0, and it departs after its timetable departure time, pi = 0, resulting in pipa = 0, in 
\eqref{Eq419}. If pipa = 0, then the delay is determined through 
\eqref{Eq421} and compared with the threshold in 
\eqref{Eq423}. Similarly, for a low-frequency line, if it departs after its timetable departure time, i.e., pi = 0, and arrives before its timetable departure time, i.e., ST3 = 0, it means that the bus would depart later than its timetable departure time, Esthii = 0, given 
\eqref{Eq420}.
The delay is then determined and compared to the threshold in 
\eqref{Eq425a} and \eqref{Eq425b}. For both high- and low-frequency lines, if AthsI = 1, we have to consider the passengers being penalized for the delay with in-vehicle cost. In this case, all the in-vehicle passengers are being penalized, via 
\eqref{Eq426} and \eqref{Eq427} for high- and low-frequency lines, respectively. Note that the extra demand for high-frequency lines, arriving during the delay time after the timetable departure time is not considered for delay penalization as those passengers were in fact supposed to get to the next bus. 
\begin{alignat}{2}
& \hspace*{-2.6cm} \text{For all } (n,l') \in \tpqqfh, q \in \ql: \nonumber \\ 
    & \hspace*{-1.57cm} \pipa = \ipp + \arpi 
    \label{Eq419}\\
    & \hspace*{-1.57cm} \pipa = 0 \Rightarrow \ADiff = \Adep - \Pdep   
    \label{Eq421}\\
    & \hspace*{-1.57cm} \pipa = 1 \Rightarrow \ADiff = 0  
    \label{Eq423}\\
    & \hspace*{-1.57cm} \Aths - \ADiff \leq 0 \Leftrightarrow  \Athsi = 1 
    \label{Eq425a}\\
    & \hspace*{-1.57cm} \Athsi = 1 \Rightarrow \pddo = \Ivddp -  \ssps  -  \aldnns + \Tbdns - \bdes    
    \label{Eq426}\\
& \hspace*{-2.6cm} \text{For all } (n,l') \in \tpqql, q \in \ql: \nonumber \\   
    & \hspace*{-1.57cm} \Esthii = \ipp + \STth 
    \label{Eq420}\\
    & \hspace*{-1.57cm} \Esthii = 0 \Rightarrow \ADiff = \Adep - \Pdep 
    \label{Eq422}\\
    & \hspace*{-1.57cm} \Esthii = 1 \Rightarrow \ADiff = 0  
    \label{Eq424}\\
    & \hspace*{-1.57cm} \Aths - \ADiff \leq 0 \Leftrightarrow  \Athsi = 1 
    \label{Eq425b}\\
    & \hspace*{-1.57cm} \Athsi = 1 \Rightarrow \pddo = \Ivddp -  \ssps  - \aldnns + \Tbdns  
   \label{Eq427}
\end{alignat}

\section{Solution Method}
\label{sec:SMethodologyS}
As stochastic programming models are notoriously difficult to solve to optimality, approximation methods are widely used. One of the most commonly applied approximation methods is Sample Average Approximation (SAA), which we also use in this work. In Section~\ref{SMIPS}, we presented the extensive form for the SAA of our problem, including scenario copies of recourse decisions. 
The procedure of creating the scenario set should be such that the SAA model closely approximates the original stochastic program. 
Thus, one of the main concerns regarding stochastic programming applications is generating a sufficient number of scenarios that appropriately mimic the stochasticity of random variables and their interactions without significantly increasing the computational burden or, even worse, leading to the intractability of models. 
For this issue, a group of studies investigated possible approaches to cluster scenarios to reduce the SAA method's computational cost \citep{bertsimas2022optimization,fairbrother2019problem,henrion2018problem,keutchayan2021problem,narumproblem}.
Different studies have shown the accuracy of \emph{problem-driven scenario reduction} as a practical approach to imitating the stochasticity of the system with a considerably small number of scenarios. Hence, we also implement a problem-driven scenario reduction technique in this study to substantially reduce the number of scenarios.
Finally, \emph{leveraging the loosely coupled structure} in our model, we propose to use the \emph{progressive hedging algorithm} to derive high-quality first-stage solutions to our SAA model. 
This algorithm decomposes the model into scenario-based subproblems enabling us to reach near-optimal solutions more efficiently compared to directly solving it via a (commercial) MIP solver. 
Nevertheless, as a final step, we use the solution of the PH algorithm as a \emph{warm-start solution} in the MIP solver GUROBI to solve our model to optimality, or at least attempt to find improving solutions as well as to derive an optimality gap. 
In what follows, we explain the details of 
our solution methodology.

\subsection{Problem-driven Scenario Reduction}
The two primary groups of scenario reduction approaches are distribution-based and problem-driven. We apply the latter based on the MIP clustering model proposed in~\citep{keutchayan2021problem} to provide a smaller scenario set that reflects the larger sample of possible scenarios for our stochastic model. The merit of the problem-driven approach over the distribution-based one is the utilization of specific problem-based knowledge to generate scenarios, especially in the case of complex models.
Moreover, in recent studies, problem-driven methods have yielded promising results because they take into account both the probability distribution of stochastic variables in a system and how these distributions influence the behaviour of the model.

Before applying the clustering MIP model, we need to provide some input data. First, we need the optimal solution of the stochastic model when it is solved deterministically for only one scenario. Next, for the $N = |S|$ number of scenarios in the initial sample, we evaluate the remaining scenarios using the optimal solution of one scenario, in our case, timetable departure times. In other words, $V_{i,j}$ shows the value of the objective function of the model being solved deterministically for scenario $j$ while fixing the timetable departure times based on the solution of the deterministic model of scenario $i$. So, if we evaluate the $N$ scenarios of the initial sample, we have an $N\times N$ matrix. Subsequently, based on these values, the following clustering MIP model is solved. 
\begin{alignat}{2}
    \text{min} \ & \frac{1}{|S|} \sum_{i \in S} {t_i}  \label{EqC1}\\
    \text{s.t.} \ & t_j \geq \sum_{i \in S} {x_{ij}}{V_{j,i}} - \sum_{i \in S} {x_{ij}}{V_{j,j}}
    &&\  \hspace*{1.25cm} j \in S  \label{EqC2}\\
    & t_j \geq \sum_{i \in S} {x_{ij}}{V_{j,j}} - \sum_{i \in S} {x_{ij}}{V_{j,i}}
    &&\  \hspace*{1.25cm} j \in S  \label{EqC3}\\
    & x_{ij} \leq u_j
    &&\  \hspace*{1.25cm} (i,j) \in S \times S  \label{EqC4}\\
    &  x_{jj} = u_j
    &&\  \hspace*{1.25cm} j \in S   \label{EqC5}\\
    &   \sum_{j \in S} {u_{j}} = M
    &&\ \hspace*{1.25cm}  \label{EqC7}\\
    &   \sum_{j \in S} {x_{ij}} = 1
    &&\ \hspace*{1.25cm} i \in S  \label{EqC6}\\
    &  x \in \{0,1\}^{|S| \times |S|}; u \in \{0,1\}^{|S|}; t \in \mathbb{R}^{|S|} &&\ \hspace*{3.25cm}  \label{EqC8}
\end{alignat}

The MIP model groups the scenarios into $M\ll N$ clusters, each containing one representative scenario through minimizing implementation errors. The implementation error occurs by implementing the solution of the problem with the reduced number of scenarios in the original problem being solved for all scenarios. 
The objective function of the clustering model minimizes $t_j$ variables that measure the clustering error. More specifically, the clusters are selected so as to minimize the difference between the representative scenario of each cluster and the other cluster members. 
\eqref{EqC2}-\eqref{EqC3} establish the absolute value of this the clustering error. If scenario $i$ belongs to the same cluster as the representative scenario $j$, then $x_{ij} = 1$. 
\eqref{EqC4} and \eqref{EqC5} ensure that cluster $j$ exists to be able to assign scenario $i$ to cluster $j$. If scenario $j$ is one of the representative scenarios, the variable $u_j$ equals one. 
\eqref{EqC7} ensures that the model clusters all the scenarios into $M$ clusters. In addition, 
\eqref{EqC6} is set so that each scenario is only assigned to one cluster. {{Note that this MIP clustering method necessitates that the stochastic program has relatively complete recourse; that is, for every feasible first-stage solution, it is guaranteed to have feasible second-stage solutions, which is satisfied in our problem.}} 

As previously stated, the scenario reduction MIP model requires deterministically solving the original stochastic model for each scenario in the initial sample set. However, in our case, even solving the original stochastic model for a single scenario is computationally expensive. As a result, during this step of our solution approach, we use a simpler model. We disregard the arrival pattern of local passengers, which converts our model to the buffer-based deterministic model from  \citep{ansarilari4043348novel}. Using that model, we can obtain the $V_{i,j}$ values considerably faster. {{In Section \ref{sec:Experiment}, we demonstrate the reliability of developing the reduced set of scenarios using the simplified version of our original stochastic model through numerical examples.}} 

\subsection{Progressive Hedging Algorithm}

So far we have found the reduced number of scenarios to include in our SAA model. However, another factor that causes the complexity of solving the proposed stochastic model is the number of variables, especially integer variables. Thus, even with reduced number of scenarios we need to apply another technique to solve the model efficiently. To do so, in this study, we use the Progressive Hedging (PH) method to break down the model down into smaller scenario-based subproblems.
The PH method is designed based on the augmented Lagrangian relaxation concept in which one or more constraints are relaxed and instead their violations are included in the objective function.
To apply the PH method we first add copies of the first-stage decisions for each scenario, $\Pdeppsnk$, along with the \emph{nonanticipativity constraints}: 
\begin{alignat}{2}
    & \hspace*{-0.7cm} \Pdeppsnk - \Pdep = 0 &&\ \hspace*{4.25cm} l\in L, n\in \nl, p \in \pl, s \in S \label{EqC9}
\end{alignat}
Clearly, constraints \eqref{EqC9} are the only reason that prevents us from breaking down the model to scenario-based sub-problems. Thus,
in the PH method, we relax \eqref{EqC9} and add two violation costs to the objective function, \eqref{Eq1Pe} and \eqref{Eq1Pf} for each scenario. The \eqref{Eq1Pe} refers to the Lagrangian multipliers, $\mu$, of the nonanticipativity constraints \eqref{EqC9} and the second element, \eqref{Eq1Pf}, refers to the quadratic penalty with a positive penalty parameter of $\rho$. Then, the relaxed model decomposes by scenario, yielding the subproblem given in \eqref{Eq1P} for each to be used at iteration $k$ of the algorithm.
\begin{subequations}
\label{Eq1P}
    \begin{align}
    \text{min} \ & \sum_{(n,l,l') \in \tp}{ \hspace*{0cm}} \sum_{p \in \pl} \sum_{q \in \plpqs}{(\prob)}{(\ctw)}{(\tds)}{(\ndwaitss)} + \label{Eq1Pa}\\
    & \sum_{(n,l) \in \tpqq}\hspace*{0cm}\sum_{p \in {\pl}}{(\prob)}{(\cvt/2)}{(\vtdo) (\RDiffp)} + \label{Eq1Pb}\\
    & \sum_{(n,l) \in \tpqq}\hspace*{0cm}\sum_{p \in {{\pl}}}{(\prob)}{(\cvh)}{ (\TotalBDEwp + \GBDEwp) (\TEwaitp)} +  \label{Eq1Pc}\\
    & \sum_{(n,l) \in \tpqq}\hspace*{0cm}\sum_{p \in {{\pl}}}{(\prob)}{((\cvh + \cvt)/2)}{(\pddop)}{(\TEwaitp)} + \label{Eq1Pd}\\
    & \sum_{(n,l) \in \tpqq}\hspace*{0.15cm}\sum_{p \in {\pl}}{\muup(\Pdepps - \Pdeppsa )}   +\label{Eq1Pe}\\
    & \sum_{(n,l) \in \tpqq}\hspace*{0.15cm}\sum_{p \in {\pl}}{||(\Pdepps - \Pdeppsa )||^2_2}\hspace*{0.2cm}{{\rho}/2 } \label{Eq1Pf}\\
    \text{s.t.} \  
    & \ \Pdepps \leq \hlmax  && \hspace{-3cm} l \in L , n = n_t(l), p = 1 \\
    & \ \Pdepps - \text{Pdep}^{ksn}_{l(p-1)} \geq \hlmin &&  \hspace{-3cm} l\in L, n \in TN_{(lt)}, p \in \pl, p\neq 1 \\  
    & \ \Pdepps - \text{Pdep}^{ksn}_{l(p-1)} \leq \hlmax && \hspace{-3cm} l\in L, n \in TN_{(lt)}, p \in \pl, p\neq 1 \\
    & \ \Adepq = \Pdepps && \hspace{-3cm} l\in L, n\in \nl, n = n_t(l), p \in \pl \\
    \ & \ \eqref{Eq6}-\eqref{Eq427}, \ \eqref{Eq91}-\eqref{Eq267} 
    \end{align} 
\end{subequations}
At iteration $k$, each scenario-based subproblem is solved 
to decide on $\Pdepps$ variables in a way that the expected discrepancy between $\Pdepps$ and $\Pdeppsa$ across all scenarios is minimized, where $\Pdeppsa$ is the average of $\Pdepps$ over all scenarios from the previous iteration. The steps of the PH method are presented in Algorithm \ref{alg:cap}.
\begin{algorithm}[!ht]
\caption{Progressive Hedging Algorithm}\label{alg:cap}
\small

\begin{algorithmic}
 \smallskip
\State Initialize  \hspace*{0.45cm} Define the termination criteria: relative convergence tolerance of $\theta$ or $k_{max}$\\
  \hspace*{1.8cm} iteration counter $k \leftarrow 0, \rho \leftarrow 1, \muup \leftarrow 0$, \hspace*{0.05cm} $n \in \nl$, $l \in L, p \in P_{l}$, $s \in S$
\State Step 1: 
 \hspace*{0.70cm} Solve the subproblem \eqref{Eq1P} for each scenario $s \in S$, saving optimal $\Pdepps$ values
\State Step 2:  \hspace*{0.70cm} Compute the initial solution of $\Pdeppsakz = \sum_{s \in S}{\prob \Pdeppskz}$\\
\hspace*{1.785cm} Compute initial $\muupz = \rho (\Pdeppskz - \Pdeppsakz)$, for all $s \in S$\\
\textbf{Repeat}\\
\State Step 3:  \hspace*{0.70cm} Solve  \eqref{Eq1P} for all $s \in S$, and get \Pdeppskk \\
\State Step 4:  \hspace*{0.70cm} Update  $\Pdeppsakk = \sum_{s \in S}{\prob \Pdeppskk}$
\State Step 5:  \hspace*{0.70cm} Update $\muupkk = \muup +  \rho (\Pdeppskk - \Pdeppsakk)$, for all $s \in S$\\
\State Step 6: \hspace*{0.70cm}  $k$ = $k$ + 1\\
\\
\textbf{Until} \hspace*{0.90cm} $ \sum_{s \in S}{\prob} {||\Pdeppskk - \Pdeppsakk||_2} \leq \theta$ or $k \geq k_{max}$
\end{algorithmic}

\end{algorithm}

Note that it is not guaranteed that the final solution derived from PH, with minimum or zero differences between $\Pdepps$ and $\Pdeppsa$, is the optimal solution of our SAA model. Therefore, in our study, after some preliminary examination, we define a threshold to stop the iterations of the PH algorithm and then use the best achieved $\Pdeppsa$ solution as the \emph{warm start} to solve our SAA model to optimality using the commercial solver GUROBI.

\section{Numerical Results}
\label{sec:Experiment}
In this section, we present various experiments to illustrate the primary contributions of our stochastic transfer synchronization model, namely, the detailed formulation of different waiting times and bus dwell time determination based on passengers' arrival patterns while considering stochasticity in passenger walking times between bus stops, bus running times, and passenger demand. Moreover, we discuss the efficiency and performance of our solution method. Prior to examining the findings of our experiments, we describe the input data and the characteristics for chosen transfer nodes as our instances in Section \ref{input}, followed by the explanation of how we generate sample scenarios and the reduced subset in Section \ref{scengen}.

In Section \ref{ExpTwo2}, for the numerical experiment, we first solve our stochastic model, model \textit{SM}, and compare the results with its deterministic counterpart, model \texttt{DSM}. We assess the results through one of the most common measures for stochastic models, the Value of Stochastic Solution (VSS). VSS values demonstrate the goodness of a stochastic model over a deterministic model in which random variables are substituted with their mean or another estimated value, such as 80\% percentile. VSS is the difference between objective function values when evaluating a test scenario with solutions of the deterministic model, being solved with the average of random variables, and the solutions of the stochastic model. Typically, the VSS value is reported as the mean of all test scenarios.

Subsequently, in order to demonstrate the merit of a sophisticated stochastic transfer synchronization model over simpler models, we also reformulate the deterministic model, model \texttt{DB}, from  \citep{ansarilari4043348novel} into a stochastic one, model \texttt{SDB}.
We compare the final results of different waiting costs between the stochastic model suggested in this work and the \texttt{SDB} model with regard to passenger arrival patterns and more accurate depictions of different passenger waiting times, i.e., unnecessary and delay waiting times. This comparison aims to highlight the value of our approach of introducing stochasticity into timetabling and transfer synchronization models.

We also discuss the computational performance of our solution method in Section \ref{ExpPH}.
All the experiments are implemented using Python and solved by GUROBI \citep{gurobi}, on a MacBook Pro with a 2.8 GHz (Intel Core i7) processor and 16 GB memory.

\subsection{Input Data for Experiment Instances and Parameters}
\label{input}

The proposed stochastic model requires several input data on both passenger demand and bus service operation on a trip-by-trip basis. Despite emerging technologies for capturing transit data, the available sources, such as APC, AFC and AVL, are insufficient to provide input for this study. Therefore, we derived the required data sets for our experiments from an agent-based simulation platform called Nexus that mimics real-world transit operations. This mesoscopic simulation tool, Nexus, was built at the University of Toronto \citep{srikukenthiran2017enabling}. General Transit Feed Specification (GTFS) data and travel behaviour survey data are inputs to this simulation. Similar to previous studies, Nexus considers bus dwell times \citep{rashidi2022estimating,wen2017data} and bus running times between two consecutive stops \citep{hollander2008estimation,kieu2014establishing,li2017quantifying, uno2009using} as random variables with lognormal distributions. Furthermore, Nexus provides dynamic integration of multiple simulation software for transit demand assignment, each responsible for different parts of a transit system, e.g., surface transit, subway and rail. Thus, transit users are agents moving between different software components as travelling between various origin-destination choices while their overall travel time is minimized based on a route-choice decision-making process. 
As such, Nexus can produce highly disaggregate output, including bus running times and dwell times as well as incoming bus in-vehicle demand and the number of boarding and alighting passengers.

{{As a proof of concept}}, we examine our proposed stochastic model and solution method for two instances with one transfer node. The two transfer nodes are selected from the City of Toronto Transit Network (CTTN) (see Figure~\ref{fig:Toronto} for visualization of the considered nodes). We use the simulation model of the City of Toronto, based on 2016 data, for the morning peak hours of 6 to 9 a.m. The characteristics of each transfer node, mainly the intersecting lines' headway combinations and the number of transfer pair directions, are presented in Table \ref{table:NodeInfo}. 

\renewcommand{\arraystretch}{1.25}
\begin{table}[!ht]
\caption{Characteristics of chosen transfer nodes from the CTTN}
\label{table:NodeInfo}
\centering
\small
\resizebox{0.59\textwidth}{!}{

\begin{tabular}{l c c c } \toprule
\multicolumn{1}{c}{Transfer} & Number of & \multicolumn{1}{c}{Min - Max - Average} & Number of Transfer \\
\multicolumn{1}{c}{Node ID} & Routes & \multicolumn{1}{c}{Headway (min)} & Directions \\
\cmidrule(r){1-1}
\cmidrule(lr){2-2}
\cmidrule(lr){3-3}
\cmidrule(lr){4-4}
\cmidrule(lr){4-4}
{A} & {6} & {6.1 - 17.0 - 10.8} & {7}\\
{B} & {5} & {4.6 -  11.0 - 6.5} & {10}\\
\bottomrule
\end{tabular}


}
\end{table}

It is necessary to use appropriate coefficients for different waiting time components in the objective function, as these coefficients substantially affect the credibility of solutions. Hence, we employ the various waiting time weights derived from the research on passengers' time perceptions \citep{cats2014optimizing,cats2016dynamic,gavriilidou2019reconciling}: in-vehicle unnecessary waiting time (\texttt{$c^{vt}$}$= 1.5$), transfer waiting time (\texttt{$c^{tw}$}$ = 2$), out-of-vehicle delay (\texttt{$c^{sl}$}$ = 3.27 $), {{in-vehicle delay (\texttt{$0.5(c^{st} + c^{sl})$}$ = 2.39$).}}
Additionally, for three additional parameters, time for opening and closing doors, alighting time and boarding time per passenger, we use the estimated values for CTTN, 7.43, 1.12 and 1.96 seconds per person, respectively \citep{miller2018surface}. For simplicity, the following parameters are defined based on line headways:
$\betao = 30\% \hl$, $\betatt = 9\% \hl$, $\hlmin = 90\% \hl$, and $\hlmax = 110\% \hl$ except for $\Aths = \Rths = 1$ minute and
except for {{$\nudemand = 80\%$}}. Note that all these parameters can specifically be chosen for each line based on historical data.

\subsection{Scenario Generation and Reduction}
\label{scengen}
As the output of a simulation run, Nexus data sets include fixed values mainly suitable for deterministic experiments and analysis. 
One approach to providing data for our stochastic model is to simulate the network for different days with different random seed numbers. However, to generate the required scenarios, we adopt the standard approach applied in the stochastic programming literature, Monte Carlo simulation, to reduce computational effort. 
We derive the associated lognormal distributions for feeder and connecting bus running times data from Nexus, from which we sample the scenarios. However, for passenger walking times between bus stops and demand data, due to the complexity of using fitted distributions for each trip at each node, we randomly draw from a range of possible values for each data input. This approach has been used by \citet{cao2019optimal} for incorporating dwell time variability as uncertainty noises in their stochastic transfer synchronization model. For local passengers of high-frequency lines, we assume the Poisson process random arrival with given arrival rates from Nexus. 

We generate a sample scenario set with $N = 100$ scenarios. As recommended by \citet{keutchayan2021problem} to cluster the initial sample set to 2-3\% of the total number of scenarios, the scenario reduction MIP model selects 
three scenarios as the representatives. 
As previously mentioned, due to the long computational time, we do not use the proposed stochastic model to provide the input values, i.e., $V_{i,j}$ for the scenario reduction MIP model. Instead, we use the stochastic version of the deterministic model from  \citep{ansarilari4043348novel}, model \texttt{SDB}, as the simpler version of model \texttt{SM}. 
The output and calculation time of our scenario-reduction step are presented in Table \ref{table:SceRe}.
To test and validate the solution of our stochastic model, we also generate 500 new scenarios. 

\renewcommand{\arraystretch}{1.25}
\begin{table}[!ht]
\caption{Output and calculation time of the scenario-reduction step}
\label{table:SceRe}
\centering
\small
\resizebox{0.809\textwidth}{!}{


\begin{tabular}{l c c c c} \toprule
\multicolumn{1}{c}{Transfer} & Number of scenarios & \multicolumn{1}{c}{Time for preparing} & Time to solve \\
\multicolumn{1}{c}{node ID} & in each cluster & \multicolumn{1}{c}{{$V_{i,j}$}  \text{values (hour)}} & the scenario-reduction MIP (min) \\
\cmidrule(r){1-1}
\cmidrule(lr){2-2}
\cmidrule(lr){3-3}
\cmidrule(lr){4-4}
\cmidrule(lr){4-4}
{A} & {37, 42, 21} & {6.34} & {6}\\
{B} & {22, 57, 21} & {6.22} & {14}\\
\bottomrule
\end{tabular}


}
\end{table}

\subsection{Comparison of Deterministic and Stochastic Models} 
\label{ExpTwo2}
As mentioned above, in this experiment, we have two main goals: (i) examining the necessity of incorporating details into a transfer synchronization model, and (ii) verifying the value of stochastic transfer synchronization formulation over a deterministic one. To do so, we evaluate the 500 test scenarios individually with four sets of timetable departure times, i.e., solutions of models \texttt{SM}, \texttt{DSM}, \texttt{DB}, and \texttt{SDB}. 

Tables \ref{table:Node745} and \ref{table:Node745B} present the average values over 500 scenarios for the two transfer node instances. We report the values of the objective function components, the total value of the objective function (with weights), as well as the average total delay and unncessary in-vehicle time to help with analysis. Model \texttt{SM} provides timetable departure times which lead to the minimum objective function values and passenger transfer waiting times amongst all the models. 

\renewcommand{\arraystretch}{1.25}
\begin{table}[!ht]
\caption{Model comparison - Average over 500 scenarios - Node A}
\label{table:Node745}
\centering
\small
\resizebox{0.89\textwidth}{!}{\begin{tabular}{l l c c c c} \toprule
Objective item & \texttt{SM} & \texttt{DSM} & \texttt{SDB} & \texttt{DB}\\
\cmidrule(r){1-1}
\cmidrule(lr){2-2}
\cmidrule(lr){3-3}
\cmidrule(lr){4-4}
\cmidrule(lr){5-5}
{Total objective function value (person.min)} &
{38,215} &
{63,681} &
{40,570} &
{70,572}\\
{Total transfer waiting time * demand (person$\cdot$min)} &
{37,209} &
{60,432} &
{39,545} &
{69,285}\\
{Total delay time (min)} &
{29} &
{1} &
{28} &
{47}\\
{Total delay time * demand (person$\cdot$min)} &
{127} &
{1.43} &
{121} &
{238} \\
{Total unnecessary in-vehicle time (min)} &
{179} &
{467} &
{193} &
{218}\\
{Total unnecessary in-vehicle time * demand (person$\cdot$min)} &
{879} &
{3,248} &
{904} &
{1,077}\\
\bottomrule
\end{tabular}

}
\end{table}

\begin{table}[!ht]
\caption{Model comparison - Average over 500 scenarios - Node B}
\label{table:Node745B}
\centering
\small
\resizebox{0.89\textwidth}{!}{
\begin{tabular}{l l c c c c} \toprule
Objective item & \texttt{SM} & \texttt{DSM} & \texttt{SDB} & \texttt{DB}\\
\cmidrule(r){1-1}
\cmidrule(lr){2-2}
\cmidrule(lr){3-3}
\cmidrule(lr){4-4}
\cmidrule(lr){5-5}
{Total objective function value (person.min)} &
{120,383} &
{130,749} &
{133,539} &
{128,164}\\
{Total transfer waiting time * demand (person$\cdot$min)} &
{119,257} &
{129,860} &
{132,448} &
{128,164}\\
{Total delay time (min)} &
{145} &
{186} &
{134} &
{139}\\
{Total delay time * demand (person$\cdot$min)} &
{482} &
{362} &
{470} &
{219}\\
{Total unnecessary in-vehicle time (min)} &
{91} &
{93} &
{91} &
{61}\\
{Total unnecessary in-vehicle time * demand (person$\cdot$min)} &
{644} &
{527} &
{621} &
{309}\\
\bottomrule
\end{tabular}
}
\end{table}

Furthermore, the obtained average values are smaller when \texttt{DSM} is applied instead of \texttt{DB} for transfer node A. Although for transfer node B,  \texttt{DB} performs better compared to  \texttt{SDB} and \texttt{DSM}, the differences are not notable. Nonetheless, the results of transfer node B indicate the effect of transfer node characteristics on outcomes, which has been discussed in  \citep{ansarilari4043348novel}. 
	
{We note that we have established a 1-minute threshold for evaluating the delay penalty and unnecessary in-vehicle time. Consequently, the numbers displayed in Tables \ref{table:Node745} and \ref{table:Node745B} are derived from the records that are above the 1-minute threshold. This is one of the reasons why the transfer waiting time fraction of the total objective function values is significantly larger than the other waiting time fractions across all outcomes. As seen in the results, deterministic models have a tendency to generate timetables with longer unnecessary in-vehicle time. This behaviour is comparable to assigning large values to slack time variables in order to prevent missing connections at the expense of penalizing passengers who are already on board. However, in stochastic models, unnecessary in-vehicle times are not only shorter than in those of the deterministic models but they are also assigned more strategically so that they can be used to reduce transfer waiting times, similar to the concept of transfer buffer times in  \citep{ansarilari4043348novel}.

The relative disparities in unnecessary in-vehicle times and delay times between nodes A and B are another noteworthy result. In node A, unnecessary in-vehicle times are longer than delay times, whereas, in node B, the opposite is true. This is likely owing to the variances between headway combinations, distributions of bus running times, and demand variability at transfer nodes A and B. Moreover, in node A, the average headway is larger due to two low-frequency lines, whereas in node B, it is closer to high-frequency lines. The substantially reduced transfer waiting times in node A demonstrate the necessity and effectiveness of our models in transfer synchronization for low-frequency lines, where transfer synchronization is vital. }

Table \ref{table:formulation} shows that a thorough model is preferable to a simplified one for solving the transfer synchronization problem.  
Regardless of solving the transfer synchronization problem in stochastic or deterministic settings, this comparison highlights that the accuracy of the results would be undermined if the details of the transfer process and timetabling are not correctly incorporated into models.

\renewcommand{\arraystretch}{1.25}
\begin{table}[!ht]
\caption{Analysis of quality of model formulation}
\label{table:formulation}
\centering
\small
\resizebox{0.49\textwidth}{!}{


\begin{tabular}{l c c c} \toprule
\multicolumn{1}{c}{Transfer node ID} & \texttt{SM} over \texttt{SDB} & \texttt{DSM} over \texttt{DB}\\
\cmidrule(r){1-1}
\cmidrule(lr){2-2}
\cmidrule(lr){3-3}
\cmidrule(lr){4-4}
\cmidrule(lr){4-4}
{A} & {5.81\%} & {14.37\%}\\
{B} & {9.85\%} & {-2.02\%}\\
\bottomrule
\end{tabular}


}
\end{table}

To determine the value of the stochastic model over deterministic modelling approach for the transfer synchronization problem, we calculate the VSS value both for the proposed stochastic model \texttt{SM} and the model \texttt{SDB}, as shown in Table \ref{table:VSSValuesSM}. 
The positive VSS values for model \texttt{SM} for both the transfer nodes justify the extra effort and indicate the value of applying the stochastic model \texttt{SM} over its deterministic counterpart, model \texttt{DSM}. 
For example, for transfer node A, the 40\% VSS implies that using the timetable departure times derived from the model \texttt{SM} decreased the average of objective function values of test scenarios by 40\% compared to using the timetable departure times derived from the model \texttt{DSM}.
Note that the model used as the evaluation framework is model \texttt{DSM}, which represents a more realistic testbed compared to model \texttt{DB}. That is why the VSS values for model \texttt{SDB} is negative. 
Also, the higher value of 42\% of model \texttt{SDB} does not indicate that model \texttt{SM} performs worse. {{In fact, comparison of the VSS values of the two stochastic models \texttt{SM} and \texttt{SDB}, Tables \ref{table:formulation} and \ref{table:VSSValuesSM}, demonstrates that, in the case of adopting a stochastic model for the transfer synchronization problem, it is preferable to use a detailed model that more accurately mimics the system's stochasticity.}} As mentioned above, model \texttt{SM} performs the best for all transfer nodes. 

\renewcommand{\arraystretch}{1.25}
\begin{table}[!ht]
\caption{Analysis of quality of stochastic solution}
\label{table:VSSValuesSM}
\centering
\small
\resizebox{0.59\textwidth}{!}{


\begin{tabular}{l c c c c} \toprule
\multicolumn{1}{c}{Transfer node ID} & VSS for model \textit{SM} & VSS for model \textit{SDB}\\
\cmidrule(r){1-1}
\cmidrule(lr){2-2}
\cmidrule(lr){3-3}
\cmidrule(lr){4-4}
{A} & {40\%} & {42\%}\\
{B} & {8\%} & {-4\%}\\
\bottomrule
\end{tabular}


}
\end{table}

\subsection{Computational Performance of the Solution Method}
\label{ExpPH}
The total calculation time for the PH method involves two main parts: solving the subproblems for each representative scenario and applying Algorithm \ref{alg:cap}. Table \ref{table:PHTime} presents the total calculation times for solving the subproblems, applying Algorithm \ref{alg:cap}, and solving the stochastic model. 

\renewcommand{\arraystretch}{1.25}
\begin{table}[!ht]
\caption{Solution Method Performance}
\label{table:PHTime}
\centering
\small
\resizebox{0.89\textwidth}{!}{
\begin{tabular}{l l c} \toprule
Node & Task Item & Time (hr)\\
\cmidrule(r){1-1}
\cmidrule(lr){2-2}
\cmidrule(lr){3-3}
{} &
{Solving the subproblems for each representative scenario} &
{0.90}\\
{A} &
{Total iteration time until the reaching the stopping criterion} &
{1.83}\\
{} &
{Solve the stochastic model with all scenarios until reaching the stopping criterion} &
{3.00}\\
\midrule
\midrule
{} &
{Solving the subproblems for each representative scenario} &
{0.40}\\
{B} &
{Total iteration time until the reaching the stopping criterion} &
{0.26}\\
{}&
{Solve the stochastic model with all scenarios until reaching the stopping criterion} &
{0.01}\\
\bottomrule
\end{tabular}
}
\end{table}

\renewcommand{\arraystretch}{1.25}
\begin{table}[!ht]
\caption{Solution Method Output values}
\label{table:PHTask}
\centering
\small
\resizebox{0.89\textwidth}{!}{
\begin{tabular}{l l c} \toprule
Node & Item & Value\\
\cmidrule(r){1-1}
\cmidrule(lr){2-2}
\cmidrule(lr){3-3}
{} &
{Initial difference of \hspace*{0.03cm}$ \sum_{s \in S}{\prob} {||\Pdeppskk - \Pdeppsakk||_2}$ } &
{3236}\\
{} &
{Final difference of \hspace*{0.03cm} $ \sum_{s \in S}{\prob} {||\Pdeppskk - \Pdeppsakk||_2}$} &
{696}\\
{} &
{Number of iterations} &
{8}\\
{A}&
{Met stopping criterion for Algorithm \ref{alg:cap}} &
{Time}\\
{} &
{Initial gap of solving the stochastic model with warm-start solution from PH} &
{27.4\%}\\
{} &
{Final gap of solving the stochastic model until reaching the stopping criterion} &
{17.1\%}\\
{}&
{Met stopping criterion for solving the stochastic model}&
{Time}\\
\midrule
\midrule
{} &
{Initial difference of \hspace*{0.03cm} $ \sum_{s \in S}{\prob} {||\Pdeppskk - \Pdeppsakk||_2}$ } &
{1502}\\
{} &
{Final difference of \hspace*{0.03cm}$ \sum_{s \in S}{\prob} {||\Pdeppskk - \Pdeppsakk||_2}$} &
{299}\\
{} &
{Number of iterations} &
{7}\\
{B}&
{Met stopping criterion for Algorithm \ref{alg:cap}} &
{$\theta$}\\
{} &
{Initial gap of solving the full model with warm-start solution from PH} &
{6.32\%}\\
{} &
{Final gap of solving the full model until reaching the stopping criterion} &
{0.0\%}\\
{}&
{Met stopping criterion for solving the stochastic model}&
{Final gap}\\
\bottomrule
\end{tabular}
}
\end{table}


The stopping criteria for Algorithm \ref{alg:cap} are the maximum number of iterations, $k_{max}=15$; time limit of 5 hours, and $\theta = 600$ seconds in this study. The first met stopping criterion is reported with other outputs of the solution method in Table \ref{table:PHTask} for each transfer node. The final $\Pdeppsk$ values derived from Algorithm \ref{alg:cap} are then used as a warm start to solve the stochastic model including all the representative scenarios. The initial and final gaps are also presented in Table \ref{table:PHTask}. The stopping criteria for the stochastic model are a final gap of less than $15\%$ and a total calculation time limit of 3 hours. Lastly, we note that the commercial solver GUROBI could not provide any solutions given the same calculation time budget.

\section{Conclusion}
\label{sec:SConclusion}
Most previous studies focused on bus running times variability to provide synchronized timetables. However, stochastic timetabling and transfer synchronization literature has not effectively considered bus dwell time and demand uncertainty. 
In this paper, we propose a sophisticated model that jointly considers the variability of passenger walking times between bus stops, bus running times and bus dwell times, as well as uncertainty in passenger demand and their arrival patterns at bus stops to reduce transfer waiting times, delay times, and unnecessary in-vehicle times. The main output of the proposed two-stage stochastic programming model is timetable departure times which would be published and made available to passengers. 

We estimate the number of boarding passengers, both local and transferring, considering their arrival times at the bus stop relative to the arrival time of their bus and its timetable departure time to calculate dwell times. We consider distribution-based arrival of local passengers for low-frequency lines in which passengers minimize their waiting time based on timetable departure times. On the other hand, local passengers for high-frequency lines are assumed to arrive randomly at bus stops. 
Passengers are subject to experience an unnecessary in-vehicle time when a bus is detained without any service, i.e., being ready to depart but must wait until its timetable departure time. However, if a bus needs a longer dwell time to provide service to passengers and, as a result, exceeds its timetable departure time, or if a bus arrives after its timetable departure time despite having a short dwell time, the bus departure time would be later than its timetable departure time and passengers would be penalized for the delay.

We choose two single transfer nodes from CCTN for our experiments. 
We solve a sample average approximation of our model using the progressive hedging algorithm and a problem-based scenario reduction strategy. 
To evaluate the advantages of a stochastic approach over a deterministic one for the transfer synchronization problem, we analyze the outcomes using the standard measure of VSS. The adoption of stochastic models is justified by the high VSS values for both transfer nodes, 
In addition, by comparing the proposed stochastic model to a simpler model, we assess the need to include passenger arrival pattern and dwell time determination in the transfer synchronization problem and timetabling. In this regard, the deterministic transfer synchronization model provided in \citep{ansarilari4043348novel} is reformulated to a stochastic model. Results reveal that the proposed stochastic transfer synchronization model can decrease passenger waiting times substantially more than deterministic and simple stochastic transfer synchronization models.
Our computational experiment on two single transfer nodes in the City of Toronto serve as a proof-of-concept for the potential benefits of the offered stochastic model. 

For future work, some extensions of our model can be considered to generate potentially higher-quality timetables. For instance, in our model, we consider trips in single directions, ignoring possible effects on layover times at terminals and assuming buses are always available to start their trips in the opposite direction according to the timetable. Consideration of trip cycle times and vehicle assignments might improve the applicability of transfer synchronization models in practice. The objective function can be extended to reflect the total costs incurred by all stakeholders, i.e., all types of transit passengers and operators, or maximize the net benefits for all. 
Our model does not consider the potential impact of the produced timetable on demand levels and patterns. Future studies can explore the joint formulation of travel behaviour and timetabling models in an integrated framework, in particular passengers route/mode change. 
Moreover, real-time transfer synchronization strategies, such as the bus holding, can be incorporated into the stochastic model as the second-stage implementable recourse decisions. Further methodological developments would be of crucial interest for future research to be able to address practical-size instances as well as all of the potential extensions of the stochastic model. Accordingly, studying joint optimization problems, such as jointly designing  route layouts, route frequencies and bus sizes in a transfer synchronization model, would also be beneficial.

\section*{Disclosure statement}
No potential conflict of interest was reported by the authors.

\section*{Data Availability Statement}
Data is available upon request.


\singlespacing
\small
\bibliographystyle{elsarticle-harv} 
\bibliography{buffer_ref}

\clearpage
\newpage
\normalsize
\appendix
\onehalfspacing
\section*{APPENDIX}\label{sec:appendix}
\setcounter{figure}{0}
\setcounter{table}{0}

\section{Notation}

\renewcommand{\arraystretch}{0.995}
\begin{table}[!ht]
\caption{Notation: Sets}
\label{t:notation}
\small
\centering
\resizebox{0.89\textwidth}{!}{\begin{tabular}{ll} 
\toprule
\ns & \hspace*{-0.20cm} Set of chosen transfer nodes of the whole network; $n$ $\in$ \ns  \\
$NT$ & \hspace*{-0.20cm} Set of terminals and chosen transfer nodes of the whole network; $n$ $\in$ NT  \\
$L$ & \hspace*{-0.20cm} Set of lines passing at least one of the chosen transfer nodes in the network; $l$ $\in$ $L$ \\
$L_{Lf}$ & \hspace*{-0.20cm} Set of low-frequency lines passing at least one of the chosen transfer nodes in the network; \\
& \hspace*{-0.20cm} $l$ $\in$ $L_{Lf}$ \\
$L_{Hf}$ & \hspace*{-0.20cm} Set of high-frequency lines passing at least one of the chosen transfer nodes in the network;\\& \hspace*{-0.20cm} $l$ $\in$ $L_{Hf}$ \\
\nl & \hspace*{-0.20cm} Set of transfer nodes along line $l$ $\in$ $L$; $n$ $\in$ \nl\\
$TN_{(lt)}$ & \hspace*{-0.20cm} Set of transfer nodes and starting terminal along line $l$ $\in$ $L$; $n$ $\in TN_{(lt)}$ \\
\tp & \hspace*{-0.20cm} Set of transfer pairs of lines at transfer nodes; $(n,l,l') \in \tp$, $n$ $\in$ \nl \hspace*{0.04cm} $\cap$ \nlq,  $l, l^\prime$  $\in$ $L$, $l \neq l^\prime$\\*[0.1cm]
\tph & \hspace*{-0.20cm} Set of transfer pairs of lines at transfer nodes, when the connection is a high-frequency line;\\
& \hspace*{-0.20cm} $(n,l,l') \in \tph$, $n$ $\in$ \nl \hspace*{0.04cm} $\cap$ \nlq, $l^\prime$ $\in$ $L_{Hf}$ $l$  $\in$ $L$, $l \neq l^\prime$\\*[0.1cm]
\tpl & \hspace*{-0.20cm} Set of transfer pairs of lines at transfer nodes when the connection is a low-frequency line;\\
& \hspace*{-0.20cm} $(n,l,l') \in \tpl$, $n$ $\in$ \nl \hspace*{0.04cm} $\cap$ \nlq,  $l^\prime$ $\in$ $L_{Lf}$ $l$   $\in$ $L$, $l \neq l^\prime$\\*[0.1cm]
$PL^{n}$ & \hspace*{-0.20cm} Set of lines passing a terminal or a transfer node; $l \in PL^{n}$,  $n\in TN_{(lt)}$\\
\tpqq & \hspace*{-0.20cm} Set of connecting lines at a terminal or transfer node $n \in$ \ns; $(n,l') \in$ \tpqq\\
\tpqqfh & \hspace*{-0.20cm} Set of high-frequency connecting lines at a terminal or transfer node $n \in$ \ns; $(n,l') \in$ \tpqqfh\\
\tpqqfl& \hspace*{-0.20cm} Set of low-frequency connecting lines at a terminal or transfer node $n \in$ \ns; $(n,l') \in$ \tpqqfl\\
\pl & \hspace*{-0.20cm}  Index set of feeder buses of line $l \in L$ during the planning horizon; $p \in \pl = \{1,2,\hdots, \floor{T/\hl} \}$ \\*[0.1cm]
\plpq & \hspace*{-0.20cm} Index set of connecting buses of line $l' \in L$ during the planning horizon; \\
& \hspace*{-0.20cm} $ q \in \plpq = \{ 1,2,\hdots, \floor{T/\hlq} \} $\\*[0.1cm]
\plpqn & \hspace*{-0.20cm} Index set of connecting buses of line $l'$ during the planning horizon without the last bus; \\
\plpqs & \hspace*{-0.20cm} Index set of connecting buses of line $l' \in L$ during the planning horizon, \plpq, and the first trip \\
& \hspace*{-0.20cm} of the next horizon; $q \in \plpqs = \plpq \cup q^*_{l'}$  \\*[0.1cm]
$S$ & \hspace*{-0.20cm} Set of scenarios of  bus running times between transfer nodes, passenger walking times between   \\
& \hspace*{-0.20cm}  bus stops, and number of different types of passengers for trip of each line in the network; $s$ $\in$ $S$ \\
\bottomrule
\end{tabular}

}
\end{table}

\begin{table}[!ht]
\caption{Notation: Parameters}
\label{t:notationP}
\small
\centering
\resizebox{0.89\textwidth}{!}{\begin{tabular}{ll} 
\toprule
$T$ & \hspace*{-0.40cm} Length of the planning time horizon\\ 
\hl & \hspace*{-0.40cm} Headway of line $l$ during the planning horizon; $l$ $\in$ $L$ \\*[0.05cm]
\hlmin & \hspace*{-0.40cm} Minimum allowable headway of line $l$ during the planning horizon; $l$ $\in$ $L$ \\*[0.05cm]
\hlmax & \hspace*{-0.40cm} Maximum allowable headway of line $l$ during the planning horizon; $l$ $\in$ $L$ \\*[0.05cm]
$n_t(l)$ & \hspace*{-0.40cm} The starting terminal node index of line $l$; $l$ $\in$ $L$ \\*[0.1cm]
$n_{prev} $ & \hspace*{-0.40cm} The transfer node prior to node $n$ in the selected nodes sequence of a line, or terminal  \\
&\hspace*{-0.40cm} if $n$ is the first node of the line; $l$ $\in$ $L$, $n_{prev}$ $\in$ \nl  \\
\ttds & \hspace*{-0.40cm} Bus running time between two consecutive transfer nodes, or terminal to the first transfer \\
&\hspace*{-0.40cm} node, for line $l$ from the previous node $n_{prev}$ to its  consecutive one $n$ in the sequence in  \\
&\hspace*{-0.40cm} scenario $s$; $l$ $\in$ $L$, $n$ $\in$ \nl, $s$ $\in$ $S$\\
\aws & \hspace*{-0.40cm}  Walking time between stop locations of line $l$ and $l'$ at node $n$ in scenario $s$;\\
& \hspace*{-0.40cm} $(n,l,l') \in \tp$, $s$ $\in$ $S$ \\*[0.1cm]
\tds & \hspace*{-0.40cm} Number of transferring passengers in vehicle $p$ of line $l$ who want to transfer to \\*[0.1cm]
& \hspace*{-0.40cm} line $l'$ at node $n$  in scenario $s$;  $(n,l,l') \in \tp$, $p$ $\in$ \pl, $s$ $\in$ $S$ \\
\TldLs & \hspace*{-0.40cm} Total number of local boarding passengers of a low-frequency line, bus $q$ of line $l'$  \\
& \hspace*{-0.40cm} at node $n$ in scenario $s$; $l'$ $\in$ $L_{LF}$, $q$ $\in$ \plpq, $n$ $\in$ \nlq, $s$ $\in$ $S$ \\
\alds & \hspace*{-0.40cm} Number of alighting passengers from bus $q$ of line $l'$ at node $n$ in scenario $s$; \\
& \hspace*{-0.40cm} $l'$ $\in$ $L$, $q$ $\in$ \plpq, $n$ $\in$ \nlq, $s$ $\in$ $S$\\
\ssps & \hspace*{-0.40cm} Net number of people alighting and boarding at the stops between the two consecutive  \\
& \hspace*{-0.40cm} selected transfer nodes $n_{prev}$ and $n$ for bus $q$ of line $l'$ in scenario $s$; \\
& \hspace*{-0.40cm} $l'$ $\in$ $L$, $q$ $\in$ \plpq, $n_{prev},n$ $\in$ \nlq, $s$ $\in$ $S$\\*[0.05cm]
\lams & \hspace*{-0.40cm} Arrival rate (person/time) of local passenger demand for high-frequency line $l'$ at node $n$ \\
& \hspace*{-0.40cm} in scenario $s$; $l'$ $\in$ $L_{Hf}$, $n$ $\in$ \nlq, $s$ $\in$ $S$\\*[0.12cm]
$D$ & \hspace*{-0.40cm} The total demand for a low-frequency line \\
$\nudemand$ & \hspace*{-0.40cm} The percentage of the total demand for a low-frequency line arrives as the first group\\
 & \hspace*{-0.40cm} of local passengers\\
$\betao$ & \hspace*{-0.40cm} The percentage of a low-frequency line $h_{l'}$ when the first group of local passengers\\
& \hspace*{-0.40cm} arrive in scenario $s$ \\
$\betatt$  & \hspace*{-0.40cm} The percentage of a low-frequency line $h_{l'}$ when the second group of local passengers\\
& \hspace*{-0.40cm} arrive in scenario $s$\\
\qS & \hspace*{-0.40cm} {{The first trip of line $l'$ in the next planning horizon; $l' \in L$}}\\
\ctw & \hspace*{-0.40cm} Relative penalty cost of transfer waiting time\\
\cvt & \hspace*{-0.40cm}  Relative penalty cost of unnecessary in-vehicle time\\
\cvh & \hspace*{-0.40cm}  Relative out-of-vehicle penalty cost due to delayed bus\\
$a^t$ & \hspace*{-0.40cm}  Alighting time required per passenger\\
$b^t$ & \hspace*{-0.40cm}  Boarding time required per passenger\\
\prob & \hspace*{-0.40cm} The probability of occurrence of scenario, $s$ $\in$ $S$\\
\bottomrule
\end{tabular}}
\end{table}

\begin{table}[!ht]
\caption{Notation: Decision variables for both high- and low-frequency lines}
\label{t:variablesAll}
\small
\centering
\resizebox{0.89\textwidth}{!}{\begin{tabular}{ll} \toprule
\textbf{For each $ n \in NT_{lt}, l$ $\in$ $PL^{n}$, $p$ $\in$ \pl}: &\\
\midrule
\Pdep & \hspace*{-4.80cm} Timetable departure time of bus $p$ of line $l$ from node $n$\\
\midrule
\midrule
\textbf{For each $ n \in NT_{lt}, l$ $\in$ $PL^{n}$, $p$ $\in$ \pl, $s$ $\in$ $S$}: &\\
\midrule
\Adepq & \hspace*{-4.80cm} Departure time of bus $p$ of line $l$ at node $n$\\
\midrule
\midrule
\textbf{For each $ n \in NT, l$ $\in$ $PL^{n}$, $p$ $\in$ \pl, $s$ $\in$ $S$}: &\\
\midrule
\Aarr & \hspace*{-4.80cm} Arrival time of bus $p$ of line $l$ at node $n$ \\
\ippp & \hspace*{-4.80cm} 1: if departure time of bus $p$ of line $l$ at node $n$ equals its timetable departure time,\\
& \hspace*{-4.80cm} 0: otherwise\\
\Tbdnsp & \hspace*{-4.80cm} Total number of passengers who board on bus $p$ of line $l$ at node $n$ \\
\ivdsp & \hspace*{-4.80cm} Number of in-vehicle passengers of vehicle $p$ of line $l$ when it arrives at node $n$\\
\RDiffp  & \hspace*{-4.80cm} The unnecessary stopping time of bus $p$ of line $l$ at node $n$ \\
\Ewaitp & \hspace*{-4.80cm} Equals the difference between the timetable departure time and the arrival of bus $p$  \\
& \hspace*{-4.80cm} of line $l$ at node $n$ when the bus arrives after its timetable departure time\\
\TEwaitp & \hspace*{-4.80cm} Equals \Ewait if \Ewaitp exceeds the threshold, otherwise, equals zero\\
\TotalBDEwp & \hspace*{-4.80cm} Local demand arriving between the timetable departure time and departure time of   \\
& \hspace*{-4.80cm} bus $p$ of line $l$ at node $n$ with high-frequency when the bus arrives after its timetable \\
& \hspace*{-4.80cm} departure time; or equals to total local demand for bus $p$ of line $l$ at node $n$  \\
& \hspace*{-4.80cm} with low-frequency who are penalized due to the delayed arrival of the bus at Zone 4\\
\GBDEwp & \hspace*{-4.80cm} Successful transfer Type 1 who had to wait after the timetable departure time for the \\
& \hspace*{-4.80cm} arrival of their bus and being penalized for the delayed departure time of\\
& \hspace*{-4.80cm} the bus $p$ of line $l$ at node $n$\\
\pddop & \hspace*{-4.80cm} Number of passengers being penalized due to delayed departure time when bus $p$ of  \\
& \hspace*{-4.80cm} line $l$ at node $n$ arrives before its timetable departure time\\
\vtdop & \hspace*{-4.80cm} The in-vehicle demand being penalized for unnecessary stopping time for\\
& \hspace*{-4.80cm} a bus $p$ of line $l$ at node $n$\\
$\text{Ivdd}_{lp}^{sn}$ & \hspace*{-4.80cm} In-vehicle demand of bus $p$ of line $l$ when it arrives at node $n$ \\
\midrule
\midrule
\textbf{For each $(n,l,l') \in \tp, p \in \pl$, $q$ $\in$ \plpqs, $s$ $\in$ $S$}: &\\
\midrule
\ndwaitss & \hspace*{-4.80cm} Transfer waiting time of passengers from bus $p$ of line $l$ transferring to bus $q$ of line  \\
& \hspace*{-4.80cm} $l'$ at node, i.e., the exact time people wait for their connecting bus to arrive   \\
& \hspace*{-4.80cm} (if passengers arrive during the service time of bus $q$, this value is zero)\\
\twaits & \hspace*{-4.80cm} Transfer waiting time of passengers from bus $p$ of line $l$ transferring to bus $q$ of line  \\
& \hspace*{-4.80cm} $l'$ at node $n$ with the help of dwell transfer buffer time, i.e., successful transfer Type 2 \\
\iwaits & \hspace*{-4.80cm} 1: if the difference of $\ndwait$ and $\twait$ is exactly equal to $\tb$ for   \\
& \hspace*{-4.80cm} a high-frequency connecting line  or $\tbso$ for a low-frequency connecting line, \\
& \hspace*{-4.80cm} 0: otherwise \\
\bottomrule
\end{tabular}}
\end{table}

\begin{table}[!ht]
\caption{Notation: Time-related decision variables for high-frequency connecting lines}
\label{t:variables}
\small
\centering
\resizebox{0.89\textwidth}{!}{\begin{tabular}{ll} \toprule
\textbf{For each $ n \in NT, l'$ $\in$ $L_{Hf}$, $q$ $\in$ \plpq, $s$ $\in$ $S$}: &\\
\midrule
\tbs & \hspace*{-5.80cm} Dwell transfer buffer time of bus $q$ of line $l'$ at node $n$\\
\ttbos & \hspace*{-5.80cm} Transfer buffer time for bus $q$ of line $l'$ at node $n$ (can be negative)\\
\tbbs & \hspace*{-5.80cm} Non-negative transfer buffer time for bus $q$ of line $l'$ at node $n$ (equal \ttbos if it is positive)\\
\ibs & \hspace*{-5.80cm} 0: if \tbb is zero, 1 otherwise; for bus $q$ of line $l'$ at node $n$\\
\servos & \hspace*{-5.80cm} Service time for the first group of local passengers and Type 1 of successfully transferred \\
& \hspace*{-5.80cm}  passengers to bus $q$ of line $l'$ at node $n$\\*[0.1cm]
\servts & \hspace*{-5.80cm} Service time for Type 2 of successfully transferred passengers to bus $q$ of line $l'$ at node $n$\\*[0.1cm]
\servfs & \hspace*{-5.80cm} Service time for Type 3 of successfully transferred passengers to bus $q$ of line $l'$ at node $n$\\*[0.1cm]
\servths & \hspace*{-5.80cm} Service time for the second group of local passengers for bus $q$ of line $l'$ at node $n$,\\
& \hspace*{-5.80cm} when the bus arrives at Zone 1\\*[0.1cm]
\serves & \hspace*{-5.80cm} Service time for the second group of local passengers for bus $q$ of line $l'$ at node $n$,\\
& \hspace*{-5.80cm} when the bus arrives at Zone 2\\*[0.1cm]
$\deltahh$ & \hspace*{-5.80cm} The difference between timetable departure time and arrival time for bus $q$ of line $l'$  \\
& \hspace*{-5.80cm} at node $n$, if the bus arrives at Zone 1. 0: otherwise.\\
\arpi & \hspace*{-5.80cm} 1: if bus  $q$ of line $l'$ at node $n$ arrives after or the same time as its timetable departure\\
& \hspace*{-5.80cm} time, 0: otherwise  \\
\alti & \hspace*{-5.80cm} 1: if alighting time of bus $q$ of line $l'$ at node $n$ is equal or larger than $\deltahh$ bus \\
& \hspace*{-5.80cm} $q$ of line $l'$ at node $n$, 0: otherwise \\
\dwb  & \hspace*{-5.80cm} Equals to the larger value of alighting time and $\deltahh$ for bus $q$ of line $l'$ at node $n$\\
\tbop & \hspace*{-5.80cm} Transfer buffer time upper bound for bus $q$ of line $l'$ at node $n$\\
\TIs & \hspace*{-5.80cm} 1: if service time needed for boarding passengers and \tbbs is less than or equal to  \\
& \hspace*{-5.80cm} the time needed for alighting passengers for bus $q$ of line $l'$ at node $n$, 0: otherwise \\
\dwsi & \hspace*{-5.80cm} Equals to the maximum value of alighting time and boarding passengers time plus \tbbs\\
\midrule
\midrule
\textbf{For each $(n,l,l') \in \tph, p \in \pl$, $q$ $\in$ \plpqs, $s$ $\in$ $S$}: &\\
\midrule
\yys  & \hspace*{-5.80cm} 1: if arrival time of transferring passengers from bus $p$ of line $l$ is less than arrival time of\\
& \hspace*{-5.80cm} bus $q$ of line $l'$ at node $n$ plus \tbs, 0: otherwise\\
\yyyys  & \hspace*{-5.80cm} 1: if there is a successful transfer of Type 3 from feeder bus $p$ of line $l$ to connecting \\
& \hspace*{-5.80cm} bus $q$ of line $l'$ at node $n$, 0: otherwise\\
\zys  & \hspace*{-5.80cm} 1: the first eligible connecting bus $q$ of line $l'$ for transferring passengers of feeder \\
& \hspace*{-5.80cm} vehicle $p$ of line $l$ at node $n$; 0: either not eligible or not the first eligible connecting bus \\
\sys  & \hspace*{-5.80cm} 1: if there is a successful transfer of Type 1 or 2 from feeder bus $p$ of line $l$ to connecting  \\
& \hspace*{-5.80cm} line $l'$;  0: if a successful transfer of Type 3 from feeder bus $p$ of line $l$ to connecting line $l'$\\
\tys  & \hspace*{-5.80cm} 1: if there is a successful transfer, Type 1 or 2 from feeder bus $p$ of line $l$ to connecting\\
& \hspace*{-5.80cm} bus $q$ of line $l'$ at node $n$, 0: otherwise \\
\bottomrule
\end{tabular}}
\end{table}

\begin{table}[!ht]
\caption{Notation: Demand-related decision variables for high-frequency lines}
\label{t:variablesDemand}
\small
\centering
\resizebox{0.89\textwidth}{!}{\begin{tabular}{ll} \toprule
\textbf{ For each $ n \in N, l'$ $\in$ \tpqqfh, $q$ $\in$ \plpq, {{$s$ $\in$ $S$}}}  & \\
\midrule
\bdos & \hspace*{-5.30cm} Total number of the first group of local and transferring passengers who get  \\
& \hspace*{-5.30cm} successful transfer of Type 1 to bus $q$ of line $l'$ at node $n$\\
\bdts & \hspace*{-5.30cm} Total number of transferring passengers who get successful transfer of Type 2 to\\
 & \hspace*{-5.30cm} bus $q$ of line $l'$ at node $n$ \\
\bdths & \hspace*{-5.30cm} Total number of the second group of local passengers of bus $q$ of line $l'$ \\
 & \hspace*{-5.30cm} at node $n$ arriving in Zone 1 \\
\bdfs & \hspace*{-5.30cm} Total number of transferring passengers who get successful transfer of Type 3 to\\
 & \hspace*{-5.30cm} bus $q$ of line $l'$ at node $n$ \\
\bdes & \hspace*{-5.30cm} Total number of the second group of local passengers to bus $q$ of line $l'$ at node $n$\\
 & \hspace*{-5.30cm} arriving in Zone 2 \\
\bottomrule
\end{tabular}}
\end{table}

\begin{table}[!ht]
\caption{Notation: Time-related decision variables for low-frequency lines - I}
\label{t:variablesLO}
\small
\centering
\resizebox{0.89\textwidth}{!}{\begin{tabular}{ll} \toprule
\textbf{For each $ n \in NT, l'$ $\in$ $L_{Lf}$, $q$ $\in$ \plpq, $s$ $\in$ $S$}: &\\
\midrule
\STo & \hspace*{-5.25cm} 0: if bus $q$ of line $l'$ at node $n$ arrives at Zone 1, 1: otherwise\\
\STt & \hspace*{-5.25cm} 0: if bus $q$ of line $l'$ at node $n$ arrives at Zone 1 or 2, 1: otherwise\\
\STot & \hspace*{-5.25cm} 1: if bus $q$ of line $l'$ at node $n$ arrives at Zone 2, 0: otherwise\\
\STtt & \hspace*{-5.25cm} 1: if bus $q$ of line $l'$ at node $n$ arrives at Zone 3, 0: otherwise\\
\STth & \hspace*{-5.25cm} 1: if bus $q$ of line $l'$ at node $n$ arrives at Zone 4, 0: otherwise\\
\tbso & \hspace*{-5.25cm} Dwell transfer buffer time of bus $q$ of line $l'$ at node $n$  for the first \\
& \hspace*{-5.25cm} possible successful transfer Type 2\\
\ttboso & \hspace*{-5.25cm} First transfer buffer time for bus $q$ of line $l'$ at node $n$ (can be negative)\\
\ttbosop & \hspace*{-5.25cm} Non-negative first transfer buffer time for bus $q$ of line $l'$ at node $n$ \\
& \hspace*{-5.25cm} (equals \ttboso if it is positive)\\
\ttbost & \hspace*{-5.25cm} Second transfer buffer time for bus $q$ of line $l'$ at node $n$ arrives at Zone 1 (can be negative)\\
\ttbosth & \hspace*{-5.25cm} Third transfer buffer time for bus $q$ of line $l'$ at node $n$ in scenario $s$ arrive at Zone 1; \\
& \hspace*{-5.25cm} or second transfer buffer time for bus $q$ of line $l'$ at node $n$ arrive at Zone 2 (can be negative)\\
\servos & \hspace*{-5.25cm} Service time for Type 1 of successfully transferred passenger to bus $q$ of line $l'$ at node $n$ if \\*[0.1cm]
& \hspace*{-5.25cm}  the bus arrives at Zone 1; or service time for the first group of local passengers and Type 1  \\
& \hspace*{-5.25cm} of successfully transferred passengers to bus $q$ of line $l'$ at node $n$ if the bus arrives at Zone 2;   \\
& \hspace*{-5.25cm} or service time for all local passengers and Type 1 of successfully transferred passenger to \\
& \hspace*{-5.25cm} bus $q$ of line $l'$ at node $n$  if the bus arrives at Zone 3 or four\\*[0.1cm]
\servts & \hspace*{-5.25cm} Service time for the first group of successfully transferred passengers with Type 2 to bus $q$ of\\
& \hspace*{-5.25cm} line $l'$ at node $n$\\
\servLo & \hspace*{-5.25cm} Service time for the first group of local passengers arriving at Zone 2 for bus $q$ of line $l'$\\
& \hspace*{-5.25cm} at node $n$ \\
\servLt & \hspace*{-5.25cm} Service time for the second group of local passengers arriving at Zone 3 for bus $q$ of\\
& \hspace*{-5.25cm} line $l'$ at node $n$\\
\servtsq & \hspace*{-5.25cm} Service time for the second group of successfully transferred passengers Type 2 in Zone 2  \\
& \hspace*{-5.25cm} for bus $q$ of line $l'$ at node $n$ arrives at Zone 1\\
\servthsq & \hspace*{-5.25cm} Service time for the second group of successfully transferred passengers Type 3 in Zone 2  \\
& \hspace*{-5.25cm} for bus $q$ of line $l'$ at node $n$  arrives at Zone 1\\
\servtsg & \hspace*{-5.25cm} Service time for the third group of successfully transferred passengers Type 2 in Zone 3 for\\
& \hspace*{-5.25cm} bus $q$ of line $l'$ at node $n$ arriving at Zone 1; or service time for the second group of \\
& \hspace*{-5.25cm} transferred passengers Type 2 in Zone 3 for bus $q$ of line $l'$ at node $n$ arriving Zone 2\\
\servthsg & \hspace*{-5.25cm} Service time for the third group of successfully transferred passengers Type 3 in Zone 3 for \\
& \hspace*{-5.25cm} bus $q$ of line $l'$ at node $n$  arrives at Zone 1; or service time for second group of successfully  \\
& \hspace*{-5.25cm} transferred passengers Type 3 in Zone 3 for bus $q$ of line $l'$ at node $n$ arrives at Zone 2\\
\servELf & \hspace*{-5.25cm} Excess dwell time after the timetable departure time if bus $q$ of line $l'$ at node $n$ does not\\
& \hspace*{-5.25cm} arrive at Zone 4; or equals the total dwell time if bus $q$ of line $l'$ at node $n$ arrives at Zone 4\\
\alightb & \hspace*{-5.25cm} Alighting budget time for bus $q$ of line $l'$ at node $n$\\
\alightE & \hspace*{-5.25cm} Excess alighting time if the alighting budget is not enough for bus $q$ of line $l'$ at node $n$\\
\alightbi & \hspace*{-5.25cm} 1: if the alighting budget of bus $q$ of line $l'$ at node $n$ is enough, 0: otherwise\\
\dwIE & \hspace*{-5.25cm} Equals the larger value of \alightE \hspace*{0.03cm} and \servELf for bus $q$ of line $l'$ at node $n$\\
\TIs & \hspace*{-5.25cm} 0: if \alightE \hspace*{0.03cm}is less than \servELf \hspace*{0.03cm}for bus $q$ of line $l'$ at node $n$, 0: otherwise\\
\bottomrule
\end{tabular}}
\end{table}

\begin{table}[!ht]
\caption{Notation: Time-related decision variables for low-frequency lines - II}
\label{t:variablesLOT}
\small
\centering
\resizebox{0.859\textwidth}{!}{\begin{tabular}{ll} \toprule
\textbf{For each $ n \in NT, l'$ $\in$ $L_{Lf}$, $q$ $\in$ \plpq, $s$ $\in$ $S$}: &\\
\midrule
\dwbozo & \hspace*{-4.975cm} First stopping budget to determine upper limit of first transfer buffer time for bus $q$  \\
& \hspace*{-4.975cm} of line $l'$ at node $n$ arrives at Zone 1\\
\dwbiozo & \hspace*{-4.975cm} 1: if the \dwbozo is less than or equal to \servos plus \servts for bus $q$ of line $l'$ at \\
& \hspace*{-4.975cm} node $n$ arrives at Zone 1; 0: otherwise\\
\dwbozt & \hspace*{-4.975cm} First stopping budget to determine upper limit of first transfer buffer time for bus $q$  \\
& \hspace*{-4.975cm} of line $l'$ at node $n$ arrives at Zone 2\\
\dwbiozt & \hspace*{-4.975cm} 1: if the \dwbozt is less than or equal to \servos plus \servts for bus $q$ \\
& \hspace*{-4.975cm} of line $l'$ at node $n$ arrives at Zone 2; 0: otherwise\\
\dwbozth & \hspace*{-4.975cm} First stopping budget to determine upper limit of first transfer buffer time for bus  \\
& \hspace*{-4.975cm} $q$ of line $l'$ at node $n$ arrives at Zone 3\\
\dwbiozth & \hspace*{-4.975cm} 1: if the \dwbozth is less than or equal to \servos plus \servts for bus $q$ of line $l'$ at \\
& \hspace*{-4.975cm} node $n$arrives at Zone 3; 0: otherwise\\
\dwbozto & \hspace*{-4.975cm} Second stopping budget to determine upper limit of second transfer buffer time for  \\
& \hspace*{-4.975cm} bus $q$ of line $l'$  at node $n$ arrives at Zone 1\\
\dwboztho & \hspace*{-4.975cm} Third stopping budget to determine upper limit of second transfer buffer time for \\
& \hspace*{-4.975cm} bus $q$ of line $l'$  at node $n$ arrives at Zone 1\\
\dwboztht & \hspace*{-4.975cm} Second stopping budget to determine upper limit of second transfer buffer time for  \\
& \hspace*{-4.975cm} bus $q$ of line $l'$ at node $n$ arrives at Zone 2\\
\dwbozthf & \hspace*{-4.975cm} Stopping budget to determine upper limit of second transfer buffer time at Zone 3 \\
& \hspace*{-4.975cm} interval for bus $q$ of line $l'$ at node $n$ arrives at Zone 2 or 1\\
\dwbioztho & \hspace*{-4.975cm} 1: if the \dwbozthf is less than or equal to \Est plus \servLt plus \servtsg for bus  \\
& \hspace*{-4.975cm} $q$ of line $l'$ at node $n$ arrives at Zone 3; 0: otherwise\\
\tbopozo & \hspace*{-4.975cm} Upper limit of first transfer buffer time for bus $q$ of line $l'$ at node $n$ arrives at Zone 1\\
\tbopozt & \hspace*{-4.975cm} Upper limit of first transfer buffer time for bus $q$ of line $l'$ at node $n$ arrives at Zone 2\\
\tbopozth & \hspace*{-4.975cm} Upper limit of first transfer buffer time for bus $q$ of line $l'$ at node $n$ arrives at Zone 3\\
\tbopozt & \hspace*{-4.975cm} Upper limit of second transfer buffer time for bus $q$ of line $l'$ at node $n$ arrives at Zone 1\\
\tbopoztot & \hspace*{-4.975cm} Upper limit of \ttbost for bus $q$ of line $l'$ at node $n$ arrives at Zone 1\\
\tbopoztoth & \hspace*{-4.975cm} Upper limit of \ttbosth for bus $q$ of line $l'$ at node $n$ arrives at Zone 1\\
\Eso & \hspace*{-4.975cm} Excess time passing the Zone 2 boundary for bus $q$ of line $l'$ at node $n$ arrives at Zone 1\\
\Estzo & \hspace*{-4.975cm} Excess time passing the Zone 3 boundary for bus $q$ of line $l'$ at node $n$ arrives at Zone 1\\
\Estzt & \hspace*{-4.975cm} Excess time passing the Zone 3 boundary for bus $q$ of line $l'$ at node $n$ arrives at Zone 2\\
\Est & \hspace*{-4.975cm} Excess time passing the Zone 3 boundary for bus $q$ of line $l'$ at node $n$\\
\Esthzot & \hspace*{-4.975cm} Excess time passing the Zone 4 boundary for bus $q$ of line $l'$ at node $n$ arrives at Zone 1/2 \\
\Esthzoth & \hspace*{-4.975cm} Excess time passing the Zone 4 boundary for bus $q$ of line $l'$ at node $n$ arrives at Zone 3\\
\Esth & \hspace*{-4.975cm} Excess time passing the Zone 4 boundary for bus $q$ of line $l'$ at node $n$\\
\Esoi & \hspace*{-4.975cm} 1: if the first stopping budget is enough and not pass the Zone 2 boundary for bus \\
& \hspace*{-4.975cm} $q$ of line $l'$ at node $n$ arrives at Zone 1, 0: otherwise\\
\Estzoi & \hspace*{-4.975cm} 1: if the second stopping budget is enough and not pass the Zone 3 boundary for bus \\
& \hspace*{-4.975cm} $q$ of line $l'$ at node arrives at Zone 1, 0: otherwise\\
\Estzti & \hspace*{-4.975cm} 1: if the second stopping budget is enough and not pass the Zone 3 boundary for bus \\
& \hspace*{-4.975cm} $q$ of line $l'$ at node $n$ arrives at Zone 2, 0: otherwise\\
\Esthzoti & \hspace*{-4.975cm} 1: if the second stopping budget is enough and not pass the Zone 4 boundary for bus \\
& \hspace*{-4.975cm} $q$ of line $l'$ at node  $n$ arrives at Zone 1 or 2, 0: otherwise\\
\Esthzothi & \hspace*{-4.975cm} 1: if the second stopping budget is enough and not pass the Zone 4 boundary for bus \\
& \hspace*{-4.975cm} $q$ of line $l'$ at node $n$ arrives at Zone 3, 0: otherwise\\
\bottomrule
\end{tabular}

}
\end{table}

\begin{table}[!ht]
\caption{Notation: Time-related decision variables for low-frequency lines - III}
\label{t:variablesLT}
\small
\centering
\resizebox{0.89\textwidth}{!}{\begin{tabular}{ll} 
\toprule
\textbf{For each $(n,l,l') \in \tph, p \in \pl$, $q$ $\in$ \plpqs, $s$ $\in$ $S$}: &\\
\midrule
\yyso  & \hspace*{-6.180cm} 1: if arrival time of transferring passengers from bus $p$ of line $l$ is less than or equal arrival \\
& \hspace*{-6.180cm} time of bus $q$ of line $l'$ at node $n$ plus \tbs, 0: otherwise\\
\yyyyso  & \hspace*{-6.180cm} 1: if there is a successful transfer of Type 3 (the first one) from feeder bus $p$ of line\\
& \hspace*{-6.180cm} $l$ to connecting bus $q$ of line $l'$ at node $n$ \\
\yyst & \hspace*{-6.180cm} 1: if successful transfer Type 2 occurs at Zone 2 for the second time for a connecting bus \\
& \hspace*{-6.180cm} $q$ of line $l'$ at node $n$ arrives at Zone 1\\
\yqo  & \hspace*{-6.180cm} 1: if arrival time of transferring passengers from bus $p$ of line $l$ is in Zone 2,  0: otherwise\\
\yqt  & \hspace*{-6.180cm} 1: if arrival time of transferring passengers from bus $p$ of line $l$ is less than Zone 2 boundary \\
& \hspace*{-6.180cm} plus the service time of the first group of local passengers and excess time from Zone 1 for  \\
& \hspace*{-6.180cm} the connecting bus $q$ of line $l'$ at node $n$, 0: otherwise\\
\yqot & \hspace*{-6.180cm} 1: if successful transfer Type 2 occurs at Zone 2 for the second time for the bus arrives \\
& \hspace*{-6.180cm} at Zone 1, 0: otherwise\\
\yyyystc & \hspace*{-6.180cm} 1: if successful transfer Type 3 can occur at Zone 2 for a connecting bus $q$ of line $l'$  \\
& \hspace*{-6.180cm} at node $n$, 0: otherwise.\\
\yyyyst & \hspace*{-6.180cm} 1: if successful transfer Type 3 can occur at Zone 2 for the second time for a connecting \\
& \hspace*{-6.180cm} bus $q$ of line $l'$ at node $n$ arrives at Zone 1 0: otherwise.\\
\ygo  & \hspace*{-6.180cm} 1: if arrival time of transferring passengers from bus $p$ of line $l$ is in Zone 3, 0: otherwise\\
\ygt  & \hspace*{-6.180cm} 1: if arrival time of transferring passengers from bus $p$ of line $l$ is less than Zone 3 boundary \\
& \hspace*{-6.180cm} plus the service time of all local passengers and excess time from Zone 2 for the  \\
& \hspace*{-6.180cm} connecting bus bus $q$ of line $l'$ at node $n$, 0: otherwise\\
\ygot & \hspace*{-6.180cm} 1: if successful transfer Type 2 can occur at Zone 3 for the third time for the bus arrives \\
& \hspace*{-6.180cm} at Zone 1, or for the second time if the bus arrives at Zone 2, 0: otherwise\\
\yysthzo  & \hspace*{-6.180cm} 1: if successful transfer Type 2 occurs at Zone 3 for the third time for a connecting bus $q$ of \\
& \hspace*{-6.180cm} line $l'$ at node $n$ arrives at Zone 1 \\
\yysthzt  & \hspace*{-6.180cm} 1: if successful transfer Type 2 occurs at Zone 3 for the third time for a connecting bus $q$ of \\
& \hspace*{-6.180cm} line $l'$ at node $n$ arrives at Zone 1 \\ 
\yyyysthc & \hspace*{-6.180cm} 1: if successful transfer Type 3 can occur at Zone 3 for for a connecting bus $q$ of line $l'$\\
& \hspace*{-6.180cm} at node $n$ arrives 0: otherwise.\\
\yyyysthzo & \hspace*{-6.180cm} 1: if successful transfer Type 3 can occur at Zone 3 for the third time for a connecting \\
& \hspace*{-6.180cm} bus $q$ of line $l'$ at node $n$ arrives at Zone 1 \\
\yyyysthzt & \hspace*{-6.180cm} 1: if successful transfer Type 3 can occur at Zone 3 for the second time for a connecting \\
& \hspace*{-6.180cm} bus $q$ of line $l'$ at node $n$ arrives at Zone 2, 0: otherwise\\
\zyso  & \hspace*{-6.180cm} 1: the first eligible connecting bus $q$ of line $l'$ for transferring passengers of feeder vehicle $p$ of\\
& \hspace*{-6.180cm} line $l$ at node $n$; 0: either not eligible or not the first eligible connecting bus\\
\syso  & \hspace*{-6.180cm} 1: if there is a successful transfer of Type 1 or 2 for the first time from feeder bus $p$ of line $l$  \\
& \hspace*{-6.180cm} to connecting line $l'$, 0: if there is a successful transfer of Type 3 from feeder bus $p$ of line $l$ \\
 & \hspace*{-6.180cm}to connecting line  $l'$\\
\tyso  & \hspace*{-6.180cm} 1: if there is a successful transfer, Type 1 or 2 for the first time from feeder bus $p$ of line $l$ to  \\
& \hspace*{-6.180cm} connecting bus $q$ of line $l'$ at node $n$, 0: otherwise \\
\bottomrule
\end{tabular}}
\end{table}

\begin{table}[!ht]
\caption{Notation: Demand-related decision variables for low-frequency lines}
\label{t:variablesDemandL}
\small
\centering
\resizebox{0.89\textwidth}{!}{\begin{tabular}{ll} \toprule
\textbf{ For each $ n \in N, l'$ $\in$ \tpqqfl, $q$ $\in$ \plpq, {{$s$ $\in$ $S$}}}  & \\
\midrule
\midrule
\bdos  & \hspace*{-5.30cm} Number of Type 1 of successfully transferred passengers to bus $q$ of line $l'$ at node $n$  \\
& \hspace*{-5.30cm} if the bus arrives at Zone 1; Or number of first group of local passengers and Type 1 \\
& \hspace*{-5.30cm} of successfully transferred passenger to bus $q$ of line $l'$ at node $n$ if the bus arrives at  \\
& \hspace*{-5.30cm} Zone 2; or number of all local passengers and Type 1 of successfully transferred  \\
& \hspace*{-5.30cm} passenger to bus $q$ of line $l'$ at node $n$ if the bus arrives at Zone 3 or 4\\*[0.1cm]
\bdts & \hspace*{-5.30cm} Number of first group of successfully transferred passengers with Type 2 to bus $q$ of\\
& \hspace*{-5.30cm} line $l'$ at node $n$\\
\bdq & \hspace*{-5.30cm} Number of second group of successfully transferred passengers Type 2 in Zone 2  \\
& \hspace*{-5.30cm} for connecting bus $q$ of line $l'$ at node $n$ arrives at Zone 1\\
\bdthsq & \hspace*{-5.30cm} Number of second group of successfully transferred passengers Type 3 in Zone 2  \\
& \hspace*{-5.30cm} for connecting bus $q$ of line $l'$ at node $n$ arrives at Zone 1\\
\bdg & \hspace*{-5.30cm} Number of third group of successfully transferred passengers Type 2 in Zone 3 for \\
& \hspace*{-5.30cm} connecting bus $q$ of line $l'$ at node $n$ arrives at Zone 1. or number of the second group  \\
& \hspace*{-5.30cm} of successfully transferred passengers passengers Type 2 in Zone 3 for connecting\\
& \hspace*{-5.30cm} bus $q$ of line $l'$ at node $n$ arrives at Zone 2\\
\bdthsg & \hspace*{-5.30cm} Number of third group of successfully transferred passengers Type 3 in Zone 3 for  \\
& \hspace*{-5.30cm} connecting bus $q$ of line $l'$ at node $n$ arrives at Zone 1\\
& \hspace*{-5.30cm}Number of second group of successfully transferred passengers Type 3 in Zone 3 for \\
& \hspace*{-5.30cm} connecting bus $q$ of line $l'$ at node $n$ arrives at Zone 2\\
\bottomrule
\end{tabular}}
\end{table}

\begin{table}[!ht]
\caption{Notation: Waiting-related decision variables for both high- and low-frequency lines}
\label{t:variablesDelayy}
\small
\centering
\resizebox{0.89\textwidth}{!}{\begin{tabular}{ll} \toprule
\textbf{For each $ n \in NT, l$ $\in$ $L$, $q$ $\in$ \plpq, $s$ $\in$ $S$}: &\\
\midrule
\Rths & \hspace*{-5.30cm} The threshold for unnecessary stopping time\\
\Rthsi & \hspace*{-5.30cm} 1: if the unnecessary stopping time exceeds the threshold for a bus $q$ of line $l'$ \\
& \hspace*{-5.30cm} at node $n$, 0: otherwise\\
\Ewait & \hspace*{-5.30cm} The difference between the timetable departure time and the arrival of bus $q$ \\
& \hspace*{-5.30cm} of line $l'$ at node $n$ when the bus arrives after its timetable departure time \\
\EwaitII & \hspace*{-5.30cm} 1: if \Ewait is equals or more than the threshold, 0: otherwise\\
\Aths & \hspace*{-5.30cm} The threshold for delayed departure time\\\
\hspace*{-0.18cm} \Athsi & \hspace*{-5.30cm} 1: if the difference between the delayed departure time and the timetable departure  \\
& \hspace*{-5.30cm} time of bus $q$ of  line $l'$ at node $n$ exceeds the threshold, 0: otherwise\\
\TYE & \hspace*{-5.30cm} 1: successful transfer Type 1 when the connecting bus $q$ of line $l'$ at node $n$ arrives  \\
& \hspace*{-5.30cm} after its timetable departure time, 0:otherwise.\\
\STTYE & \hspace*{-5.30cm} 1: if successful transfer Type 1 waits more than the threshold after the timetable departure \\
& \hspace*{-5.30cm}  time when bus $q$ of line $l'$ at node $n$ arrives after its timetable departure time, 0:otherwise\\
\midrule
\midrule
\textbf{For each $ n \in NT, l$ $\in$ $L_{Lf}$, $q$ $\in$ \plpq, $s$ $\in$ $S$}: &\\
\midrule
\dwIzo  & \hspace*{-5.10cm} Total boarding time required for a  bus $q$ of line $l'$ at node $n$ arrives at Zone 1\\
\dwIzt  & \hspace*{-5.10cm} Total boarding time required for a bus $q$ of line $l'$ at node $n$ arrives at Zone 2\\
\dwIzth  & \hspace*{-5.10cm} Total boarding time required for a bus $q$ of line $l'$ at node $n$ arrives at Zone 3\\
\dwsi  & \hspace*{-5.10cm} Total boarding time of a bus $q$ of line $l'$ at node $n$ regardless of arrival Zone\\
\Esthii & \hspace*{-5.10cm} 1: if a low-frequency bus $q$ of line $l'$ at node $n$ arrives before its timetable departure time \\
& \hspace*{-5.10cm} yet leaves after its timetable departure time, 0: otherwise\\
\midrule
\midrule
\textbf{For each $ n \in NT, l$ $\in$ $L_{Hf}$, $q$ $\in$ \plpq, $s$ $\in$ $S$}: &\\
\midrule
\pipa & \hspace*{-5.30cm} 1: if a high-frequency bus $q$ of line $l'$ at node $n$ arrives before its timetable departure \\
& \hspace*{-5.30cm} time yet leaves after its timetable departure time, 0: otherwise\\
\bottomrule
\end{tabular}

}
\end{table}

\clearpage

\section{Model Details for Low-frequency Connecting Lines}
\label{appsec:lowfreq}
As explained in Section \ref{subsubsec:zonedefns}, we establish four zones for the arrival of a low-frequency bus. Regarding the local passengers' arrival patterns used for two zones, based on the historical data, we assume that we can convert the distribution of local passengers' arrival at a bus stop into two parts. The first group of passengers would arrive $\betao \hlq$ before the timetable departure time, while the second group would arrive $\betat \hlq$ before the timetable departure time, where $\betat \hlq$ is smaller than $\betao \hlq$. These parameters can be selected independently of headways and based on available data for each line/trip. 
\eqref{Eq83}-\eqref{Eq90} are specified to distinguish the four zones. If ST1 equals zero in 
\eqref{Eq83}, bus $q$ of line $l$ would arrive in Zone 1 when there are no local passengers. Otherwise, if ST1 equals one, it simply indicates that the bus is in Zone 2, 3, or 4. Likewise, if ST2 equals one, the bus would arrive in Zone 3 or 4. Thus, if ST1 equals one and ST2 equals zero, ST12 equals one indicating that the bus would arrive in Zone 2. If ST3 equals zero in 
\eqref{Eq87} and ST2 equals one, then ST23 equals one in 
\eqref{Eq90}, so the bus would arrive in Zone 3. Finally, if ST3 equals one, the bus would arrive in Zone 4 after its timetable departure time.

\medskip
For all $(n,l') \in \tpqql, q \in \ql$:
\vspace*{-0.25cm}
\begin{alignat}{2}
    & \hspace*{-0.7cm} (\Pdep - \Aarrq) \leq (\betao \hlq)  \Leftrightarrow \STo = 1 
    \label{Eq83}\\
    & \hspace*{-0.7cm} (\Pdep - \Aarrq) \leq \betat \hlq \Leftrightarrow \STt = 1 
    \label{Eq85}\\
    & \hspace*{-0.7cm} \Pdep \leq \Aarrq \Leftrightarrow \STth = 1 
    \label{Eq87}\\
    & \hspace*{-0.7cm} \STot = \STo - \STt 
    \label{Eq89}\\
    & \hspace*{-0.7cm} \STtt = \STt - \STth 
    \label{Eq90}
\end{alignat}



Now that we know the arrival zone of a connecting bus of a low-frequency line, we need to determine the number of passengers while accounting for their relative arrival time in comparison to the arrival time of  bus $q$ on line $l'$ and its timetable departure time. Next, we investigate the three different waiting times and the bus departure time.
{It is worth noting that the approach to determining service times and identifying possible successful transfers for low-frequency lines is fairly similar to the approach for high-frequency lines which we described in Section \ref{subsec:highfreq}. The only differences are how we handle local passengers' dwell time, dwell transfer buffer upper limit and the transfer buffer upper limit, which vary depending on which zone a bus arrives in.} Regardless, if any group of transferring passengers arrives before the connecting bus, iwait = 1 as indicated for high-frequency lines, they should be taken into account in GBD1.
Although successful transfer Types 1 and 2, in 
\eqref{Eq91}, can occur regardless of arrival zone of a bus, the upper limit for dwell transfer buffer, dtb, would vary depending on the arrival zone. According to 
\eqref{Eq93}, dtb should be less than serv1. Nonetheless, depending on the bus arrival zone, serv1 would cover service time for different groups of passengers, GBD1. 
\eqref{Eq94} is specified so that the right value of generated demand for serv1 is determined based on the values of zone-based binary variables. If, ST1 = 1, no local passengers are included. If ST12 = 1, a bus would arrive in Zone 2, and the first group of local passengers are included in serv1, 
in \eqref{Eq95}. Similarly, if a bus arrives in Zone 3 or Zone 4, all local passengers would have arrived prior to the bus arrival and are therefore included in GBD1. The number of GBD2, successfully transferring passengers via Type 2, and their service time, serv2, are determined by 
\eqref{Eq96} and \eqref{Eq97}. Note that determining iwait and TY for successful transfer Types 1 and 2 will be explained later in this section, though the approach is almost identical to that of high-frequency lines.
\begin{alignat}{2}
    & \hspace*{-0.9cm} \text{For all } (n, l,l')\in \tpl, q \in \ql, p \in \pl: \nonumber \\
    & \hspace*{-0.2cm}(\Aarr + \aws) \leq (\Aarrq + \tbso) \Leftrightarrow \yyso = 1 
    \label{Eq91} \\[0.15cm]
    & \hspace*{-0.9cm} \text{For all } (n,l') \in \tpqql, q \in \ql: \nonumber \\
   &\hspace*{-0.2cm} \tbso \leq \servos 
   \label{Eq93} \\
    & \hspace*{-0.2cm} \bdos = \sum_{l\in L: (n,l,l') \in \tpl}\sum_{p \in \pl}{\iwaits \tds} + (\STot \gamma D) + (\STtt D) + (\STth D)  
  \label{Eq94}\\
   & \hspace*{-0.2cm} \servos = b^t  \bdos  
   \label{Eq95}\\
    & \hspace*{-0.2cm} \bdts = \sum_{l\in L: (n,l,l') \in \tpl}\sum_{p \in \pl}{(\tys - \iwaits) \tds}  
    \label{Eq96}\\
   & \hspace*{-0.2cm} \servts = b^t  \bdts  
   \label{Eq97}
\end{alignat}


\eqref{Eq98} identifies possible successful transfers of Type 3. This type of successful transfer can happen regardless of the zone a bus arrives in. However, the bus arrival zone determines the upper limit for transfer buffer time, tb1. To consider the relative arrival time marker of the first group of local passengers while deciding regarding successful Type 3, we assume that the arrival of the first group of possible successful transfer Type 3 for a connecting bus arriving in Zone 1 would be before the arrival of the first group of local passengers, i.e., Zone 1. As a result, the first dwell budget for a bus arriving at Zone 1, dwb1z1, is computed using 
\eqref{Eq100}-\eqref{Eq104}. In 
\eqref{Eq106}-\eqref{Eq108}, the upper limit of transfer buffer time for that bus, tbO1z1, is derived by comparing dwb1z1 to serv1 and serv2. 

On the other hand, if a bus arrives in Zone 2, the first dwell budget, dwb1z2, is defined by 
\eqref{Eq110}-\eqref{Eq114}, while the upper limit of transfer buffer time, tbO1z2, is determined by 
\eqref{Eq116}-\eqref{Eq118}. Likely, 
\eqref{Eq120}-\eqref{Eq124} define a bus first dwell budget arriving at Zone 3, dwb1z3, while its upper transfer limit, tbO1z3, is derived via 
\eqref{Eq126}-\eqref{Eq128}. It is worth noting that the upper limit of transfer buffer for a bus arriving in Zone 4 is zero. Thus, 
\eqref{Eq130} determine the upper limit for tb1 in 
\eqref{Eq98}. Constraints \eqref{Eq131} show the minimum negative value for tb1, the same concept described above for high-frequency lines. Accordingly, GBD3 and serv3 are determined through 
\eqref{Eq132} and \eqref{Eq133}, respectively.
\begin{alignat}{2}
    & \hspace*{-1.05cm} \text{For all } (n, l,l')\in \tpl, q \in \ql, p \in \pl: \nonumber \\
    & \hspace*{-0.5cm}(\Aarr + \aws) = (\Aarrq + \servos + \servts +\ttboso ) \Leftrightarrow \yyyyso = 1 
  \label{Eq98} \\[0.15cm]
      & \hspace*{-1.05cm} \text{For all } (n, l')\in \tpqql, q \in \ql: \nonumber \\
  & \hspace*{-0.5cm} \STo = 0 \Rightarrow \dwbozo = (\Pdepq - \betao \hlq) - \Aarrq  
   \label{Eq100} \\
  & \hspace*{-0.5cm} \STo = 1 \Rightarrow \dwbozo = 0  
  \label{Eq102} \\
  & \hspace*{-0.5cm} \dwbozo \leq (\servos + \servts) \Leftrightarrow \dwbiozo = 1 
  \label{Eq104} \\
   & \hspace*{-0.5cm} \dwbiozo = 0 \Rightarrow \tbopozo = \dwbozo - (\servos + \servts)   
   \label{Eq106} \\
   & \hspace*{-0.5cm} \dwbiozo = 1 \Rightarrow \tbopozo = 0   
   \label{Eq108} \\
  & \hspace*{-0.5cm} \STot = 1 \Rightarrow \dwbozt = (\Pdep - \betat \hlq) - \Aarrq    
   \label{Eq110} \\
  & \hspace*{-0.5cm} \STot = 0 \Rightarrow \dwbozt = 0     
  \label{Eq112} \\
  & \hspace*{-0.5cm} \dwbozt \leq (\servos + \servts) \Leftrightarrow \dwbiozt = 1 
  \label{Eq114} \\
   & \hspace*{-0.5cm} \dwbiozt = 0 \Rightarrow \tbopozt = \dwbozt - (\servos + \servts)   
  \label{Eq116} \\
   & \hspace*{-0.5cm} \dwbiozt = 1 \Rightarrow \tbopozt = 0  
   \label{Eq118} \\
  & \hspace*{-0.5cm}  \STtt = 1 \Rightarrow \dwbozth = (\Pdep - \Aarrq)   
  \label{Eq120} \\
  & \hspace*{-0.5cm} \STtt = 0 \Rightarrow \dwbozth = 0    
  \label{Eq122} \\
    & \hspace*{-0.5cm} \dwbozth \leq (\servos + \servts) \Leftrightarrow \dwbiozth = 1 
   \label{Eq124} \\*[0.15cm]
   & \hspace*{-0.5cm} \dwbiozth = 0 \Rightarrow \tbopozth = \dwbozth - (\servos + \servts)   
   \label{Eq126} \\
   & \hspace*{-0.5cm} \dwbiozth = 1 \Rightarrow \tbopozth = 0   
   \label{Eq128} \\
    & \hspace*{-0.5cm} \ttboso \leq \tbopozo + \tbopozt + \tbopozth  
    \label{Eq130} \\
   & \hspace*{-0.5cm} -\servts \leq \ttboso 
   \label{Eq131}\\
    & \hspace*{-0.5cm} 
    \ttbosop = \max \{0,\ttboso\} 
    \label{Eq129}\\
   & \hspace*{-0.5cm} \bdths = \sum_{l\in L: (n,l,l') \in \tpl}\sum_{p \in \pl}{\yyyyso \tds}  
   \label{Eq132}\\
   & \hspace*{-0.5cm} \servths = b^t  \bdths 
   \label{Eq133}
\end{alignat}
As previously stated, one of the challenges in determining dwell time for low-frequency lines is accounting for the arrival patterns of local passengers, specifically when a bus arrives at Zone 1 or 2. To address this challenge, we specify the upper limit of tb1 for successful transfer Type 3 as dwell time budget which is the difference between the arrival of the next group of local passengers and the arrival of the bus, dwb1z1 and dwb1z2. First, let us assume a bus arrives at Zone 1. So far, we have identified the required service times for successful transfer Types 1, 2, and 3.
We calculate a time marker by adding all of these service times. If it is more than the first dwell budget, dwb1z1, we pass the point when the first group of local passengers arrives. Thus, 
\eqref{Eq134}-\eqref{Eq138} determine this excess time, Es1. As a result, the time marker is now in Zone 2. 

\medskip
For all $(n, l')\in \tpl, q \in \ql$:
\vspace*{-0.25cm}
\begin{alignat}{2}
  & \hspace*{-0.7cm} (\servos + \servts + \tbbso + \servths) \leq \dwbozo \Leftrightarrow \Esoi = 1  
  \label{Eq134} \\
  & \hspace*{-0.7cm} \Esoi = 0 \Rightarrow \Eso = (\servos + \servts + \tbbso + \servths) - \dwbozo  
  \label{Eq136} \\
  & \hspace*{-0.7cm} \Esoi = 1 \Rightarrow 
  \Eso = 0  
  \label{Eq138}
\end{alignat}
At this time marker, we have the first group of local passengers. Then we verify if a group of transferring passengers arrives between the start of Zone 2 and the end of excess service time from Zone 1, Es1, plus the service time of the first group of local passengers, servL1. If yes, successful transfer Type 2 can occur again. 
\eqref{Eq140}-\eqref{Eq149} are defined for this aim.
Note that to formulate this semi-type two of a successful transfer, we have to formulate it differently. If Q12 equals one, the feeder would arrive within the above-described time interval, according to 
\eqref{Eq144}. However, if the bus does arrive in Zone 2, the successful transfer is, in fact, the original Type 2 successful transfer that we have previously established. As a consequence, to avoid duplication, the final successful transfer decision variable is YY2. YY2 is determined not just by the value of Q12, but also by the zone in which the bus originally arrived. As a result of 
\eqref{Eq145}-\eqref{Eq147}, if ST1 equals zero and Q12 equals one, YY2 equals one. Otherwise, YY2 should be zero. The associated demand for YY2 and its service time are determined through 
\eqref{Eq148} and \eqref{Eq149}, respectively. Note that the dwell time of the first group of passengers, servL1, is determined through 
\eqref{Eq150}.
\begin{alignat}{2}
& \hspace*{-1.5cm} \text{For all } (n, l,l')\in \tpl, q \in \ql, p \in \pl: \nonumber \\  
   & \hspace*{-0.7cm}(\Aarr + \aws) \geq (\Pdep - \betao \hlq) \Leftrightarrow \yqo = 0  
  \label{Eq140} \\
  & \hspace*{-0.7cm}(\Aarr + \aws) \leq (\Pdep - \betao \hlq + \servLo + \Eso) \Leftrightarrow \yqt = 1 
   \label{Eq142} \\
   & \hspace*{-0.7cm} \yqot = \yqt - \yqo 
   \label{Eq144} \\
   & \hspace*{-0.7cm} \yyst \leq \yqot + \STo 
   \label{Eq145} \\
   & \hspace*{-0.7cm} \yyst \leq 1 - \STo 
   \label{Eq146} \\
   & \hspace*{-0.7cm} \yyst \geq \yqot - \STo 
   \label{Eq147} \\[0.15cm]
    & \hspace*{-1.5cm} \text{For all } (n, l')\in \tpqql, q \in \ql: \nonumber \\
   & \hspace*{-0.7cm} \bdq = \sum_{l\in L: (n,l,l') \in \tp}\sum_{p \in \pl}{\yyst \tds}  
   \label{Eq148}\\
    & \hspace*{-0.7cm} \servtsq = b^t  \bdq 
    \label{Eq149}\\
    & \hspace*{-0.7cm} \servLo = b^t \nudemand D 
    \label{Eq150}
\end{alignat}
Similarly, semi-type 3 successful transfers are possible as well, provided by 
\eqref{Eq151}-\eqref{Eq155}. As previously stated, we can have the second transfer buffer time by taking into account the Es1, servL1, and serv2q while defining the second dwell budget for a bus arriving at Zone 1, dwb2z1. Note that dwb2z1 equals to the duration of Zone 2, described by 
\eqref{Eq156}. The upper limit of the transfer buffer would be determined using the same rationale described above in 
\eqref{Eq100}-\eqref{Eq130}. Constraints \eqref{Eq167} also ensure the negative limit for tb2. The number of successful transfers and their required service time are determined through 
\eqref{Eq168} and \eqref{Eq169}, respectively.
\begin{alignat}{2}
& \hspace*{-0.5cm} \text{For all } (n, l,l')\in \tpl, q \in \ql, p \in \pl: \nonumber \\ 
  & \hspace*{-0.2cm}\Aarr + \aws = \Pdep - \betao \hlq + \servLo + \servtsq + \Eso + \ttbost 
  \Leftrightarrow \yyyystc = 1 
   \label{Eq151} \\
  & \hspace*{-0.2cm} \yyyyst \leq \yyyystc + \STo 
  \label{Eq153} \\
  & \hspace*{-0.2cm} \yyyyst \geq \yyyystc - \STo 
  \label{Eq154} \\
  & \hspace*{-0.2cm} \yyyyst \leq 1 - \STo 
  \label{Eq155} \\[0.15cm]
& \hspace*{-0.5cm} \text{For all } (n,l') \in \tpqqfl, q \in \ql: \nonumber \\ 
  & \hspace*{-0.2cm} \STo = 0 \Rightarrow \dwbozto = (\betao \hlq - \betat \hlq)   
  \label{Eq156} \\
  & \hspace*{-0.2cm} \STo = 1 \Rightarrow \dwbozto = 0    
  \label{Eq158} \\
  & \hspace*{-0.2cm} \dwbozto \leq (\Eso + \servLo + \servtsq) \Leftrightarrow \dwbiozto = 1 
   \label{Eq160} \\
   & \hspace*{-0.2cm} \dwbiozto = 0 \Rightarrow \tbopoztot = \dwbozto - (\Eso + \servLo + \servtsq)   
   \label{Eq162} \\
   & \hspace*{-0.2cm} \dwbiozto = 1 \Rightarrow \tbopoztot = 0   
   \label{Eq164} \\
   & \hspace*{-0.2cm}  \ttbost \leq \tbopoztot 
   \label{Eq166} \\
    & \hspace*{-0.2cm} 
    \ttbostp = \max\{0,\ttbost\} 
    \label{Eq165}\\
   & \hspace*{-0.2cm}  -\servtsq  \leq \ttbost 
   \label{Eq167} \\
  & \hspace*{-0.2cm} \bdthsq = \sum_{l\in L: (n,l,l') \in \tpl}\sum_{p \in \pl}{\yyyyst \tds} 
  \label{Eq168}\\
  & \hspace*{-0.2cm} \servthsq = b^t  \bdthsq  
  \label{Eq169}
\end{alignat}
So far, the time marker tracking for a bus arriving in Zone 1 is reached to the arrival of the second group of transferring passengers. Similar to Es1, we must assess whether or not, subject to the service time of passengers arriving at the bus during Zone 2, the time marker for the bus crosses the boundary of Zone 3. 
\eqref{Eq170}-\eqref{Eq174} determine Es2z1. However, we should also apply the same assessment to a bus that arrives in Zone 2. For this bus, after the first group of local passengers and any potential transferring passengers, we determine whether or not its time marker passes the time when the second group of passengers will arrive, Es2z2, through 
\eqref{Eq176}-\eqref{Eq180}. Finally, 
\eqref{Eq182}-\eqref{Eq186} determine the final value of Es2 based on arrival zone of the bus. 

\medskip
For all $(n,l') \in \tpqql, q \in \ql$:
\vspace*{-0.25cm}
\begin{alignat}{2}
    & \hspace*{-0.7cm} (\servLo + \servtsq + \Eso + \tbbst + \servthsq) \leq \dwbozto \Leftrightarrow  \Estzoi = 1 
    \label{Eq170} \\
   & \hspace*{-0.7cm} \Estzoi = 0 \Rightarrow \Estzo = (\servos + \servts + \tbbst + \servths) - \dwbozto  
   \label{Eq172} \\
   & \hspace*{-0.7cm} \Estzoi = 1 \Rightarrow \Estzo = 0   
   \label{Eq174} \\
   & \hspace*{-0.7cm} (\servos + \servts + \tbbso + \servths) \leq \dwbozt \Leftrightarrow \Estzti = 1 
   \label{Eq176} \\
   & \hspace*{-0.7cm} \Estzti = 0 \Rightarrow \Estzt = (\servos + \servts + \tbbso + \servths) - \dwbozt   
   \label{Eq178} \\
   & \hspace*{-0.7cm} \Estzti = 1 \Rightarrow \Estzt = 0   
   \label{Eq180} \\
   & \hspace*{-0.7cm} 
   \STo = 0 \Rightarrow \Est = \Estzo
   \label{Eq182} \\
   & \hspace*{-0.7cm} 
   \STot = 1 \Rightarrow \Est = \Estzt
   \label{Eq184} \\
   & \hspace*{-0.7cm} 1 - \STo = \STot = 0 \Rightarrow \Est = 0 
   \label{Eq186}
\end{alignat}
As previously stated, if a bus arrives in Zone 3, the related generated demand and successful transfers are determined by 
\eqref{Eq91}-\eqref{Eq133}. On the other hand, if a bus arrives in Zone 1 or 2 there are two types of successful transfers that can occur in Zone 3. The concept is the same as 
\eqref{Eq151}-\eqref{Eq169}, with the exception that the bus arrival zone might now be Zone 1 or Zone 2. Considering the service time for the second group of local passengers, servL2, via 
\eqref{Eq187}, constraints \eqref{Eq188}-\eqref{Eq192} identify semi-successful transfer Type 2. The final value for that successful transfer, Y3z1 and Y3z2, would then be decided by 
\eqref{Eq193}-\eqref{Eq199} depending on the bus arrival Zone avoiding duplication.
The corresponding demand and the required service time are also determined through 
\eqref{Eq200} and \eqref{Eq201}. 
\begin{alignat}{2}
& \hspace*{-1.2cm} \text{For all } (n,l') \in \tpqqfl, q \in \ql: \nonumber \\ 
   & \hspace*{-0.7cm} \servLt = b^t (1 - \nudemand) D 
   \label{Eq187}\\[0.15cm]
& \hspace*{-1.2cm} \text{For all } (n, l,l')\in \tpl, q \in \ql, p \in \pl: \nonumber \\ 
   & \hspace*{-0.7cm}(\Aarr + \aws) \geq (\Pdep - \betat \hlq) \Leftrightarrow \ygo = 0
   \label{Eq188} \\
  & \hspace*{-0.7cm}(\Aarr + \aws) \leq (\Pdep - \betat \hlq + \servLt + \Est) \Leftrightarrow \ygt = 1 
   \label{Eq190} \\
    & \hspace*{-0.7cm} \ygot = \ygt - \ygo 
    \label{Eq192} \\
   & \hspace*{-0.7cm} \yysthzo \leq \ygot + \STo 
   \label{Eq193} \\
   & \hspace*{-0.7cm} \yysthzo \leq 1 - \STo 
   \label{Eq194} \\
   & \hspace*{-0.7cm} \yysthzo \geq \ygot - \STo 
   \label{Eq195} \\
   & \hspace*{-0.7cm} \yysthzt \leq \STot 
   \label{Eq196} \\
   & \hspace*{-0.7cm} \yysthzt \leq \ygot + \STot 
   \label{Eq197} \\
   & \hspace*{-0.7cm} \yysthzt \geq \ygot + ( \STot - 1) 
   \label{Eq198} \\
   & \hspace*{-0.7cm} \yysthzt \leq \ygot + (1 - \STot) 
   \label{Eq199} \\[0.15cm]
   & \hspace*{-1.2cm} \text{For all } (n,l') \in \tpqqfl, q \in \ql: \nonumber \\ 
   & \hspace*{-0.7cm} \bdg = \sum_{l\in L: (n,l,l') \in \tpl}\sum_{p \in \pl}{(\yysthzo + \yysthzt) \tds}  
   \label{Eq200}\\
   & \hspace*{-0.7cm} \servtsg = b^t  \bdg 
   \label{Eq201}
\end{alignat}
\eqref{Eq202} is also used to evaluate semi-successful Type 3. As a result, the final values of YY3z1 and YY3z2 are through 
\eqref{Eq204}-\eqref{Eq210}. The upper limit for the important element of semi-successful Type 3, tb3, should be calculated carefully given that now it depends on whether the bus arrival is Zone 1 or 2. The dwell budget are determined based on the bus arrival zone through 
\eqref{Eq213}-\eqref{Eq222}. Then considering Es2, servL2, and serv2g, the upper limit for tb3 is determined via 
\eqref{Eq224}-\eqref{Eq228}. Constraints \eqref{Eq229} also ensure the lower limit for tb3.
\begin{alignat}{2}
& \hspace*{-0.4cm} \text{For all } (n, l,l')\in \tpl, q \in \ql, p \in \pl: \nonumber \\ 
   & \hspace*{0cm} \Aarr \hspace*{-0.1cm} + \hspace*{-0.1cm} \aws = \Pdep \hspace*{-0.1cm} - \hspace*{-0.1cm} \betat \hlq \hspace*{-0.1cm} + \hspace*{-0.1cm} \servLt \hspace*{-0.1cm} + \hspace*{-0.1cm} \servtsg \hspace*{-0.1cm} + \hspace*{-0.1cm} \Est \hspace*{-0.1cm} + \hspace*{-0.1cm} \ttbosth 
    \Leftrightarrow \yyyysthc = 1 
    \label{Eq202} \\
  & \hspace*{-0cm} \yyyysthzo \leq \yyyysthc + \STo 
  \label{Eq204} \\
  & \hspace*{-0cm} \yyyysthzo \geq \yyyysthc - \STo 
  \label{Eq205} \\
  & \hspace*{-0cm} \yyyysthzo \leq 1 - \STo 
  \label{Eq206} \\
  & \hspace*{-0cm} \yyyysthzt \leq \STot 
  \label{Eq207} \\
  & \hspace*{-0cm} \yyyysthzt \leq \yyyysthc + \STot 
  \label{Eq208} \\
  & \hspace*{-0cm} \yyyysthzt \geq \yyyysthc + (\STot - 1) 
  \label{Eq209} \\
  & \hspace*{-0cm} \yyyysthzt \leq \yyyysthc + (1 - \STot) 
  \label{Eq210} \\[0.15cm]
  & \hspace*{-0.4cm} \text{For all } (n,l') \in \tpqqfl, q \in \ql: \nonumber \\ 
   & \hspace*{-0cm} \bdthsg = \sum_{l\in L: (n,l,l') \in \tpl}\sum_{p \in \pl}{(\yyyysthzo + \yyyysthzt) \tds} 
    \label{Eq211}\\
     & \hspace*{-0cm} \servthsg = b^t  \bdthsg  
     \label{Eq212} \\
     & \hspace*{-0cm} \STo = 0 \Rightarrow \dwboztho = (\betat \hlq)   
     \label{Eq213} \\
    & \hspace*{-0cm} \STo = 1 \Rightarrow \dwboztho = 0   
    \label{Eq215} \\
    & \hspace*{-0cm} \STot = 1 \Rightarrow \dwboztht = (\betat \hlq)  
    \label{Eq217} \\
    & \hspace*{-0cm} \STot = 0 \Rightarrow \dwboztht = 0    
    \label{Eq219} \\
   & \hspace*{-0cm}  \dwbozthf = \dwboztht + \dwboztho  
   \label{Eq221} \\
   & \hspace*{-0cm} \dwbozthf \leq (\Est + \servLt + \servtsg) \Leftrightarrow \dwbioztho = 1 
    \label{Eq222} \\
   & \hspace*{-0cm} \dwbioztho = 0 \Rightarrow \tbopoztoth = \dwbozthf - (\Est + \servLt + \servtsg)  
   \label{Eq224} \\
   & \hspace*{-0cm} \dwbioztho = 1 \Rightarrow \tbopoztoth = 0   
   \label{Eq226} \\
   & \hspace*{-0cm}  \ttbosth \leq \tbopoztoth 
   \label{Eq228} \\
   & \hspace*{-0cm}  -\servtsg \leq \ttbosth 
   \label{Eq229}\\
  & \hspace*{-0cm} 
  \ttbosthp = \max\{0, \ttbosth \} 
  \label{Eq227}
\end{alignat}
After identifying all possible successful transfers based on a bus arrival zone, we must ensure that the first connecting bus is assigned to a feeder. To do so, variables ZY, SY, and TY are defined using the same logic as for high-frequency lines, 
\eqref{Eq230} - \eqref{Eq234}. The distinction here is that all potential successful transfers are considered while determining SY variables. Additionally, 
\eqref{Eq235} ensures that only one connecting bus is assigned to a feeder bus. 

\medskip
\hspace*{-0.88cm} For all $(n, l,l')\in \tpl, p \in \pl$:
\vspace*{-0.25cm}
\begin{alignat}{2}
  &\hspace*{-1.17cm} \zyso = \yyso  
  && q = 1
  \label{Eq230} \\
  &\hspace*{-1.17cm} \zyso = \yyso - \yynso  
  && q \in \ql, q\neq 1
  \label{Eq231} \\
  &\hspace*{-1.17cm} \syso = \sum_{q \in \ql}{(\zyso)}-\sum_{q \in \plpn}{(\yyyyso + \yysthzo + \yysthzt+ } \quad 
  &&\ \nonumber\\
  & \hspace*{2cm}
  { \yyyysthzo + \yyst + \yyyyst + \yyyysthzt)} 
  \label{Eq232}\\
  &\hspace*{-1.17cm} \tyso \leq \syso  
  && q \in \ql
  \label{Eq233} \\
  &\hspace*{-1.17cm} \tyso \geq (\syso - 1) + \zyso 
  && q \in \ql
  \label{Eq234} \\
  & \hspace*{-1.17cm}\sum_{q \in \ql}{(\tyso + \yyyyso + \yysthzo + \yysthzt  +  }  &&\ \nonumber\\
  & \hspace*{-0.15cm}  
  { \yyyysthzo + \yyyysthzt + \yyst + \yyyyst) + \ttys} = 1 
  \label{Eq235}
\end{alignat}
To establish if a bus can depart based on its timetable departure time, we calculate whether the total of dwell times considering the arrival zone of a bus would surpass the timetable departure time. The excess dwell time of a bus with an arrival Zone of 1 or 2, ES3z12, that would pass the timetable departure time is calculated by 
\eqref{Eq236}-\eqref{Eq240}. However, the extra time for a bus with an arrival Zone of 3, Es3z3, must be computed individually using 
\eqref{Eq242}-\eqref{Eq246}. Constraints \eqref{Eq248}-\eqref{Eq254}, which account for the bus arrival zone, are then used to compute the final value for Es3.
\begin{alignat}{2}
  & \hspace*{-0.22cm} \text{For all } (n,l') \in \tpqqfl, q \in \ql: \nonumber \\ 
   & \hspace*{0cm} \servLt + \servtsg + \Est + \tbbsth + \servthsg \leq \dwbozthf   
   \Leftrightarrow \Esthzoti = 1  
   \label{Eq236} \\
   & \hspace*{-0cm} \Esthzoti = 0 \Rightarrow \Esthzot = \servLt + \servtsg + \Est + \tbbsth + \servthsg 
   - \dwbozthf  
   \label{Eq238} \\
   & \hspace*{-0cm} \Esthzoti = 1 \Rightarrow \Esthzot = 0   
   \label{Eq240} \\
   & \hspace*{-0cm} \servos + \servts + \tbbso + \servths \leq \dwbozth \Leftrightarrow \Esthzothi = 1 
   \label{Eq242} \\
   & \hspace*{-0cm} \Esthzothi = 0 \Rightarrow \Esthzoth = \servos + \servts + \tbbso + \servths - \dwbozth  
   \label{Eq244} \\
   & \hspace*{-0cm} \Esthzothi = 1 \Rightarrow \Esthzoth = 0   
   \label{Eq246} \\
    & \hspace*{-0cm} 
    \STo = 0 \Rightarrow \Esth = \Esthzot
    \label{Eq248} \\
   & \hspace*{-0cm} 
   \STot = 1 \Rightarrow \Esth = \Esthzot
   \label{Eq250} \\
   & \hspace*{-0cm} 
   \STtt = 1 \Rightarrow \Esth = \Esthzoth 
   \label{Eq252} \\
   & \hspace*{-0cm} 
   (1 - \STo) = (\STot) = (\STtt) = 0 \Rightarrow \Esth = 0
   \label{Eq254}
\end{alignat}


On the other hand, if a bus arrives at Zone 4, i.e., after its timetable departure time, its required dwell time, the sum of serv1, serv2, and serv3 would certainly pass the published departure time. Thus, the extra dwell time required for a bus that occurs after the associated timetable departure time, servELf, is determined by 
\eqref{Eq255}-\eqref{Eq257}. Note that, if the bus arrives at Zone 4, i.e., ST3 = 1, serv3 is for the passengers who have successful transfer Type 3 but with transfer buffer time of zero. In order to determine the actual departure time of a bus, similar to high-frequency lines, we also need to consider the required alighting time for a bus. For a low-frequency line, to account for all possible arrival zones, we define alighting time budget via 
\eqref{Eq258}-\eqref{Eq259}. Then we compare if the required dwell time surpass the alighting time budget and if yes how longer it is in 
\eqref{Eq260}-\eqref{Eq262}. Finally, by comparing the excess time required for boarding and alighting, through 
\eqref{Eq263}-\eqref{Eq265} , and the departure time of the bus is determined through 
\eqref{Eq266}-\eqref{Eq267}.
\begin{alignat}{2}
  & \hspace*{-2.5cm} \text{For all } (n,l') \in \tpqqfl, q \in \ql: \nonumber \\ 
    & \hspace*{-1.7cm} \STth = 1 \Rightarrow \servELf = (\servos + \servts + \servths)    
   \label{Eq255} \\
   & \hspace*{-1.7cm} \STth = 0 \Rightarrow \servELf = \Esth   
   \label{Eq257} \\
    & \hspace*{-1.7cm} \STth = 0 \Rightarrow \alightb = \Adep - \Aarrq   
    \label{Eq258} \\
    & \hspace*{-1.7cm} \STth = 1 \Rightarrow \alightb = 0  
    \label{Eq259} \\
    & \hspace*{-1.7cm} \alightb \geq (\aldnns) a^t \Leftrightarrow \alightbi = 1 
    \label{Eq260} \\
    & \hspace*{-1.7cm} \alightbi = 0 \Rightarrow \alightE = (\aldnns) a^t - \alightb   
    \label{Eq261} \\
    & \hspace*{-1.7cm} \alightbi = 1 \Rightarrow \alightE = 0   
    \label{Eq262} \\
   &\hspace*{-1.7cm} \alightE \geq \servELf \Leftrightarrow \TIse = 1 
   \label{Eq263}\\
   & \hspace*{-1.7cm} \TIse = 1 \Rightarrow \dwIE = \alightE    
   \label{Eq264}\\
    & \hspace*{-1.7cm} \TIse = 0 \Rightarrow \dwIE = \servELf  
    \label{Eq265}\\
   & \hspace*{-1.7cm} \STth = 0 \Rightarrow \Adep = \Pdep + \dwIE   
   \label{Eq266}\\
    & \hspace*{-1.7cm} \STth = 1 \Rightarrow \Adep = \Aarrq + \dwIE   
    \label{Eq267}
\end{alignat}

\clearpage

\section{Details on Model Parameters}
\begin{figure}[!ht]
\begin{center}
\centerline{\includegraphics[width=0.75\linewidth]{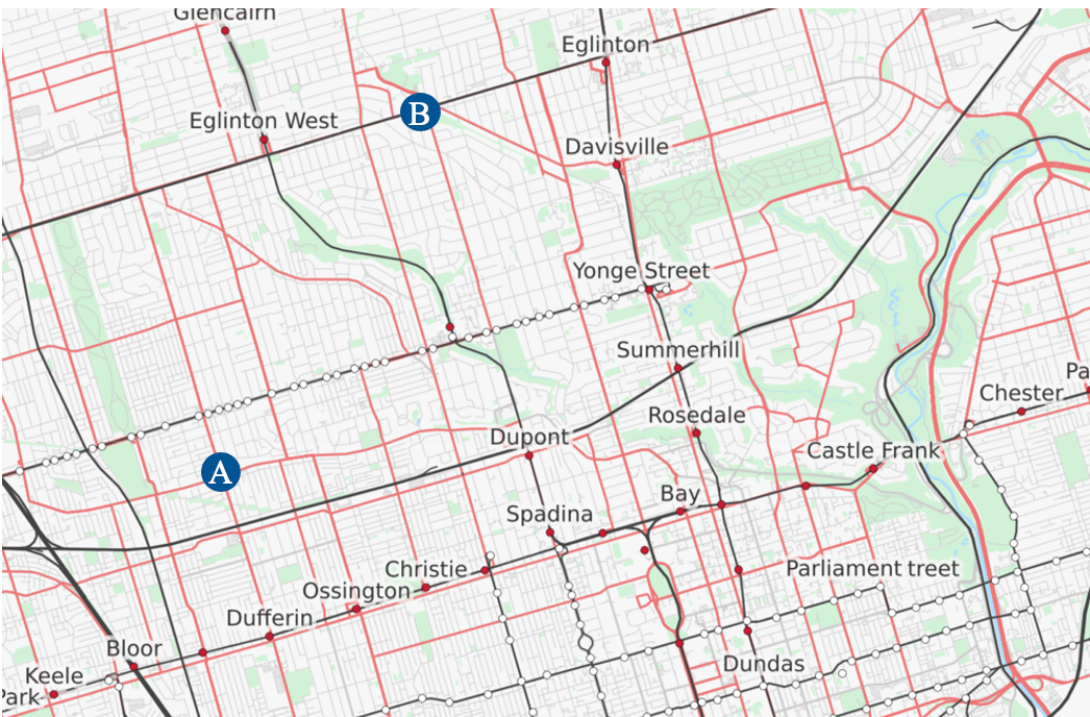}}
\caption{Two transfer nodes in Toronto}
\label{fig:Toronto}
\end{center}
\end{figure}

Figure~\ref{fig:Toronto} provides exact coordinates of the two transfer nodes located in the City of Toronto.

\end{document}